\documentclass[oneside,11pt]{article}
\usepackage[utf8]{inputenc}
\usepackage[T1]{fontenc}

\usepackage{amsthm}
\usepackage{amsmath}
\usepackage{amsfonts}
\usepackage{amssymb}

\usepackage{graphicx} 
\usepackage{graphbox}
\usepackage{float}
\usepackage{booktabs}
\usepackage{multirow}
\usepackage{rotating}
\usepackage{array}
\usepackage{xcolor}
\usepackage{subcaption}
\usepackage[ruled]{algorithm2e}
\usepackage{tabularx}
\usepackage{caption}

\usepackage{comment}

\newcolumntype{C}[1]{>{\centering\arraybackslash}p{#1}}
\newcolumntype{R}[1]{>{\raggedleft\let\newline\\\arraybackslash\hspace{0pt}}m{#1}}
\newcolumntype{L}[1]{>{\raggedright\let\newline\\\arraybackslash\hspace{0pt}}m{#1}}

\usepackage{breakcites}

\usepackage{nowidow}

\usepackage{enumitem}

\usepackage[top=1.5cm,bottom=2cm,right=2.5cm,left=2.5cm]{geometry}
\usepackage{titling}

\usepackage{natbib}

\usepackage{ifplatform}

\usepackage{textcomp}
\usepackage{bbm}

\ifwindows
\fi

\allowdisplaybreaks

\newcommand{\lp}{\left(}
\newcommand{\rp}{\right)}
\newcommand{\lc}{\left[}
\newcommand{\rc}{\right]}

\newcommand{\R}{\mathbb{R}}

\newcommand{\rd}{\mathrm{d}}

\newcommand{\bx}{\boldsymbol{x}}
\newcommand{\be}{\boldsymbol{e}}

\newcommand{\by}{\boldsymbol{y}}
\newcommand{\bX}{\boldsymbol{X}}
\newcommand{\bY}{\boldsymbol{Y}}
\newcommand{\bZ}{\boldsymbol{Z}}

\newcommand{\bmu}{\boldsymbol\mu}

\newcommand{\bv}{\boldsymbol{v}}

\newcommand{\bz}{\boldsymbol{z}}

\newcommand{\bp}{\boldsymbol{p}}

\newcommand{\zero}{\mathbf{0}}

\newcommand{\bxi}{\boldsymbol\xi}

\newcommand{\ba}{\boldsymbol\alpha}

\newcommand{\bzeta}{\boldsymbol\zeta}

\newcommand{\btheta}{\boldsymbol\theta}

\newcommand{\bSigma}{\boldsymbol\Sigma}

\newcommand{\bHcal}{\boldsymbol{\mathcal{H}}}

\newcommand{\bB}{\boldsymbol{B}}

\newcommand{\bA}{\boldsymbol{A}}
\newcommand{\bC}{\boldsymbol{C}}

\newcommand{\bI}{\boldsymbol{I}}

\DeclareMathAlphabet\mathbfcal{OMS}{cmsy}{b}{n}

\newcommand{\etab}{\boldsymbol\eta}

\newcommand{\xib}{\boldsymbol\xi}

\newcommand{\sigmad}{\sigma_{d}}

\newcommand{\Sd}{\mathbb{S}^{d}}

\newcommand{\lrp}[1]{\left(#1\right)}
\newcommand{\lrc}[1]{\left[#1\right]}
\newcommand{\lrb}[1]{\left\{#1\right\}}
\newcommand{\lrpbig}[1]{\big(#1\big)}

\newcommand{\norm}[1]{\left|\left| #1\right|\right|}

\newcommand{\tr}[1]{\mathrm{tr}\left[#1\right]}
\newcommand{\vect}{\mathrm{vec}\,}

\newcommand{\om}[1]{\omega_{#1}}

\DeclareFontFamily{OT1}{pzc}{}
\DeclareFontShape{OT1}{pzc}{m}{it}{<-> s * [1.10] pzcmi7t}{}
\DeclareMathAlphabet{\mathpzc}{OT1}{pzc}{m}{it}

\newtheorem{theorem}{Theorem}[section]
\newtheorem{corollary}{Corollary}[section]

\newtheorem{proposition}{Proposition}[section]
\newtheorem{lemma}{Lemma}[section]

\newcommand{\pct}{\%}

\newcommand{\defin}{:=}

\newcommand{\bw}{\boldsymbol{w}}

\newcommand{\bpsi}{\boldsymbol\psi}
\newcommand{\bXi}{\boldsymbol\Xi}

\newcommand{\vartb}{\boldsymbol\vartheta}
\newcommand{\Sdo}{\mathbb{S}^{1}}

\newcommand{\Sdm}{\mathbb{S}^{d-1}}
\newcommand{\Sdmm}{\mathbb{S}^{d-2}}
\newcommand{\D}{\mathsf D}
\newcommand{\Do}[1]{\mathsf D^{\otimes #1}}
\newcommand{\E}[1]{\mathbb E\lc #1\rc}
\newcommand{\indef}{=:}

\usepackage[pdftex,final,bookmarksnumbered,bookmarksopen=false,breaklinks,colorlinks]{hyperref}
\hypersetup{
	final,
	bookmarksnumbered=true,
	bookmarksopen=true,
	bookmarksopenlevel=0,
	unicode=false,
	pdftoolbar=true,
	pdfmenubar=true,
	pdffitwindow=false,
	pdftitle={Blessing of dimensionality in cross-validated bandwidth selection on the sphere},
	pdfdisplaydoctitle=true,
	pdftoolbar=true,
	pdfmenubar=true,
	pdflang={English},
	pdfauthor={Jose E. Chacon, Eduardo Garcia-Portugues, Andrea Meilan-Vila},
	pdfsubject={arXiv paper},
	pdfcreator={José E. Chacón},
	pdfproducer={José E. Chacón},
	pdfkeywords={Directional data}{High-dimensional data}{Nonparametric statistics}{Smoothing},
	pdfnewwindow=true,
	breaklinks=true,
	hidelinks,
	linkcolor=black,
	citecolor=black,
	filecolor=magenta,
	urlcolor=cyan
}

\newif\ifmain
\maintrue
\newif\ifsupplement
\supplementtrue

\newif\iffigstabs
\figstabstrue

\begin{document}

\ifmain

\title{Blessing of dimensionality in cross-validated bandwidth selection on the sphere}
\setlength{\droptitle}{-1cm}
\predate{}%
\postdate{}%
\date{}

\author{José E. Chacón$^{1,3}$, Eduardo Garc\'ia-Portugu\'es$^{2}$, and Andrea Meil\'an-Vila$^{2}$}
\footnotetext[1]{Department of Mathematics, Universidad de Extremadura (Spain).}
\footnotetext[2]{Department of Statistics, Universidad Carlos III de Madrid (Spain).}
\footnotetext[3]{Corresponding author. e-mail: \href{mailto:jechacon@unex.es}{jechacon@unex.es}.}
\maketitle

\begin{abstract}
 We study the asymptotic behavior of least-squares cross-validation bandwidth selection in kernel density estimation on the $d$-dimensional hypersphere, $d\geq 1$. We show that the exact rate of convergence with respect to the optimal bandwidth minimizing the mean integrated squared error, shown to exist under mild non-uniformity conditions, is $n^{-d/(2d+8)}$, thus approaching the $n^{-1/2}$ parametric rate as $d$ grows. This ``blessing of dimensionality'' in bandwidth selection offers theoretical support for utilizing the conceptually simpler cross-validation selector over plug-in techniques for larger dimensions $d$. We compare this result for bandwidth estimation on the $d$-dimensional Euclidean space through explicit expressions for the asymptotic variance functionals.  Numerical experiments corroborate the speed of this convergence in an array of scenarios and dimensions, precisely illustrating the tipping dimension where cross-validation outperforms plug-in approaches.
\end{abstract}
\begin{flushleft}
\end{flushleft}

\section{Introduction}
\label{sec:intro}

%
Selecting an appropriate smoothing parameter is a central challenge in nonparametric density estimation, as it directly governs the trade-off between bias and variance in the resulting estimator. Among data-driven approaches, least-squares cross-validation (CV), introduced by \citet{Rudemo1982} and \citet{Bowman1984}, has received significant attention due to its appealing asymptotic properties: the sequence of bandwidths it produces yields asymptotically optimal density estimates under rather mild assumptions \citep{hall1983large, stone1984asymptotically}.

In the Euclidean setting, the relative rate of convergence of the multivariate CV bandwidth selector, with respect to the integrated squared error (ISE)-optimal bandwidth was first noted by \citet{Marron1986}, anticipating the results of \citet[][Section~2.1]{HM87}. Subsequent research focused on relative rates of convergence with respect to the mean integrated squared error (MISE)-optimal bandwidth. In the univariate case, these rates were derived by \citet{Park1990} as a direct consequence of the analysis in \cite{HM87}, showing that the rates relative to the MISE-optimal and ISE-optimal bandwidths coincide. This result was later extended to the multivariate case by \cite{Jones1992}, who also highlighted the remarkable fact that the convergence rate improves with increasing dimension, approaching a square-root rate (the fastest possible, as shown in \citealp{Hall1991}) as the dimension approaches infinity. Indeed, \cite{Savchuk2010} noted that cross-validation appears to perform particularly well in scenarios where the density estimation problem becomes more challenging. The behavior differs when the relative rate of convergence is measured with respect to the asymptotic MISE (AMISE)-optimal bandwidth. In the univariate case, \cite{scott1987biased} showed this rate is the same as that for the MISE-optimal bandwidth. In the multivariate setting, \cite{Sain1994} reported that the relative rate with respect to the AMISE-optimal bandwidth also coincides with the MISE-based rate. However, \cite{Duong2005} later corrected this result, showing that for dimensions $d>4$,  the rate slows down due to the discrepancy between the MISE-optimal and the AMISE-optimal bandwidths.

These developments motivate the study of cross-validation in non-Euclidean spaces, particularly on the hypersphere, $\Sd=\{\bx\in\mathbb{R}^{d+1}:\norm{\bx}^2=\bx^\top\bx=1\}$, $d\geq 1$, where directional data naturally arise. Special cases include circular data ($d=1$), such as animal orientations or wave directions \citep{Wang2014}. Spherical data ($d=2$) frequently appear in medicine or astronomy \citep{Marinucci2008}, while higher-dimensional directional data are encountered in areas such as text mining \citep{Banerjee2005} and genetics \citep{Eisen1998}. Comprehensive overviews of directional statistics are provided by \cite{Mardia1999a} and \cite{Ley2017a}, whereas \cite{Pewsey2021} survey recent advances. Several bandwidth selection methods have been proposed for kernel density estimation on directional data. For circular data, \cite{Taylor2008} introduced a plug-in bandwidth selector, and \cite{Oliveira2012} developed a mixture-based alternative. \cite{DiMarzio2011} proposed a CV approach for data on the torus, i.e., the Cartesian product of circles. On the hypersphere $\Sd$, \cite{Garcia-Portugues2013a} proposed a rule-of-thumb bandwidth selector under the von Mises--Fisher model, while \cite{Hall1987} derived two CV bandwidth selectors. More recently, \cite{Tsuruta2020} studied the properties of the plug-in \citep{DiMarzio2011} and the CV \citep{Hall1987} bandwidth selectors for the kernel density estimator on $\mathbb{S}^1$, establishing their asymptotic normality and demonstrating that the convergence rates are $n^{-5/14}$ for the plug-in selector and $n^{-1/10}$ for the CV selector. However, in higher dimensions, the theoretical properties of the CV selector remain largely unexplored.

This work aims to develop the asymptotic theory of the CV bandwidth selector for the kernel density estimator on the hypersphere $\Sd$. Our main contributions are as follows: (\textit{i}) we prove the existence of the MISE-optimal bandwidth under mild non-uniformity conditions; (\textit{ii}) we show the consistency of the CV selector and derive its exact relative rate of convergence with respect to the MISE-optimal bandwidth, showing it is $n^{-d/(2d+8)}$, which approaches the $n^{-1/2}$ parametric rate as $d$ increases, and generalizing the $n^{-1/10}$ rate obtained by \cite{Tsuruta2020} for $\mathbb{S}^1$; (\textit{iii}) we provide the explicit expression for the asymptotic variance, allowing comparison with the Euclidean case and illustrating the ``blessing of dimensionality'' in the spherical setting; and (\textit{iv}) we perform numerical experiments validating the theoretical results and identifying the regimes where CV outperforms plug-in selectors.

The rest of this paper is organized as follows. Section \ref{sec:pre} introduces preliminaries on kernel density estimation on the sphere and its MISE, and establishes the existence of the MISE-optimal bandwidth. Section \ref{sec:cross} proves the consistency of the cross-validation bandwidth selector for spherical data, provides its exact convergence rate with respect to the MISE, and offers theoretical comparisons with the Euclidean case. Section \ref{sec:num} presents numerical experiments that validate the theoretical findings and compare the convergence rates of the cross-validation selector with those from plug-in methods. The paper concludes with a discussion in Section \ref{sec:disc}. Proofs are relegated to the Supplementary Material (SM).

\section{Preliminaries}
\label{sec:pre}

\subsection{Kernel density estimation}
\label{sec:kde}

%
Let $f$ be a probability density function (pdf) on $\Sd=\{\bx\in\mathbb{R}^{d+1}:\norm{\bx}^2=\bx^\top\bx=1\}$, $d\geq 1$, with respect to the surface area measure $\sigmad$ on $\Sd$. We denote by $\om{d}\defin\sigmad(\Sd)=2\pi^{(d+1)/2}/\Gamma((d+1)/2)$ the surface area of $\Sd$.

Let $\bX_1,\ldots,\bX_n$ be an independent and identically distributed (iid) sample from $f$. Let $\bx\in\Sd$ and set $h\in\R_+$. The kernel density estimator (kde) of $f$ at $\bx$ is defined~as
\begin{align}
    \hat{f}(\bx;h)&\defin\frac{1}{n}\sum_{i=1}^nL_{h}(\bx,\bX_i),\quad L_{h}(\bx,\by)\defin c_{d,L}(h)L\lrp{\frac{1-\bx^\top \by}{h^2}},\label{eq:estimator}
\end{align}
where the normalized kernel $L_h:\Sd\times\Sd\to\R_{\geq 0}$ is based on the kernel $L:\R_{\geq 0}\to\R_{\geq 0}$.
The normalizing constant of $L$ is
\begin{align}
    c_{d,L}(h)^{-1}\defin&\int_{\Sd}L\lrp{\frac{1-\bx^\top \by}{h^2}}\,\sigmad(\rd\bx)
    = h^d \lambda_{h,d}(L)\sim h^d\lambda_{d}(L),\label{eq:equiv}
\end{align}
with $\lambda_{h,d}(L)\defin\om{d-1}\int_{0}^{2h^{-2}}L\lrp{s} s^{d/2-1} (2-sh^2)^{d/2-1}\,\rd s$ and
\begin{align}
    \lambda_{d}(L)\defin 2^{d/2-1}\om{d-1}\int_{0}^\infty L(s) s^{d/2-1} \,\rd s \label{eq:lambdad}
\end{align}
\citep[see][]{Bai1988}. In \eqref{eq:equiv}, $a_n\sim b_n$ denotes $a_n/b_n\to1$ as $n\to\infty$. %

A widely used choice for the kernel in \eqref{eq:estimator} is the von Mises--Fisher (vMF) kernel, defined as $L_{\mathrm{vMF}}(t)\defin e^{-t}, t\ge 0$. This kernel is closely connected to the vMF distribution, whose density is given by
\begin{align}
  \bx\in\Sd\mapsto f_{\mathrm{vMF}}(\bx;\bmu,\kappa)\defin c^\mathrm{vMF}_{d}(\kappa) e^{\kappa\bx^\top\bmu},\quad c_{d}^\mathrm{vMF}(\kappa)\defin\frac{\kappa^{(d-1)/2}}{(2\pi)^{(d+1)/2}\mathcal{I}_{(d-1)/2}(\kappa)}, \label{eq:vmf}
\end{align}
with $\mathcal{I}_\nu$ denoting the modified Bessel function of the first kind of order $\nu$. When considering the vMF kernel, the estimator \eqref{eq:estimator} corresponds to a mixture of von Mises--Fisher densities:
\begin{align}\label{eq:mixture}
    \hat{f}(\bx;h)=\frac{1}{n}\sum_{i=1}^n f_{\mathrm{vMF}}(\bx;\bX_i,1/h^2).
\end{align}

\subsection{Mean integrated squared error}
\label{sec:mise}

The most common way to quantify the performance of the kde \eqref{eq:estimator} is through the mean integrated squared error (MISE), defined as $\mathrm{MISE}(h)\equiv\mathrm{MISE}\{\hat f(\cdot;h)\}\defin\mathbb E\int_{\Sd} \{\hat{f}(\bx;h)-f(\bx)\}^2\,\sigma_d(\rd \bx)$. The MISE can be expressed as ${\rm MISE}(h)={\rm IV}(h)+{\rm ISB}(h)$, where
$$
{\rm IV}(h)\defin\int_{\Sd}\mathbb{V}\mathrm{ar}\{\hat{f}(\bx;h)\}\,\sigma_d(\rd \bx)\quad\text{ and }\quad  {\rm ISB}(h)\defin\int_{\Sd}[\mathbb E\{\hat f(\bx;h)\}-f(\bx)]^2\,\sigma_d(\rd\bx)
$$
represent the integrated variance and the integrated squared bias, respectively. Closed-form analytical expressions for both terms are given next.

To begin with, we define the convolution for a kernel $L_h:\Sd\times\Sd\to\mathbb{R}$ and a density $f:\Sd\to\mathbb{R}$ \citep{Klemela2000} as
\begin{align}
    (L_h*f)(\bx)\defin&\;\int_{\Sd} L_h(\bx,\by)f(\by) \,\sigma_d(\rd \by)\label{eq:Lhf}.%
\end{align}
Then, it readily follows that $\mathbb E\{\hat f(\bx;h)\}=(L_h*f)(\bx)$, so that
\begin{align}
{\rm ISB}(h)&= \int_{\Sd}\{(L_h*f)(\bx)-f(\bx)\}^2\,\sigma_d(\rd\bx)\nonumber\\
&= \int_{\Sd}(L_h*f)(\bx)^2\,\sigma_d(\rd\bx)-2R_{L_h}(f)+R(f),\label{eq:isb}
\end{align}
where we are using the notations $R(f)\defin \int_{\Sd}f(\bx)^2\,\sigmad(\rd \bx)$ and
\begin{align}
     R_{L_h}(f)&\defin \int_{\Sd} (L_h*f)(\bx) f(\bx)\,\sigma_d(\rd \bx)\nonumber\\
     &=\int_{\Sd} \int_{\Sd} L_h(\bx,\by) f(\bx) f(\by)\,\sigma_d(\rd \bx)\,\sigma_d(\rd \by)\nonumber\\%
    &=\mathbb E\{L_h(\bX_1,\bX_2)\}.\label{eq:expker}
\end{align}

Regarding the first term in \eqref{eq:isb}, noting that $L_h(\by,\cdot)$ integrates to one for any fixed $\by\in\Sd$. In direct analogy with \eqref{eq:Lhf}, define
\begin{align}
    \tilde{L}_h(\bx,\by)\defin &\; \{L_h * L_h(\by, \cdot)\}(\bx)=\int_{\Sd} L_h(\bx,\bz) L_h(\by,\bz)\,\sigma_d(\rd \bz).
    \label{eq:Ltilde}
\end{align}
With a slight abuse of notation, sometimes we will also denote $\tilde L_h$ by $L_h*L_h$. Both $\tilde{L}_h(\bx,\cdot)$ and $\tilde{L}_h(\cdot,\by)$ integrate to one on $\Sd$, but $\tilde{L}_h$ is not guaranteed to be a normalized kernel, since generally it cannot be written as a normalized function of $(1-\bx^\top\by)/h^2$. Nevertheless, using the rotation invariance property of $\sigma_d$ it is possible to show that $\tilde L_h(\bx,\by)$ is indeed a function of $\bx^\top\by$, so that it is symmetric in its arguments.

Moreover,
\begin{align*}
    \int_{\Sd} &(L_h*f)(\bx)^2\,\sigma_d(\rd \bx)\\
    &=\int_{\Sd}\int_{\Sd} \tilde L_h(\bx,\by)f(\by)f(\bz)\,\sigma_d(\rd \by)\,\sigma_d(\rd \bz)=R_{\tilde L_h}(f),
\end{align*}
hence, we obtain the compact form
$%
    {\rm ISB}(h)=R_{\tilde L_h}(f)-2R_{L_h}(f)+R(f).
$%

For the integrated variance, note that $\mathbb{V}\mathrm{ar}\{\hat f(\bx;h)\}=n^{-1}\mathbb E\{L_h(x,\bX_i)^2\}-n^{-1}[\mathbb E\{L_h(x,\bX_i)\}]^2$, so that
\begin{align}{\rm IV}(h)=n^{-1}\int_{\Sd}\int_{\Sd}L_h(\bx,\by)^2f(\by)\,\sigma_d(\rd \bx)\,\sigma_d(\rd \by)-n^{-1}\int_{\Sd}(L_h*f)(\bx)^2\,\sigma_d(\rd \bx).\label{eq:IV}\end{align}
The integral in the second term on the right-hand side of \eqref{eq:IV} coincides with the first term in \eqref{eq:isb}.  For the first term in \eqref{eq:IV} we have
\begin{align*}
    \int_{\Sd}\int_{\Sd}L_h(\bx,\by)^2f(\by)\,\sigma_d(\rd \bx)\,\sigma_d(\rd \by)
    =&\;c_{d,L}(h)^2 \int_{\Sd} c_{d,L^2}(h)^{-1} f(\by)\,\sigma_d(\rd \by)\\
    =&\;h^{-d}v_{h,d}(L),
\end{align*}
where $v_{h,d}(L)\defin\lambda_{h,d}(L^2)\lambda_{h,d}(L)^{-2}$. Hence, the integrated variance \eqref{eq:IV} can be simply written as ${\rm IV}(h)=n^{-1}h^{-d}v_{h,d}(L)-n^{-1}R_{\tilde L_h}(f).$

Combining the IV and ISB expressions, we obtain
\begin{align}
{\rm MISE}(h)&=n^{-1}h^{-d}v_{h,d}(L)+(1-n^{-1})R_{\tilde L_h}(f)-2R_{L_h}(f)+R(f).\label{eq:mise2}
\end{align}
With respect to this criterion, the optimal bandwidth $h_{\rm MISE}$ is defined as the minimizer of ${\rm MISE}(h)$.

While \eqref{eq:mise2} provides an exact analytical expression for the MISE, to better elucidate the role of the bandwidth, it is common to resort to asymptotic approximations. The asymptotic form of the MISE relies on the following assumptions:
\begin{enumerate}[label=\textbf{A\arabic{*}}., ref=\textbf{A\arabic{*}}]
    \item The radial extension $\bar{f}\colon\mathbb R^{d+1}\setminus\{\mathbf {0}\}\to\mathbb R_{\geq0}$, defined by $\bar f(\bx)\defin f(\bx/\|\bx\|)$, is bounded, twice continuously differentiable, and all its second partial derivatives are bounded and square integrable. \label{A1}
    \item The kernel $L:\R_{\geq 0}\to\R_{\geq 0}$ is bounded, integrable, and such that $0<\lambda_d(L^k)<\infty$ for $k=1,2$ and $\beta_{d,2}(L)<\infty$, where $\beta_{d,j}(L) \defin 2^{d/2}\omega_{d-1}\int_0^\infty L(s)s^{j/2+d/2-1}\,\rd s$.\label{A2}
    \item $h\equiv h_n$ is a sequence of positive numbers such that $h_n\to 0$ and $nh_n^d\to \infty$ as $n\to\infty$.\label{A3}
\end{enumerate}
Under the previous conditions, it is possible to express
\begin{align}\label{eq:MhAMISE}
{\rm MISE}(h)={\rm AMISE}(h)+o\big(n^{-1}h^{-d}+h^4\big),
\end{align}
with
\begin{align}
   {\rm AMISE}(h)=n^{-1}h^{-d}v_{d}(L)+h^4b_d(L)^2R\lrpbig{\nabla^2 \bar{f}},\label{eq:amise}
\end{align}
where $b_{d}\left(L\right)\defin\beta_{d,2}(L)/\{d\lambda_{d}(L)\}$, $v_{d}(L)\defin\lambda_{d}(L^2)/\lambda_{d}(L)^2$ and, for the Laplacian $\nabla^2 \bar{f}$, we denote $R\lrpbig{\nabla^2 \bar{f}}\defin \int_{\Sd}\big\{\nabla^2 \bar{f}(\bx)\big\}^2\,\sigma_d(\rd\bx)$ \citep[see][]{Garcia-Portugues2013b,Garcia-Portugues:polysphere}.

Whereas $h_{\rm MISE}$ does not admit a closed-form expression in general, it is easy to prove that the unique minimizer of ${\rm AMISE}(h)$ is given by
\begin{align}
h_0=c_0n^{-1/(d+4)},\quad c_0=\lrb{\frac{dv_d(L)}{4b_d(L)^2R(\nabla^2 \bar f)}}^{1/(d+4)}. \label{eq:hamise}
\end{align}
Moreover, it can be shown that those two bandwidths are asymptotically equivalent, in the sense that $h_0/h_{\rm MISE}\to1$ as $n\to\infty$. More precisely, reasoning as in \cite{Marron1987}, under slightly stronger assumptions, the relative rate of convergence in the previous approximation can be found to be
\begin{align}
    \frac{h_0-h_\mathrm{MISE}}{h_\mathrm{MISE}}=O(n^{-2/(d+4)}). \label{eq:Ohamise}
\end{align}

\subsection{Existence of the optimal bandwidth}
\label{sec:existence}

Alternatively, the smoothing parameter can be reparameterized in terms of a concentration parameter $\nu=1/h$. This allows using $\nu=0$, which yields the uniform density as the density estimate, resulting in a perfect, zero-error estimate if the true density is indeed uniform. Hence, the function that measures the error in terms of $\nu$ is denoted ${\rm MISE}2(\nu)\defin{\rm MISE}(1/\nu)$ if $\nu>0$ and ${\rm MISE}2(0)\defin\int_{\Sd}\{\omega_d^{-1}-f(\bx)\}^2\,\sigma_d(\rd \bx)$.

A natural preliminary question is whether an optimal bandwidth can be found; that is, if there exists a value of $\nu_{\rm MISE}\geq0$ that minimizes ${\rm MISE}2(\nu).$ For circular data ($d=1$), \cite{Tenreiro2024} showed that such a minimizer may fail to exist for a related estimator. This stands in contrast to the Euclidean setting, where \cite{Chacon2007} proved that an optimal bandwidth always exists when $f$ is square integrable. Our first result in this section shows that, for the kernel density estimator \eqref{eq:estimator}, existence is also guaranteed under the minimal assumption that $f$ is square integrable.

\begin{theorem}\label{thm:exist1}
    Assume that $f$ is square integrable and $L$ is continuous, $L(0)>0$ and satisfies \ref{A2}. Then, there exists $\nu_{\rm MISE}\geq0$ such that ${\rm MISE}2(\nu_{\rm MISE})\leq{\rm MISE}2(\nu)$, for all $\nu\geq0$.
\end{theorem}

While the previous result ensures the existence of an optimal bandwidth, it does not preclude the possibility that the optimal bandwidth may be degenerate, that is, $\nu_{\rm MISE}=0$. That is the optimal choice for the uniform distribution, though not exclusively so (see Section \ref{sec:num:MvMF}). However, for expressions such as \eqref{eq:Ohamise} to make sense, the optimum needs to be non-degenerate (i.e., $h_{\rm MISE}<\infty$), at least for large enough $n$. The next result provides sufficient conditions under which this holds.

Let $\bmu_k(f)\defin\int_{\mathbb{S}^d}\bx^{\otimes k}f(\bx)\,\sigma_d(\rd \bx)$ denote the $k$th raw moment of $f$, where $\bx^{\otimes k}=\bigotimes_{i=1}^k\bx$ is the $k$th Kronecker power of $\bx$ \citep{Holmquist1988}, and write $m_k(f)\defin\|\bmu_k(f)\|^2$ for its square norm. For the uniform density on $\mathbb{S}^d$, abbreviate their moments and square norms to $\boldsymbol\zeta_k$ and $z_k$, respectively.

\begin{theorem}\label{thm:numise}
    Assume that $L$ is $4$-times continuously differentiable at zero, with $L(0)> 0$.
    \begin{enumerate}[label=\textit{(\alph{*})}, ref=\textit{(\alph{*})}]
        \item If $m_1(f)>0$ and $L'(0)<0$, then for all $n\in\mathbb N$ there exists $\nu_{\rm MISE}>0$ minimizing ${\rm MISE}2(\nu)$.\label{thm1:a}
        \item If $m_1(f)=0$, but $m_2(f)>z_2$  and $L''(0)>0$, then
        for every
        $$
        n>\frac{L'(0)^2z_2}{L(0)L''(0)\{m_2(f)-z_2\}}
        $$
        there exists $\nu_{\rm MISE}>0$ minimizing ${\rm MISE}2(\nu)$.\label{thm1:b}
        \item Assume $m_1(f)=0$, $m_2(f)=z_2$, $m_3(f)=0$ and $m_4(f)>z_4$, and that $L=L_{\mathrm{vMF}}$. Then, there is $n_0\in\mathbb N$ such that for every $n>n_0$
        there exists $\nu_{\rm MISE}>0$ minimizing ${\rm MISE}2(\nu)$.\label{thm1:c}
    \end{enumerate}
\end{theorem}

Theorem \ref{thm:numise} ensures the existence of a non-degenerate optimal bandwidth for a broad class of densities (at least for large enough $n$). Nevertheless, part \ref{thm1:a} does not apply to antipodally symmetric distributions, since they satisfy $\bmu_1(f)=\mathbf 0$. However, it can be shown that for any density $f$ on $\Sd$, the inequality $\|\bmu_2(f)\|^2\geq1/(d+1)$ holds, with the uniform density reaching the lower bound. So, even if $f$ is antipodally symmetric, Theorem \ref{thm:numise} guarantees the existence of a non-degenerate optimal bandwidth as long as the norm of its second-order moment is not the same as that of the uniform distribution. Similarly, even if the first three moments of $f$ match those of the uniform distribution, a non-degenerate optimal bandwidth exists if the fourth moment differs.

In any case, the uniform distribution is not the only one that reaches the lower bounds for $m_2(f)$ and $m_4(f)$. It is possible to construct densities on $\Sd$ different from the uniform, whose second- and fourth-order moments coincide with $z_2$ and $z_4$, respectively. For example, the spherical cardioid distribution \citep{Garcia-Portugues:cardioid} with order larger than four. So, for any such distribution, Theorem \ref{thm:numise} is not enough to guarantee the existence of a non-degenerate optimal bandwidth. We conjecture that such a bandwidth exists (for large enough $n$) whenever $f$ is not the uniform density. However, we will not pursue this generalization in the present work.

\section{Cross-validation bandwidth selector}
\label{sec:cross}

The problem of automatic bandwidth selection involves identifying a criterion that provides a reliable estimate of the MISE, and subsequently selecting the bandwidth that minimizes that criterion. A common simplification of the problem consists of replacing the MISE objective function \eqref{eq:mise2} with
\begin{align}\label{eq:M}
    M(h)\defin n^{-1}h^{-d}v_{h,d}(L)+R_{K_h}(f)
\end{align}
where %
$$
K_h(\bx,\by)\defin\tilde{L}_h(\bx,\by)-2L_h(\bx,\by).
$$
The differences between ${\rm MISE}(h)$ and $M(h)$ are twofold: first, the term $R(f)$ is omitted (because it does not depend on $h$); second, the approximation $1-n^{-1}\sim1$ is employed, which is equivalent to retaining only the first term in the integrated variance. Then, the cross-validation criterion is obtained by estimating the unknown expectation in \eqref{eq:M} using a $U$-statistic  (see also \eqref{eq:expker}), leading to
\begin{align}
    {\rm CV}(h)=n^{-1}h^{-d}v_{h,d}(L)+\textstyle{\binom{n}{2}}^{-1}\displaystyle\sum_{1\leq i<j\leq n}K_h(\bX_i,\bX_j).\label{eq:cv}
\end{align}
It is straightforward to verify that this criterion is an unbiased estimator of $M(h)$. Our goal is to study the properties of $\hat h_{\rm CV}$, the bandwidth that minimizes ${\rm CV}(h)$.

An alternative (although equivalent) approach was proposed in \citet{Hall1987} to motivate the cross-validation bandwidth $\hat h_{\rm CV}$. The performance of the kde $\hat f$ can be alternatively measured through the integrated squared error (ISE), defined as $\mathrm{ISE}(h)\equiv\mathrm{ISE}\{\hat f(\cdot;h)\}\defin\int_{\Sd}\{\hat f(\bx;h)-f(\bx)\}^2\,\sigma_d(\rd\bx)$. This is a random quantity, which evaluates the squared error for the observed sample, and its minimizer is given by $\hat h_{\rm ISE}=\arg\min_{h>0}{\rm ISE}(h)$. It can be shown that
\begin{align}\label{eq:Rfhat}
    \int_{\Sd}\hat f(\bx;h)^2\,\sigma_d(\rd\bx)=n^{-1}h^{-d}v_{h,d}+2n^{-2}\sum_{1\leq i<j\leq n}\tilde L_h(\bX_i,\bX_j).
\end{align}
If the cross term $\int_{\Sd}\hat f(\bx;h)f(\bx)\,\sigma_d(\rd\bx)$ is written as $\mathbb E\{\hat f(\bX_0;h)|\bX_1,\dots,\bX_n\}$, with $\bX_0$ a random variable with density $f$ that is independent of $\bX_1,\dots,\bX_n$, then a cross-validated estimator for this term is given by $n^{-1}\sum_{i=1}^n\hat f_i(\bX_i;h)$, where $\hat f_i$ denotes the kde based on the sample leaving out $\bX_i$. Putting together $ \int_{\Sd}\hat f(\bx;h)^2\,\sigma_d(\rd\bx)-2n^{-1}\sum_{i=1}^n\hat f_i(\bX_i;h)$ as an estimate of ${\rm ISE}(h)-R(f)$ results exactly in the same criterion as in \eqref{eq:cv}, provided the approximation $n^{2}\approx n(n-1)$ is used in \eqref{eq:Rfhat}.

\subsection{Consistency}

The primary objective of this section is to demonstrate the consistency of the cross-validation bandwidth, $\hat h_{\rm CV}$. This means proving that $\hat h_{\rm CV}/h_{\rm MISE}\stackrel{\mathbb P}{\longrightarrow}1$. By denoting $h_M$ the minimizer of $M(h)$, it can be shown that $h_M/h_{\rm MISE}\to1$ as $n\to\infty$ \citep[see][]{chacon2011unconstrained}. Consequently, establishing the consistency of $\hat h_{\rm CV}$ is equivalent to proving that $\hat h_{\rm CV}/h_M\stackrel{\mathbb P}{\longrightarrow}1$. To achieve this, we will rely on fundamental results from $U$-statistics theory, as presented in \cite{lee1990u}.

Given a symmetric function of two variables $\varphi(\bx_1,\bx_2)$ (commonly called kernel) and an iid sample $\bX_1,\dots,\bX_n$ of $d$-dimensional random variables, the statistic
$$
U_n=\textstyle{\binom{n}{2}}^{-1}\displaystyle\sum_{1\leq i<j\leq n}\varphi(\bX_i,\bX_j)
$$
is a $U$-statistic of order $2$. It is an unbiased estimator of $\mathbb E\{\varphi(\bX_1,\bX_2)\}$ and, as shown in \citet[][Section 1.3]{lee1990u}, its variance can be expressed as
\begin{align}\label{eq:varU}
	\mathbb{V}\mathrm{ar}\{U_n\}=\frac{2}{n(n-1)}\mathbb{V}\mathrm{ar}\{\varphi(\bX_1,\bX_2)\}+\frac{4(n-2)}{n(n-1)}\mathbb{C}\mathrm{ov}\{\varphi(\bX_1,\bX_2),\varphi(\bX_1,\bX_3)\}.
\end{align}
Since the cross-validation criterion \eqref{eq:cv} is a $U$-statistic of order two, the above expression provides an analytical expression for its exact variance.

\begin{lemma}\label{lemma:varCV}
	For any $h>0$, the variance of the cross-validation criterion is
	$$\mathbb{V}\mathrm{ar}\{{\rm CV}(h)\}=\frac{2}{n(n-1)}\Big\{R_{(K_h)^2}(f)-R_{K_h}(f)^2\Big\}+\frac{4(n-2)}{n(n-1)}\big\{S_{K_h}(f)-R_{K_h}(f)^2\big\},$$
	where $S_{K_h}(f)\defin \int_{\Sd} (K_h*f)(\bx)^2 f(\bx)\,\sigma_d(\rd \bx)$.
\end{lemma}

Since $\mathbb E\{{\rm CV}(h)\}=M(h)$, a direct consequence of the previous result is that, for a fixed $h>0$, we have
\begin{align*}{\rm CV}(h)-M(h)=O_\mathbb{P}\{\mathbb{V}\mathrm{ar}^{1/2}[{\rm CV}(h)]\}=O_\mathbb{P}(n^{-1/2}).
\end{align*}
However, proving consistency of $\hat h_{\rm CV}$ requires a stronger approximation result. On the one hand, the previous bound must hold uniformly for $h>0$, or at least within an appropriate range of bandwidths. On the other hand, a more tightly concentrated, auxiliary pseudo-criterion is needed to better describe the behavior of $\hat h_{\rm CV}$.

To address this second requirement, we employ the augmented cross-validation (ACV) criterion introduced by \cite{scott1987biased}:
$$
{\rm ACV}(h)={\rm CV}(h)+\frac2n\sum_{i=1}^n\{f(\bX_i)-R(f)\}.
$$
This criterion is not fully data-based, because it involves the unknown density $f$. However, since the additional term does not depend on $h$, the minimizer of ${\rm ACV}(h)$ is also $\hat h_{\rm CV}$. Moreover, $\mathbb E\{{\rm ACV}(h)\}=\mathbb E\{{\rm CV}(h)\}$, so this criterion is also unbiased for $M(h)$. More importantly, as we shall demonstrate,  its variance is of a lower order in $n$. Intuitively, this is a consequence of the fact that ${\rm ACV}(h)$ corresponds to a refined version of the so-called $H$-decomposition of a $U$-statistic, which allows expressing a $U$-statistic as a sum of centered, uncorrelated $U$-statistics with lower-order variances \citep[see][Section 1.6]{lee1990u}.

The exact variance of ${\rm ACV}(h)$ involves computing the covariance between a $U$-statistic of order two and a sample average. Fortunately, this calculation is standard in $U$-statistics theory, as demonstrated in the following result.

\begin{lemma}\label{lemma:varACV}
	For any $h>0$, the variance of the augmented cross-validation criterion is
	\begin{align*}
		\mathbb{V}\mathrm{ar}\{{\rm ACV}(h)\}&=\mathbb{V}\mathrm{ar}\{{\rm CV}(h)\}+4n^{-1}\lrb{\int_{\Sd}  f^3(\bx)\,\sigma_d(\rd \bx) -R(f)^2}\\
		&\quad+8n^{-1}\lrb{\int_{\Sd} (K_h*f)(\bx)f(\bx)^2\,\sigma_d(\rd \bx)-R(f)R_{K_h}(f)}.
	\end{align*}
\end{lemma}

To illustrate that ${\rm ACV}(h)$ is more tightly concentrated around $M(h)$ than ${\rm CV}(h)$, we next derive asymptotic approximations for the variances of the two criteria.

\begin{lemma}\label{lemma:varCVa} Under assumptions \ref{A1}--\ref{A3}, we have:
\begin{enumerate}[label=(\textit{\roman{*}}), ref=(\textit{\roman{*}})]
\item The asymptotic variance of the cross-validation criterion is \label{lemma:varCva:i}
\begin{align*}
	\mathbb{V}\mathrm{ar}\{{\rm CV}(h_n)\}\sim 4n^{-1}\big\{\textstyle \int_{\Sd} f(\bx)^3\,\sigma_d(\rd \bx)-R(f)^2\big\}.
\end{align*}
\item The asymptotic variance of the augmented cross-validation criterion is \label{lemma:varCva:ii}
\begin{align*}
    \mathbb{V}\mathrm{ar}\{{\rm ACV}(h_n)\}\sim 2a_d(L)R(f)n^{-2}h_n^{-d},%
\end{align*}
where
\begin{align*}
	a_d(L)=&\; \lambda_d^{-2}(L)\Bigg[\lambda_d^{-2}(L) \tilde{\gamma}_d  \int_0^{\infty}\left\{\int_0^{\infty} L(r) r^{d/2-1}\varphi_d(L, r, s)\, \rd r\right\}^2s^{d/2-1}\, \rd s\\
&-4\lambda_d^{-1}(L)\gamma_d \int_0^{\infty}\int_0^{\infty} L(r) L(s)r^{d/2-1} s^{d/2-1} \varphi_d(L, r, s)\, \rd r\, \rd s+4\lambda_d(L^2)\Bigg],
\end{align*}
with
\begin{align}
\varphi_d(L, r, s)\defin & \begin{cases}L\left(r+s-2(r s)^{1/2}\right)+L\left(r+s+2(r s)^{1/2}\right), & d=1, \\
\int_{-1}^1\left(1-\theta^2\right)^{(d-3)/2} L\left(r+s-2 \theta(r s)^{1/2}\right) \rd \theta, & d \geq 2,\end{cases}\label{eq:phid}\end{align}
\begin{align*}
\gamma_d=  \begin{cases}1, & d=1, \\
\om{d-1} \om{d-2} 2^{d-2}, & d \geq 2,\end{cases}\quad \text{and} \quad
\tilde{\gamma}_d=  \begin{cases}2^{-1/2}, & d=1, \\
\om{d-1}\om{d-2}^22^{(3d-6)/2}, & d \geq 2 .\end{cases}
\end{align*}
\end{enumerate}
\end{lemma}

As a consequence of Lemma \ref{lemma:varCVa}, we have $\mathbb{V}\mathrm{ar}\{{\rm ACV}(h_n)\}=o[\mathbb{V}\mathrm{ar}\{{\rm CV}(h_n)\}]$ whenever $nh_n^d\to\infty$. This confirms that ${\rm ACV}(h_n)$ is indeed more concentrated around its mean $M(h_n)$ than ${\rm CV}(h_n)$, as claimed.

Next, in view of \eqref{eq:hamise}, we focus on bandwidths within the interval $[\varepsilon n^{-1/(d+4)},M n^{-1/(d+4)}]$ for arbitrary constants $0<\varepsilon<M<\infty$. For any bandwidth of the form $h_n=c n^{-1/(d+4)}$, part \ref{lemma:varCva:i} of Lemma \ref{lemma:varCVa}, combined with \eqref{eq:MhAMISE}, implies that
\begin{align*}
	{\rm CV}(cn^{-1/(d+4)})&=M(cn^{-1/(d+4)})+O_\mathbb{P}(n^{-1/2})\\
	&={\rm AMISE}(cn^{-1/(d+4)})-R(f)+o(n^{-4/(d+4)})+O_\mathbb{P}(n^{-1/2}).
\end{align*}

However, applying part \ref{lemma:varCva:ii} of Lemma \ref{lemma:varCVa} to the augmented criterion yields
\begin{align}
	{\rm ACV}(cn^{-1/(d+4)})&=M(cn^{-1/(d+4)})+O_\mathbb{P}(n^{-(d+8)/(2d+8)})\nonumber\\
	&={\rm AMISE}(cn^{-1/(d+4)})-R(f)+o(n^{-4/(d+4)})+O_\mathbb{P}(n^{-(d+8)/(2d+8)}),\label{eq:ACV}
\end{align}
where the stochastic error term $O_\mathbb{P}(n^{-(d+8)/(2d+8)})$ is now of smaller order than the approximation error $o(n^{-4/(d+4)})$ for all $d$. This sets the stage for establishing the consistency of the cross-validation bandwidth.

\begin{theorem}\label{th:consist}
Let $\hat h_{\rm CV}$ be the minimizer of ${\rm CV}(h)$ over $[\varepsilon n^{-1/(d+4)},M n^{-1/(d+4)}]$ for arbitrary constants $0<\varepsilon<M<\infty$. Then, under assumptions \ref{A1}--\ref{A3}, $\hat h_{\rm CV}/h_{\rm MISE}\stackrel{\mathbb P}{\longrightarrow}1$.
\end{theorem}

The previous result establishes the consistency of $\hat h_{\rm CV}$ when the search is restricted to the interval $[\varepsilon n^{-1/(d+4)}, M n^{-1/(d+4)}]$. This argument can be extended to a search over the entire range $(0, \infty)$, similar to the approach of \cite{devroye1989double} for a different bandwidth selector. However, proving consistency in this broader setting requires a substantially more detailed and technical analysis.

\subsection{Convergence rate}
\label{sec:conv}

The derivation of the relative rate of convergence of $\hat h_{\rm CV}$ parallels the approach used in \cite{HM87} or \cite{Park1990}. For this goal, some slightly stronger smoothness conditions are needed; specifically, it will be necessary to assume:
\begin{enumerate}[label=\textbf{A\arabic{*}*}., ref=\textbf{A\arabic{*}*}]
    \item The radial extension $\bar{f}$ is bounded, 4-times continuously differentiable, and all its 4th-order partial derivatives are bounded and square integrable. \label{A1b}
    \item The kernel $L:\R_{\geq 0}\to\R_{\geq 0}$ is bounded, integrable, twice continuously differentiable, and such that $0<\lambda_d(L^k)<\infty$ for $k=1,2$ and $\beta_{d,4}(L)<\infty$. Moreover, the function $G(t)\defin L'(t) t$ also satisfies $0<\lambda_d(G^k)<\infty$ for $k=1,2$ and $\beta_{d,4}(G)<\infty$.\label{A2b}
\end{enumerate}
With those assumptions, ${\rm CV}(h)$ is differentiable and $M(h)$ is twice differentiable with respect to $h$. Then a Taylor expansion yields
\begin{align}
    0={\rm CV}'(\hat h_{\rm CV})&=M'(\hat h_{\rm CV})+({\rm CV}-M)'(\hat h_{\rm CV})\nonumber\\
    &=M''(\tilde h)(\hat h_{\rm CV}-h_M)+({\rm CV}-M)'(\hat h_{\rm CV}),\label{eq:Tay1}
\end{align}
where $\tilde h$ is a random value that lies between $\hat h_{\rm CV}$ and $h_M$. From \eqref{eq:Tay1} we can express the relative error as
\begin{align}\label{eq:re}
    \frac{\hat h_{\rm CV}-h_M}{h_M}=-\frac{({\rm CV}-M)'(\hat h_{\rm CV})}{h_MM''(\tilde h)}.
\end{align}

Using the consistency of $\hat h_{\rm CV}$ and the equivalence between $h_0$, $h_M$ and $h_{\rm MISE}$, we have $\hat h_{\rm CV}/h_0\stackrel{\mathbb P}{\longrightarrow}1$. This implies that $\tilde h/h_0\stackrel{\mathbb P}{\longrightarrow}1$ as well, which leads to $h_M M''(\tilde h)=h_0M''(h_0)\{1+o_{\mathbb P}(1)\}$. Similarly, reasoning as in \citet[][p. 1928]{Jones1991}, we may replace $\hat h_{\rm CV}$ by $h_0$ on the right-hand side of \eqref{eq:re}, and $h_M$ by $h_{\rm MISE}$ on the left-hand side. Thus we obtain the asymptotic approximation
\begin{align}\label{eq:re2}
    \frac{\hat h_{\rm CV}-h_{\rm MISE}}{h_{\rm MISE}}=-\frac{({\rm CV}-M)'(h_0)}{h_0M''(h_0)}\{1+o_{\mathbb P}(1)\}.
\end{align}

On the other hand, from \eqref{eq:amise} it follows that

\begin{align*}
    {\rm AMISE}'(h)&=-dn^{-1}h^{-d-1}v_d(L)+4h^3b_d(L)^2R\lrpbig{\nabla^2 \bar{f}}\quad\text{ and}\\
    {\rm AMISE}''(h)&=d(d+1)h^{-d-2}v_d(L)+12h^2b_d(L)^2R\lrpbig{\nabla^2 \bar{f}}.
\end{align*}
Since $h_0=c_0n^{-1/(d+4)}$, and asymptotically $M(h)$ differs from ${\rm MISE}(h)$ only by a constant shift, then we obtain %
\begin{align}\label{eq:h0M2}
    h_0M''(h_0)\sim c_1 n^{-3/(d+4)},\quad\text{with }c_1=d(d+1)c_0^{-d-1}v_d(L)+12b_d(L)^2R(\nabla^2 \bar f)c_0^3.
\end{align}

Therefore, to characterize the asymptotic behavior of the relative error of $\hat h_{\rm CV}$, it suffices to study the random variable $({\rm CV}-M)'(h_0)$ appearing in the numerator of \eqref{eq:re2}.

Next, let us derive more explicit expressions for the derivative of $({\rm CV}-M)(h)$. Recall that
\begin{align*}
    ({\rm CV}-M)(h)=\textstyle{\binom{n}{2}}^{-1}\displaystyle\sum_{1\leq i<j\leq n}K_h(\bX_i,\bX_j)-R_{K_h}(f),
\end{align*}
where $R_{K_h}(f)=\mathbb E\{K_h(\bX_1,\bX_2)\}$.

First, define
\begin{align*}
    \nu_h(\bx,\by)\defin\frac{\partial}{\partial h}  K_h(\bx,\by)=\frac{\partial}{\partial h} \lrb{\tilde{L}_h(\bx,\by)-2L_h(\bx,\by)},
\end{align*}
with $ \tilde{L}_h(\bx,\by)=\int_{\Sd} L_h(\bx,\bz) L_h(\by,\bz)\,\sigma_d(\rd \bz)$. Since the map $h\mapsto L_h(\bx,\by)$ is continuous, and the kernel is assumed to be bounded, the Leibniz integral rule in its measure-theoretic version gives the derivative of the kernel normalizing constant:
\begin{align*}
    \frac{\partial}{\partial h} c_{d,L}(h)
    =&\;\frac{-1}{c_{d,L}(h)^{-2}} \frac{\partial }{\partial h} c_{d,L}(h)^{-1}\nonumber\\
    =&\;\frac{2c_{d,L}(h)^{2}}{h} \int_{\Sd}L'\lrp{\frac{1-\bx^\top \by}{h^2}}\frac{1-\bx^\top \by}{h^2}\,\sigmad(\rd\bx)\nonumber\\
    =&\;\frac{2c_{d,L}(h)^{2}}{h} \int_{\Sd}G\lrp{\frac{1-\bx^\top \by}{h^2}}\,\sigmad(\rd\bx)\nonumber\\
    =&\;\frac{2c_{d,L}(h)^{2}}{h c_{d,G}(h)}.
\end{align*}

Now, we compute
\begin{align}
    \frac{\partial}{\partial h} L_{h}(\bx,\by)&=\frac{\partial}{\partial h}\lrb{c_{d,L}(h) L\lrp{\frac{1-\bx^\top \by}{h^2}}}\nonumber\\
    &= \frac{2c_{d,L}(h)^{2}}{h c_{d,G}(h)}L\lrp{\frac{1-\bx^\top \by}{h^2}}-2c_{d,L}(h)L'\lrp{\frac{1-\bx^\top \by}{h^2}}\frac{1-\bx^\top \by}{h^3}\nonumber\\
    &= \frac{2c_{d,L}(h)}{h c_{d,G}(h)} \lrb{c_{d,L}(h)L\lrp{\frac{1-\bx^\top \by}{h^2}}-c_{d,G}(h)G\lrp{\frac{1-\bx^\top \by}{h^2}}}\nonumber\\
    &\indef \frac{2c_{d,L}(h)}{h c_{d,G}(h)} (L_h-G_h)(\bx,\by).
    \label{eq:derLht}
\end{align}

Next, using the Leibniz integral rule again, we have
\begin{align}
   \frac{\partial}{\partial h}  \tilde{L}_h(\bx,\by)=&\; \int_{\Sd}   \frac{\partial}{\partial h} \lrb{L_{h}(\bx,\bz)}  L_{h}(\by,\bz)+  \frac{\partial}{\partial h}\lrb{L_{h}(\by,\bz)}  L_{h}(\bx,\bz)\,\sigma_d(\rd \bz)\nonumber\\
    =&\; \frac{4c_{d,L}(h)}{h c_{d,G}(h)}\int_{\Sd} (L_h-G_h)(\bx,\bz)L_h(\by,\bz)\,\sigma_d(\rd \bz),\label{eq:Gh*Lh}
\end{align}
where in \eqref{eq:Gh*Lh} we have used that
\begin{align*}
    \int_{\Sd}  (L_h-G_h)(\bx,\bz)L_{h}(\by,\bz)\,\sigma_d(\rd \bz)=\int_{\Sd}  (L_h-G_h)(\by,\bz)L_{h}(\bx,\bz)\,\sigma_d(\rd \bz),
\end{align*}
since the previous expression is symmetric in $(\bx,\by)$ as it was shown for $\tilde L_h(\bx,\by)$ after \eqref{eq:Ltilde}.

Putting together \eqref{eq:derLht} and \eqref{eq:Gh*Lh},
\begin{align}
   \nu_h(\bx,\by)=& \frac{\partial}{\partial h}\lrb{\tilde{L}_h(\bx,\by)-2L_h(\bx,\by)}\nonumber\\
   =&\frac{4c_{d,L}(h)}{h c_{d,G}(h)}\big[\{(L_h-G_h)*L_h\}(\bx,\by)-(L_h-G_h)(\bx,\by)\big].\label{eq:nu}
\end{align}
Note that from \eqref{eq:nu}, the function $\nu_h(\bx,\by)$ can be written as $\nu_h(\bx,\by)=u_h(\bx^\top\by)$ for a certain function $u_h$. However, it is not immediately clear that $\nu_h(\bx,\by)$ is a function of $(1-\bx^\top\by)/h^2$.

Finally, using the Leibniz integral rule once again, it can be shown that $\frac{\partial}{\partial h}R_{K_h}(f)=\mathbb E\{\nu_h(\bX_1,\bX_2)\}$, which leads to
\begin{align}
    ({\rm CV}-M)'(h)=\textstyle{\binom{n}{2}}^{-1}\displaystyle\sum_{1\leq i<j\leq n}\nu_h(\bX_i,\bX_j)-\mathbb E\{\nu_h(\bX_1,\bX_2)\}.\label{eq:CVder}
\end{align}

From \eqref{eq:CVder}, it is clear that $\mathbb E\{{\rm CV}'(h)\}=M'(h)$ for all $h>0$. The next result establishes the asymptotic properties of $({\rm CV}-M)'(h_0)$.

\begin{lemma}\label{lemma:varCVd} Under assumptions \ref{A1b}--\ref{A2b} and \ref{A3}, it follows that
\begin{align*}
     \mathbb{V}\mathrm{ar}\{{\rm CV}'(h_0)\}\sim 2\sigma_0^2(L) R(f)n^{-2}h_0^{-d-2},
\end{align*}
where
\begin{align}
    \sigma^2_0(L)=&\;16\Bigg\{\frac{\lambda_d(G^2)}{\lambda_d(L)^2}-\frac{2{\gamma}_d}{\lambda_d(L)^3}\int_0^{\infty}\int_0^\infty L(u)L(v) \varphi_{12,d}(L, u, v)u^{d/2}v^{d/2} \,\rd u \,\rd v\nonumber\\
&+\frac{\tilde{\gamma}_d}{\lambda_d(L)^4}\int_0^{\infty}\lrc{\int_0^\infty L(u) \varphi_{1,d}(L, u, v)u^{d/2}\, \rd u}^2v^{d/2-1}  \,\rd v\Bigg\},\label{eq:sigma2}
\end{align}
with $\varphi_{1,d}(L, u, v)=\frac{\partial}{\partial u}\varphi_d(L, u, v)$ and $\varphi_{12,d}(L, u, v)=\frac{\partial}{\partial v}\frac{\partial}{\partial u}\varphi_d(L, u, v)$.

Moreover, $[\mathbb{V}\mathrm{ar}\{{\rm CV}'(h_0)\}]^{-1/2}\{({\rm CV}-M)'(h_0)\}\stackrel{d}{\longrightarrow}\mathcal{N}(0,1)$.
\end{lemma}

The following result presents the asymptotic behavior of the cross-validation bandwidth selector $\hat h_{\rm CV}$.
\begin{theorem}\label{th:rate}Under assumptions \ref{A1b}--\ref{A2b} and \ref{A3}, it follows that
\begin{align}
    n^{d/(2d+8)}(\hat h_{\rm CV}-h_{\rm MISE})/h_{\rm MISE}\stackrel{d}{\longrightarrow} \mathcal{N}\big(0,\sigma_d^2(L,f)\big),\label{eq:asympvar}
\end{align}
where $\sigma_d^2(L,f)=\tau_d(L)\rho_d(f)$,
with the density and kernel contributions to the asymptotic variance, respectively, given by
\begin{align}
    \rho_d(f)\defin&\; \frac{R(f)}{R(\nabla^2 \bar f)^{d/(d+4)}},\label{eq:rho_f}\\
    \tau_d(L) \defin&\; \frac{1}{[2^{d-4} d^{d+8}]^{1/(d+4)}(d+4)^2} \times \frac{\sigma^2_0(L)}{[v_d(L)^{d+8} b_d(L)^{2d}]^{1/(d+4)}},\label{eq:tau_L}
\end{align}
and $\sigma^2_0(L)$ given in \eqref{eq:sigma2}. Hence,
 \begin{align*}
    \frac{\hat h_{\rm CV}-h_{\rm MISE}}{h_{\rm MISE}}=O_\mathbb{P}(n^{-d/(2d+8)}).
\end{align*}
\end{theorem}

\begin{corollary} \label{cor:vmf:kernel}
For the vMF kernel,
\begin{align}
    \tau_d(L_\mathrm{vMF})=\frac{2^{(5d+4)/(d+4)} \pi^{2d/(d+4)}}{d^{(d+8)/(d+4)}} (1+2^{-(d/2+2)}-2(3/2)^{-(d/2+2)}) \frac{d(d+2)}{(d+4)^2}.\label{eq:tau_vmf}
\end{align}
\end{corollary}

As $d\to\infty$, $\tau_d(L_\mathrm{vMF})\sim 32\pi^2/d$, and hence the contribution of the vMF kernel to the asymptotic variance becomes increasingly smaller for large dimensions. However, this convergence is slow and non-monotone. Indeed, the function $d\mapsto\tau_d(L_\mathrm{vMF})$ increases monotonically up to its maximum at dimension $d^*=28$, after which it decreases monotonically (see Figure \ref{fig:sigma2}).

\begin{corollary} \label{cor:vmf:dens}
Denote $\rho_d(\kappa)\defin \rho_d(f_\mathrm{vMF}(\cdot;\bmu,\kappa))$ for the vMF density. Then:
\begin{align}
    \rho_d(\kappa)&=\mathcal{I}_{(d-1)/2}(2\kappa)\label{eq:rho_f:vmf}\\
    &\!\!\!\!\times \lrb{4 \pi^{2(1+1/d)} d (\mathcal{I}_{(d-1)/2}(\kappa))^{8/d}  \kappa^{2/d - 1} \left[2 d \mathcal{I}_{(d+1)/2}(2 \kappa)+(d+2) \kappa \mathcal{I}_{(d+3)/2}(2 \kappa)\right]}^{-d/(d+4)}\nonumber
\end{align}
and
\begin{enumerate}[label=(\textit{\roman{*}}), ref=(\textit{\roman{*}})]
    \item $\rho_d(\kappa)\sim\lrb{4\pi^{2} d(d+2)}^{-d/(d+4)}$ as $\kappa\to\infty$; \label{cor:vmf:dens:kappa}
    \item $\rho_d(\kappa)\sim d\lrb{4\pi^2e^2 \kappa^2}^{-1}$ as $d\to\infty$. \label{cor:vmf:dens:d}
\end{enumerate}
\end{corollary}

When $\kappa=0$, as already anticipated from \eqref{eq:rho_f}, $\rho_d(\kappa)=\infty$, and consequently the asymptotic variance in \eqref{eq:asympvar} is infinite, which is expected because in the uniform density $h_\mathrm{MISE}=\infty$. For any fixed $d>1$, numerical evaluation shows that $\kappa\mapsto \rho_d(\kappa)$ monotonically decreases toward the asymptote signaled in \ref{cor:vmf:dens:kappa} as $\kappa\to\infty$ (see Section \ref{sec:add} of the SM). However, for $d=1$, the function has a global minimum at $\kappa^*\approx1.0917$ such that $\rho_d(\kappa^*)/\rho_d(\infty)=0.8829$; i.e., there exists a finite concentration for the vMF density that reduces the asymptotic variance of an arbitrarily large concentration by $88.29\pct$. When $d\to\infty$, it follows from \ref{cor:vmf:dens:d} that $\sigma_d^2(L_\mathrm{vMF},\kappa)\sim 8/(e\kappa)^2$, hence providing the neat large-$d$ approximation:
\begin{align*}
    \mathbb{V}\mathrm{ar}\left\{\frac{\hat h_{\rm CV}-h_{\rm MISE}}{h_{\rm MISE}}\right\}\approx \frac{8}{(e\kappa)^2 n}.
\end{align*}

\begin{figure}[htb!]
    \centering
    \includegraphics[width=\linewidth]{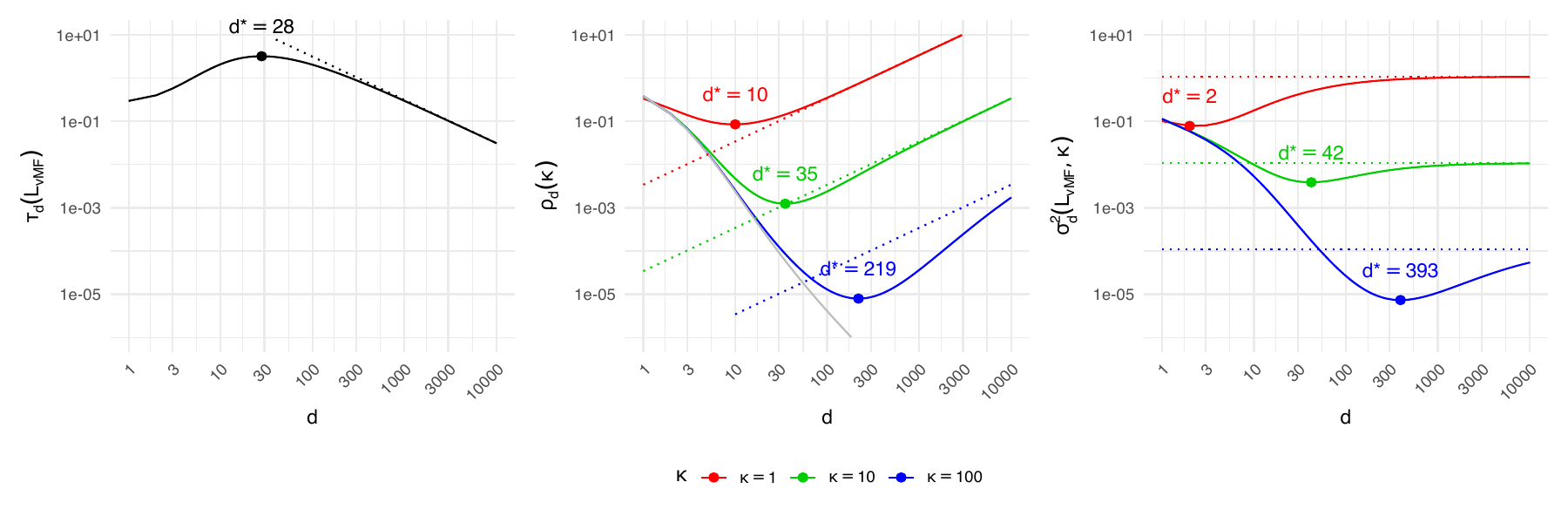}
    \begin{subfigure}{0.33\linewidth}
    \caption{ $d\mapsto\tau_d(L_\mathrm{vMF})$.\label{fig:sigma2:tau}}
    \end{subfigure}%
    \begin{subfigure}{0.33\linewidth}
    \caption{$d\mapsto\rho_d(\kappa)$.\label{fig:sigma2:rho}}
    \end{subfigure}%
    \begin{subfigure}{0.33\linewidth}
    \caption{$d\mapsto \sigma^2_d(L_\mathrm{vMF},\kappa)$.\label{fig:sigma2:sigma2}}
  \end{subfigure}%
  \caption{Asymptotic variance functionals as a function of the dimension $d$. Figure \ref{fig:sigma2:tau} shows the contribution of the vMF kernel to the asymptotic variance, while Figure \ref{fig:sigma2:rho} gives the vMF density contribution for several concentrations $\kappa$. Figure \ref{fig:sigma2:sigma2} collects the resulting asymptotic variance factor $\sigma^2_d(L_\mathrm{vMF},\kappa)$. Dotted lines represent the asymptotic approximations as $d\to\infty$. The gray curve in Figure \ref{fig:sigma2:rho} shows $d\mapsto\rho_d(\phi)$ from Section \ref{sec:comp}. Global extrema of the curves are highlighted with dots. Both axes are $\log_{10}$-scaled.\label{fig:sigma2}}
\end{figure}

Figure \ref{fig:sigma2} depicts the asymptotic variance functionals and their asymptotic approximations as a function of the dimension and concentration, showing, e.g., that $\sigma_d^2(L_\mathrm{vMF},\kappa)$ achieves its minimum at dimension $d^*=42$ when $\kappa=10$.

\subsection{Comparison with the Euclidean case}
\label{sec:comp}

As noted in Section \ref{sec:intro}, the cross-validation bandwidth selector has been extensively studied in the Euclidean case. However, in the multivariate setting, the limit distribution of its relative error has been characterized only with respect to the ISE-optimal bandwidth \citep{HM87} or the AMISE-optimal bandwidth \citep{Sain1994}. \citet{Jones1992} indeed comments on the relative rate of convergence with respect to the MISE-optimal bandwidth (without proof), but does not specify the limiting distribution of the relative error.

Using arguments analogous to those for directional data, the asymptotic behavior of the relative error of the cross-validation bandwidth with respect to the MISE-optimal bandwidth in the Euclidean setting is obtained in \eqref{eq:ANeu}.

In this section, with a slight abuse of notation, let $\bX_1,\dots,\bX_n$ be an iid sample from a density $f\colon\mathbb R^d\to\mathbb R$. Denote $\tilde f(\bx;h)=n^{-1}\sum_{i=1}^nW_h(\bx-\bX_i)$ the kde in this setting, where $W_h(\bx)=W(\bx/h)/h^d$ and $W\colon\mathbb R^d\to\mathbb R$ is a spherically symmetric kernel satisfying $\int_{\mathbb R^d}W(\bx)\,\rd\bx=1$, with finite $R(W)=\int_{\mathbb R^d}W(\bx)^2\,\rd\bx$ and $\mu_2(W)=\int_{\mathbb R^d}x_i^2W(\bx)\,\rd\bx$. Here, let $h_{\rm MISE}$ be the minimizer of $\mathbb E\int_{\mathbb R^d}\{\tilde f(\bx;h)-f(\bx)\}^2\,\rd\bx$ and consider the cross-validation criterion
\begin{align}\label{eq:cveu}
    {\rm CV}(h)=n^{-1}h^{-d}R(W)+{\textstyle\binom{n}{2}}^{-1}\sum_{1\leq i<j\leq n}(W_h*W_h-2W_h)(\bX_i-\bX_j),
\end{align}
where $*$ denotes the usual convolution product of functions on $\mathbb R^d$ \citep[see][Equation (3.7)]{Chacon2018}.
Under some regularity conditions on $f$ and $W$, the bandwidth $\tilde h_{\rm CV}$ that minimizes \eqref{eq:cveu} satisfies
\begin{align}\label{eq:ANeu}
    n^{d/(2d+8)}(\tilde h_{\rm CV}-h_{\rm MISE})/h_{\rm MISE}\stackrel{d}{\longrightarrow}{\mathcal N}(0,\sigma_{\rm CV}^2),
\end{align}
where
\begin{align}\label{eq:sig2CV}
    \sigma_{\rm CV}^2\equiv\sigma_{\rm CV}^2(W,f)=2R(\rho)R(f)(d+4)^{-2}[\mu_2(W)^{2d}R(\nabla^2 f)^d\{dR(W)\}^{d+8}]^{-1/(d+4)},
\end{align}
and $\rho(\bx)=\bx^\top\D(W*W-2W)(\bx)$, with $\D$ denoting the gradient operator on $\mathbb R^d$. A brief sketch of the proof of \eqref{eq:ANeu} is provided in Section \ref{sec:euclidean} of the SM.

For $d=1$, equation \eqref{eq:sig2CV} reduces to the one appearing in \citet[][Theorem 3.1]{Park1990} for the univariate case. Moreover, for the standard Gaussian kernel $W_\mathrm{G}(\bx)=(2\pi)^{-d/2}\exp(-\frac{1}{2}\bx^\top\bx)$ it is shown in Section \ref{sec:euclidean} of the SM that
\begin{align}\label{eq:Rrho}
    R(\rho_\mathrm{G})
    &=2^{-d}\pi^{-d/2}d(d+2)(1+2^{-(d/2+2)}-2(3/2)^{-(d/2+2)}),
\end{align}
which coincides with the formula provided in Corollary 6.4.1 of \cite{Aldershof1995} for the case $d=1$. Remarkably, equation \eqref{eq:Rrho} coincides with the expression obtained here for $\sigma_0^2(L_\mathrm{vMF})$ in the directional context, suggesting a close correspondence between the Gaussian kernel for Euclidean data and the vMF kernel for directional data. As a consequence, the contribution of the Gaussian kernel to the asymptotic variance $\sigma_{\rm CV}^2(W,f)$ coincides exactly with that of the vMF kernel given in \eqref{eq:tau_vmf}. Using $\mu_2(W_\mathrm{G})=1$ and $R(W_\mathrm{G})=(4\pi)^{-d/2}$, we obtain
\begin{align*}
    \tau_d(W_\mathrm{G})\defin&\; 2R(\rho_\mathrm{G})(d+4)^{-2}[\mu_2(W_\mathrm{G})^{2d}\{dR(W_\mathrm{G})\}^{d+8}]^{-1/(d+4)}\\
    =&\; 2^{-d+d(d+8)/(d+4)+1}\pi^{-d/2+d(d+8)/(2d+8)} \frac{d(d+2)}{(d+4)^2}\\
    &\times d^{-(d+8)/(d+4)} (1+2^{-(d/2+2)}-2(3/2)^{-(d/2+2)})\\
    =&\; \frac{2^{(5d+4)/(d+4)}\pi^{2d/(d+4)}}{d^{(d+8)/(d+4)}} (1+2^{-(d/2+2)}-2(3/2)^{-(d/2+2)})\frac{d(d+2)}{(d+4)^2} \\
    =&\;\tau_d(L_\mathrm{vMF}).
\end{align*}

The contribution to the asymptotic variance for a Gaussian density $\phi_{\sigma^2\bI_d}(\cdot-\bmu)$ follows from the second equation in \citet[page 119]{Chacon2018}, which entails that, if $\bZ\sim\mathcal{N}_d(\zero,\bI_d)$, then
\begin{align*}
  \tr{\int_{\R^d}\D^{\otimes r}\phi_{\sigma^2\bI_d}(\bx-\bmu) \D^{\otimes r}\phi_{\sigma^2\bI_d}(\bx-\bmu)^\top\,\rd \bx}=(2 \pi)^{-d / 2}|2\sigma^2\bI_d|^{-1/2}\mathbb{E}[\{\bZ^\top(2\sigma^2\bI_d)^{-1}\bZ\}^r].
\end{align*}

For $r=2$,
\begin{align*}
R(\nabla^2 \phi_{\sigma^2\bI_d}(\cdot-\bmu))
&=2^{-(d+2)} \pi^{-d / 2} \sigma^{-(d + 4)} d(d + 2)
\end{align*}
and for $r=0$,
\begin{align*}
R(\phi_{\sigma^2\bI_d}(\cdot-\bmu))=(4\pi\sigma^2)^{-d/2}.
\end{align*}
Merging these two results gives
\begin{align*}
    \rho_d(\phi)\equiv \rho_d(\phi_{\sigma^2\bI_d}(\cdot-\bmu))\defin&\;\frac{R(\phi_{\sigma^2\bI_d}(\cdot-\bmu))}{R(\nabla^2  \phi_{\sigma^2\bI_d}(\cdot-\bmu))^{d/(d+4)}}
    =(2\pi)^{-2d/(d+4)}[d(d + 2)]^{-d/(d+4)},
\end{align*}
which, unlike $\rho_d(\kappa)$ in \eqref{eq:rho_f:vmf}, does not depend on the shape parameter of the underlying density. Also, when $d\to\infty$, $\rho_d(\phi)\sim (4\pi^2 d^2)^{-1}$ and therefore $\sigma_d^2(W_\mathrm{G},\phi) \sim 
(2d)^{-3}$, which sharply contrasts with the large-$d$ stabilization for the vMF density. Figure \ref{fig:sigma2:rho} shows that $\rho_d(\phi)$ and $\rho_d(\kappa)$ are close for low dimensions, depending on $\kappa$, with larger concentrations keeping the closeness in higher dimensions. Note this result is coherent with the high-concentration Gaussian limit of the vMF, related with the fact that $\lrp{\bX\,\big|\norm{\bX}=1}\sim \mathrm{vMF}\lrp{\bmu/\norm{\bmu},\norm{\bmu}\sigma^{-2}}$ for $\bX\sim \mathcal{N}_{d+1}\lrp{\bmu,\sigma^2\bI_{d+1}}$, with $\bmu\in\R^{d+1}\backslash\{\mathbf{0}\}$ and $\sigma>0$. For smaller $\kappa$'s with respect to the dimension $d$, the difference between $\rho_d(\phi)$ and $\rho_d(\kappa)$ can be attributed to the spherical geometry inherent to the vMF distribution.

\section{Numerical experiments}
\label{sec:num}

In this section, we conduct several numerical experiments with two main goals. First, in Sections \ref{sec:num:vMF}--\ref{sec:num:MvMF}, we empirically evaluate Theorem \ref{th:rate}, i.e., that
\begin{align}
    R(n,d)\defin\frac{\hat{h}_\mathrm{CV}-h_\mathrm{MISE}}{h_\mathrm{MISE}}=O_\mathbb{P}\lrp{n^{\beta_{*}(d)}},\text{ with } \beta_*(d)\defin -\frac{d}{2d+8}, \label{eq:rate}
\end{align}
as $n$ diverges to infinity, for different data generating processes and dimensions. Second, in Sections \ref{sec:num:PI}--\ref{sec:num:dens}, we compare the convergence rates of $\hat{h}_\mathrm{CV}$ in \eqref{eq:rate} with those of plug-in approaches based on AMISE and MISE criteria, and evaluate the corresponding $\mathrm{ISE}\{\hat{f}(\cdot;\hat{h})\}$ for bandwidth selectors $\hat{h}$ obtained via cross-validation and plug-in methods.

For the first goal, we define the expectation $\mathrm{e}(n,d)\defin\mathbb E[R(n,d)]$, which is $o(1)$ according to Theorem \ref{th:consist}, and the root mean squared error (RMSE) $\mathrm{rmse}(n,d)\defin\allowbreak\sqrt{\mathbb{E}[R(n,d)^2]}$, which is $O\lrp{n^{\beta_*(d)}}$. Under \eqref{eq:rate}, it follows that there exist $n_0\geq1$ and $C>0$ such that $\mathrm{rmse}(n,d)\leq C n^{\beta_*(d)}$ for all $n\geq n_0$. Hence, in particular,
\begin{align*}
    \log_2(\mathrm{rmse}(n,d))\leq \log_2(C)+\beta_*(d)\log_2(n)
\end{align*}
for $n\geq n_0$, which in turn implies that in the population least-squares fit
\begin{align*}
    (\alpha(d),\beta(d))=\arg\min_{(a,b)\in\mathbb{R}^2} \sum_{n\geq n_0} (\log_2(\mathrm{rmse}(n,d))-a-b\log_2(n))^2
\end{align*}
the slope must satisfy $\beta(d)\leq \beta_*(d)$.

We investigate the behavior of the statistics $\{R^{(j)}(n,d)\defin(\hat{h}^{(j)}_\mathrm{CV}-h_\mathrm{MISE})/h_\mathrm{MISE}\}_{j=1}^M$ for sample sizes $n=\lfloor 2^\ell\rfloor$, $\ell\in L\defin\{5,5.5,\ldots,13\}$ and dimensions $d=1,2,\ldots,10$, using $M=10,\!000$ Monte Carlo repetitions. From these statistics, we compute the robust estimates $\hat{\mathrm{e}}(n,d)$ and $\widehat{\mathrm{rmse}}(n,d)$ that represent the plug-in mean and RMSE obtained after trimming the $5\pct$ most extreme observations of the sample. Using these statistics, we investigate in Sections \ref{sec:num:vMF}--\ref{sec:num:MvMF} three aspects: (\textit{i}) the curves $n\mapsto\hat{\mathrm{e}}(n,d)$ and $n\mapsto\widehat{\mathrm{rmse}}(n,d)$, for varying $d$; (\textit{ii}) the estimated slopes $\hat{\beta}(d)$ obtained by performing a linear regression of $\log_2(\widehat{\mathrm{rmse}}(n,d))$ onto $\log_2(n)$ for the sample $\{(\log_2(\widehat{\mathrm{rmse}}(\lfloor 2^\ell\rfloor,d)),\ell)\}_{\ell\in L}$; and (\textit{iii}) the asymptotic normality of $n^{-\beta_*(d)}R(n,d)$, for varying $d$.

Throughout the experiments, we consider $r$-mixtures of vMF densities of the form
\begin{align}
    \bx\mapsto g(\bx;\btheta)\defin\sum_{j=1}^r p_j f_{\mathrm{vMF}}(\bx;\bmu_j,\kappa_j), \label{eq:vmfmix}
\end{align}
with locations $\bmu_1,\ldots,\bmu_r\in\Sd$, concentrations $\kappa_1,\ldots,\kappa_r\geq 0$, and proportions $0\leq p_1,\ldots, p_r\leq 1$ such that $\sum_{j=1}^r p_j=1$. The parameter vector $\btheta$ concatenates all mixture parameters, resulting in $r(d+2)-1$ free parameters. We denote by $\btheta_0$ the true parameter used in the data generating process. We also consider the vMF kernel in the kernel density estimator.

\subsection{Von Mises--Fisher distribution}
\label{sec:num:vMF}

We consider a vMF distribution given in \eqref{eq:vmf} with mean direction $\bmu_0=\be_1=(1,0,\stackrel{d}{\ldots},0)^\top$ and concentration $\kappa_0=5$ as the data generating process in this section.

To accurately compute the ratio $R^{(j)}(n,d)$ it is needed to precisely determine both $h_\mathrm{MISE}=\arg\min_{h>0}\allowbreak\mathrm{MISE}\{\hat f(\cdot;h)\}$ and $\hat{h}_\mathrm{CV}=\arg\min_{h>0} \mathrm{CV}(h)$, each of which poses computational challenges. To evaluate $\mathrm{MISE}\{\hat f(\cdot;h)\}$, we used the exact expression given in Proposition 4 of \cite{Garcia-Portugues2013b} for the vMF kernel and a mixture of vMF densities \eqref{eq:vmfmix}:
\begin{align}
    \mathrm{MISE}_{\btheta}\{\hat f(\cdot;h)\}=\frac{c_d^{\mathrm{vMF}}(h^{-2})^2}{nc_d^{\mathrm{vMF}}(2h^{-2})}+\bp^\top\lrc{(1-n^{-1})\boldsymbol{\Psi}_2(h)-2\boldsymbol{\Psi}_1(h)+\boldsymbol{\Psi}_0(h)}\bp, \label{eq:exactmise}
\end{align}
where $\bp\defin(p_1,\ldots,p_r)^\top$ and the $r\times r$ matrices $\boldsymbol{\Psi}_a(h)$, $a=0,1,2$, have $ij$-entries given by
\begin{align*}
    \Psi_{0,ij}(h)&\defin\frac{c_d^{\mathrm{vMF}}(\kappa_i)c_d^{\mathrm{vMF}}(\kappa_j)}{c_d^{\mathrm{vMF}}(\|\kappa_i\bmu_i+\kappa_j\bmu_j\|)},\\
    \Psi_{1,ij}(h)&\defin\int_{\Sd}\frac{c_d^{\mathrm{vMF}}(1/h^2) c_d^{\mathrm{vMF}}(\kappa_i) f_{\mathrm{vMF}}(\by;\bmu_j,\kappa_j)}{c_d^{\mathrm{vMF}}\lrp{\|\by/h^2+\kappa_i\bmu_i\|}}\,\sigmad(\rd\by),\\
    \Psi_{2,ij}(h)&\defin\int_{\Sd}\frac{c_d^{\mathrm{vMF}}(1/h^2)^2 c_d^{\mathrm{vMF}}(\kappa_i)c_d^{\mathrm{vMF}}(\kappa_j)}{c_d^{\mathrm{vMF}}(\|\by/h^2+\kappa_i\bmu_i\|)c_d^{\mathrm{vMF}}(\|\by/h^2+\kappa_j\bmu_j\|)}\,\sigmad(\rd\by).
\end{align*}
Since \eqref{eq:exactmise} does not depend on the sample, its evaluation cost is $O(1)$ in terms of the sample size $n$. The integrals in $\boldsymbol{\Psi}_1(h)$ and $\boldsymbol{\Psi}_2(h)$ can be efficiently computed by importance sampling. First, for the integral in $\Psi_{1,ij}(h)$,
\begin{align}
    \int_{\Sd}\frac{f_{\mathrm{vMF}}(\by;\bmu_j,\kappa_j)}{c_d^{\mathrm{vMF}}\lrp{\|\by/h^2+\kappa_i\bmu_i\|}}\,\sigmad(\rd\by)\approx \frac{1}{B} \sum_{b=1}^B c_d^{\mathrm{vMF}}\lrp{\|\bY_{j,b}/h^2+\kappa_i\bmu_i\|}^{-1},\label{eq:intmc1}
\end{align}
where $\bY_{j,1},\ldots,\bY_{j,B}$ is a random sample from $f_{\mathrm{vMF}}(\cdot;\bmu_j, \kappa_j)$. For $\Psi_{2,ij}(h)$,
\begin{align}
    \int_{\Sd}[c_d^{\mathrm{vMF}}&(\|\by/h^2+\kappa_i\bmu_i\|)c_d^{\mathrm{vMF}}(\|\by/h^2+\kappa_j\bmu_j\|)]^{-1}\,\sigmad(\rd\by)\nonumber\\
    =&\;\int_{\Sd}\frac{f_{ij}(\by)}{f_{ij}(\by)c_d^{\mathrm{vMF}}(\|\by/h^2+\kappa_i\bmu_i\|)c_d^{\mathrm{vMF}}(\|\by/h^2+\kappa_j\bmu_j\|)}\,\sigmad(\rd\by)\nonumber\\
    \approx&\; \frac{1}{B} \sum_{b=1}^B \lrc{f_{ij}(\bY_{ij,b})c_d^{\mathrm{vMF}}(\|\bY_{ij,b}/h^2+\kappa_i\bmu_i\|)c_d^{\mathrm{vMF}}(\|\bY_{ij,b}/h^2+\kappa_j\bmu_j\|)}^{-1},\label{eq:intmc2}
\end{align}
where $\bY_{ij,1},\ldots,\bY_{ij,B}$ is a random sample from the mixture $f_{ij}(\cdot)\defin\frac{1}{2}f_{\mathrm{vMF}}(\cdot;\bmu_i, \kappa_i)+\frac{1}{2}f_{\mathrm{vMF}}(\cdot;\bmu_j, \kappa_j)$. Notice that the samples and computations for rows $i=1,\ldots,r$ in \eqref{eq:intmc1} are reusable for estimating \eqref{eq:intmc2}. We used Monte Carlo samples of size $B=10,\!000$, and we fixed these samples for different bandwidths. Note the evaluation of \eqref{eq:exactmise} is $O(1)$ with respect to $n$.

To compute $\hat{h}_\mathrm{CV}$ we used Proposition 10 in \cite{Garcia-Portugues:polysphere}, which provides the exact cross-validation loss \eqref{eq:cv} for the vMF kernel:
\begin{align}
  \mathrm{CV}(h)=&\;\frac{c_{d}^{\mathrm{vMF}}\big(h^{-2}\big){}^2}{n c_{d}^{\mathrm{vMF}}\big(2 h^{-2}\big)}\nonumber\\
  &-\frac{2 c_{d}^{\mathrm{vMF}}\big(h^{-2}\big)}{n^2} \sum_{1\leq i<j\leq n} \Bigg[\frac{2n}{n-1}e^{\bX_{i}^\top\bX_{j}h^{-2}}-\frac{c_{d}^{\mathrm{vMF}}\big(h^{-2}\big)}{c_{d}^{\mathrm{vMF}}\big([2(1-\bX_{i}^\top\bX_{j})]^{1/2}h^{-2}\big)}\Bigg].\label{eq:exactlscv}
\end{align}
The cross-validation loss cost is $O(n^2)$, but its evaluation can be alleviated by precomputing and reusing $\{\bX_i^\top\bX_j\}_{1\leq i<j\leq n}$ for different bandwidths.

The minimization of $h\mapsto\mathrm{MISE}\{\hat f(\cdot;h)\}$ was done by first running an initial search on the bandwidths $h_0\in\{c\times h_\mathrm{AMISE}:c=0.5,0.55,\ldots,2\}$, and then refining the obtained bandwidth $\tilde{h}_0$ with the lowest MISE by initializing a Newton-type minimization algorithm on $\tilde{h}_0$. We used base R's \texttt{nlm()} with a scaled gradient tolerance of $10^{-10}$ as a stopping criterion. To minimize the more nonlinear function $h\mapsto\mathrm{CV}(h)$, we followed an analogous procedure: performed a search for the initial grid of bandwidths $h_0\in\{0.005,0.030,\ldots,0.555\}$ and then refined the best obtained bandwidth $\hat{h}_0$ with the same Newton-type minimization algorithm. The initial grids of bandwidths were chosen after identifying the feasible ranges of bandwidths for all scenarios $(n,d)$ through simulations.

Minimizing \eqref{eq:exactmise} and \eqref{eq:exactlscv} requires repetitively evaluating the vMF normalizing constant $c_{d}^{\mathrm{vMF}}(x)$ with its modified Bessel function $\mathcal{I}_{(d-1)/2}(x)$, which is both costly and can easily overflow for large argument $x$. To notably speed up its evaluation and perform numerically stable computations, we performed a table evaluation of $x\mapsto\log(e^{-x}\mathcal{I}_{(d-1)/2}(x))$ for a dense grid in $[0,10^4]$, and then carried out a spline interpolation when $x\in[0,10^4]$. For $x>10^4$, we used the asymptotic expansion $\log(e^{-x}\mathcal{I}_{(d-1)/2}(x))\approx \operatorname{log1p}(-d(d-2)/(8x))-\log(2\pi x)/2$.

Figure \ref{fig:avg-rmse} shows the evolution of $\hat{\mathrm{e}}(n,d)$ and $\widehat{\mathrm{rmse}}(n,d)$. The estimated average converges rapidly to zero for larger dimensions, but for $d=1$ it is markedly erratic and seems to converge very slowly over the explored sample sizes. An explanation for this fact lies in the distribution of the bandwidths $\hat{h}_\mathrm{CV}$ being highly left-skewed for small $d$'s, as it can be seen in Figure \ref{fig:dens}. Indeed, the version of the plot where the trimmed mean in $\hat{\mathrm{e}}(n,d)$ is replaced with the median (see Section \ref{sec:add} of the SM) shows that the median curves approach zero from positive values and that there is a decreasing monotone relation between $d$ and the degree of left skewness. The estimated RMSE decreases approximately linearly (in $\log_2$-scale) with respect to the sample size $n$, which is evidenced in the accuracy of the linear fits. The pattern clearly shows that, monotonically, the larger the dimension $d$, the steeper the negative slope capturing the reduction of the RMSE.

Table \ref{tab:linfit:vMF} collects the outcomes of the linear fits $\log_2(\widehat{\mathrm{rmse}}(n,d))\approx\hat{\alpha}(d)+\hat{\beta}(d)\log_2(n)$ for $n\geq 128$ and $d=1,\ldots,10$. In all cases, the coefficient of determination is almost one for all dimensions. The estimated slopes $\hat{\beta}(d)$ show an interesting and somewhat unexpected behavior: they are significantly smaller than the theoretical rate $\beta_*(d)$. Indeed, the difference $\Delta=\beta_*(d)-\hat{\beta}(d)$ is positive and attains $1/10$ when $d=10$, while the relative difference $\Delta\pct=(\beta_*(d)-\hat{\beta}(d))/|\beta_*(d)|\times 100$ exhibits a slightly decaying trend from $d=1$ to $d=10$, indicating that the empirical (log)rate is approximately $27\pct$ faster than the theoretical (log)rate. One-sided $t$-test of $H_0:\beta(d)=\beta_*(d)$ against $H_1:\beta(d)< \beta_*(d)$ corroborates that the differences between $\hat{\beta}(d)$ and $\beta_*(d)$ are highly significant ($p$-values smaller than $7\times 10^{-6}$; omitted). The finding that $\hat{\beta}(d)<\beta_*(d)$ for all the explored dimensions $d$ is compatible with \eqref{eq:rate} and suggests that this rate could be tightened for finite sample sizes.

\begin{figure}[htb!]
  \centering
  \begin{subfigure}{0.5\linewidth}
    \includegraphics[width=\linewidth]{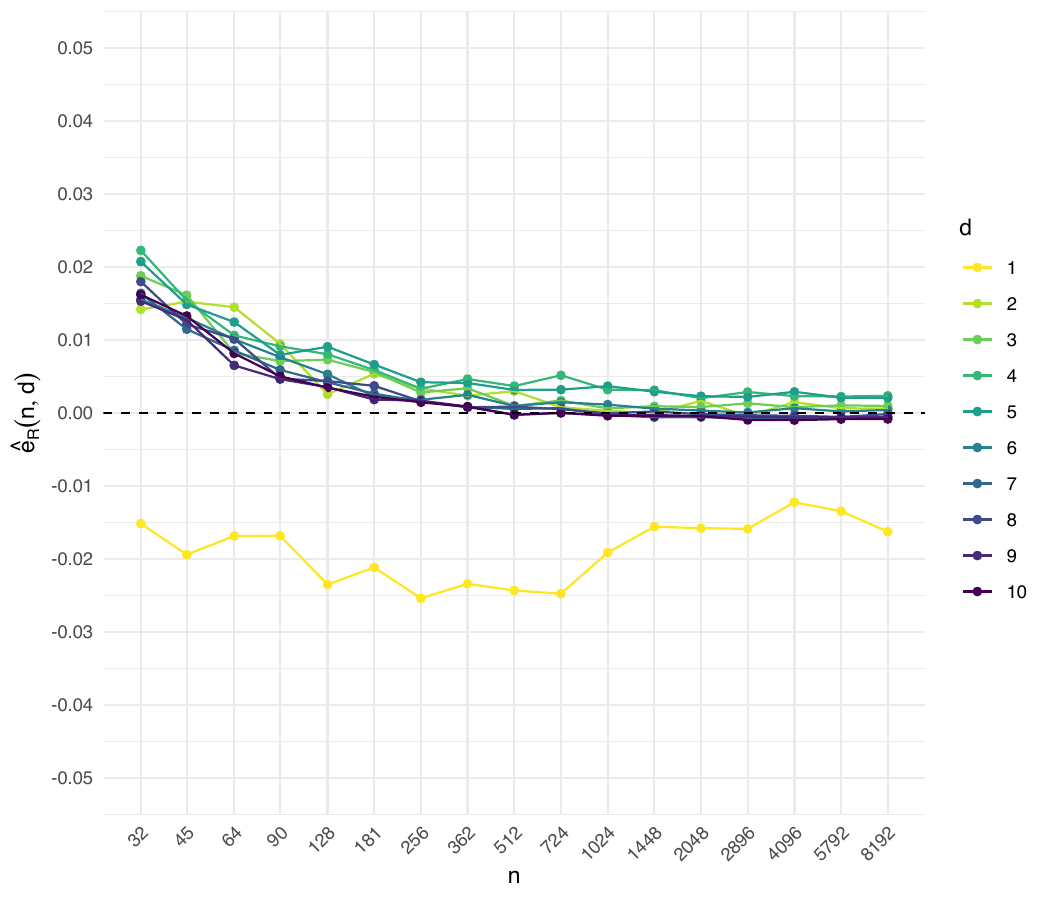}
    \caption{Averages $n\mapsto\hat{\mathrm{e}}(n,d)$.\label{fig:avg-rmse:avg}}
  \end{subfigure}%
  \begin{subfigure}{0.5\linewidth}
    \includegraphics[width=\linewidth]{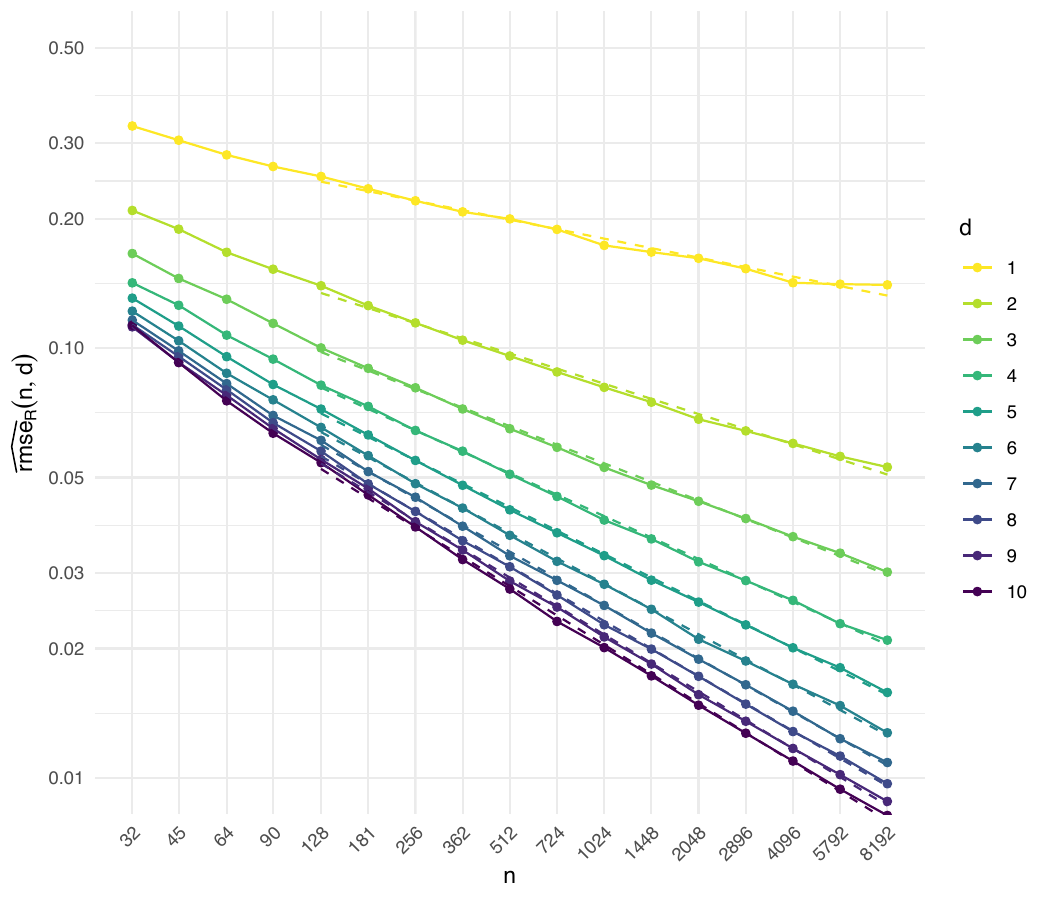}
    \caption{RMSEs $n\mapsto\widehat{\mathrm{rmse}}(n,d)$.\label{fig:avg-rmse:rmse}}
  \end{subfigure}
  \caption{Curves $n\mapsto\hat{\mathrm{e}}(n,d)$ (Figure \ref{fig:avg-rmse:avg}) and $n\mapsto\widehat{\mathrm{rmse}}(n,d)$ (Figure \ref{fig:avg-rmse:rmse}) for dimensions $d=1,2,\ldots,10$. A $\log_2$-scale is used in the horizontal axes of both panels, and also for the vertical axis of the right panel. Linear model fits for $\{(\log_2(\widehat{\mathrm{rmse}}(n,d)),\log_2(n))\}_{n\geq n_0}$ with $n_0=128$ are shown with dashed lines (see Table \ref{tab:linfit:vMF} for their summaries).\label{fig:avg-rmse}}
\end{figure}

\begin{table}[hbt!]
\centering
\small
\setlength{\tabcolsep}{4pt}
\begin{tabular}{llrrrrrrrrrr}
\toprule
$\hat{h}$ & Metric & $d = 1$ & $d = 2$ & $d = 3$ & $d = 4$ & $d = 5$ & $d = 6$ & $d = 7$ & $d = 8$ & $d = 9$ & $d = 10$\\
\midrule
CV & $\beta_*(d)$ & $-0.10$ & $-0.17$ & $-0.21$ & $-0.25$ & $-0.28$ & $-0.30$ & $-0.32$ & $-0.33$ & $-0.35$ & $-0.36$\\
& $\hat{\beta}(d)$ & $-0.15$ & $-0.23$ & $-0.29$ & $-0.33$ & $-0.36$ & $-0.39$ & $-0.41$ & $-0.42$ & $-0.44$ & $-0.45$\\
& $\Delta\pct$ & $47\pct$ & $40\pct$ & $34\pct$ & $33\pct$ & $30\pct$ & $30\pct$ & $30\pct$ & $27\pct$ & $27\pct$ & $27\pct$\\
& $R^2(d)$ & $0.98$ & $1.00$ & $1.00$ & $1.00$ & $1.00$ & $1.00$ & $1.00$ & $1.00$ & $1.00$ & $1.00$\\\midrule
AMI & $\beta_*(d)$ & $-0.40$ & $-0.33$ & $-0.29$ & $-0.25$ & $-0.22$ & $-0.20$ & $-0.18$ & $-0.17$ & $-0.15$ & $-0.14$\\
& $\hat{\beta}(d)$ & $-0.44$ & $-0.38$ & $-0.31$ & $-0.28$ & $-0.25$ & $-0.22$ & $-0.21$ & $-0.19$ & $-0.18$ & $-0.17$\\
& $\Delta\pct$ & $10\pct$ & $13\pct$ & $10\pct$ & $12\pct$ & $13\pct$ & $12\pct$ & $13\pct$ & $16\pct$ & $19\pct$ & $22\pct$\\
& $R^2(d)$ & $1.00$ & $1.00$ & $1.00$ & $1.00$ & $1.00$ & $1.00$ & $1.00$ & $1.00$ & $1.00$ & $1.00$\\\midrule
EMI & $\beta_*(d)$ & $-0.50$ & $-0.50$ & $-0.50$ & $-0.50$ & $-0.50$ & $-0.50$ & $-0.50$ & $-0.50$ & $-0.50$ & $-0.50$\\
& $\hat{\beta}(d)$ & $-0.50$ & $-0.50$ & $-0.50$ & $-0.50$ & $-0.50$ & $-0.50$ & $-0.51$ & $-0.52$ & $-0.53$ & $-0.53$\\
& $\Delta\pct$ & $0\pct$ & $0\pct$ & $0\pct$ & $1\pct$ & $0\pct$ & $0\pct$ & $2\pct$ & $4\pct$ & $5\pct$ & $6\pct$\\
& $R^2(d)$ & $1.00$ & $1.00$ & $1.00$ & $1.00$ & $1.00$ & $1.00$ & $1.00$ & $1.00$ & $1.00$ & $1.00$\\
\bottomrule
\end{tabular}
\caption{Summaries of the linear model fits $\log_2(\widehat{\mathrm{rmse}}(n,d))\approx\hat{\alpha}(d)+\hat{\beta}(d)\log_2(n)$ for $n\geq n_0=128$. The theoretical rates are $\beta_{*,\mathrm{CV}}(d)=-d/(2d+8)$, $\beta_{*,\mathrm{AMI}}(d)=-4 / (2 d + 8)$, and $\beta_{*,\mathrm{EMI}}(d)=-1/2$. The relative difference is $\Delta\pct=(\beta_*(d)-\hat{\beta}(d))/|\beta_*(d)|\times 100$ and $R^2(d)$ stands for the coefficient of determination of the linear fit.\label{tab:linfit:vMF}}
\end{table}

Finally, we evaluate the convergence of $n^{-\beta_*(d)}R(n,d)$ toward a normal distribution in Figure \ref{fig:dens}, showing the kernel density estimates computed for the Monte Carlo samples $\{n^{-\beta_*(d)}R^{(j)}(n,d)\}_{j=1}^M$, for varying $n$ and $d$. For small dimensions, especially $d=1,2$, the densities are highly left-skewed, even for large sample sizes, indicating persistent non-normality. However, for larger dimensions, the densities exhibit a rapid convergence to a normal-like shape, which is also confirmed by the normal QQ-plots in Section \ref{sec:add} of the SM. The $p$-values of the Lilliefors test of normality in Section \ref{sec:add} confirm this trend: the null hypothesis of normality is strongly rejected for dimensions $d\geq 4$ in all the explored sample sizes. However, as the dimension increases, the $p$-values across different sample sizes $n$ become increasingly uniform. The $p$-values of the Lilliefors test of normality in Section \ref{sec:add} confirm this trend: the null hypothesis of normality is strongly rejected for dimensions $d\geq 4$ in all the explored sample sizes. However, as the dimension increases, the $p$-values across different sample sizes $n$ become increasingly uniform. For example, for dimensions $d=7,8,9,10$, non-rejections of normality at the $5\pct$ significance level begin to occur for sample sizes $n$ as small as $64$ and~$128$.

\begin{figure}[htb!]
  \centering
  \includegraphics[width=\linewidth]{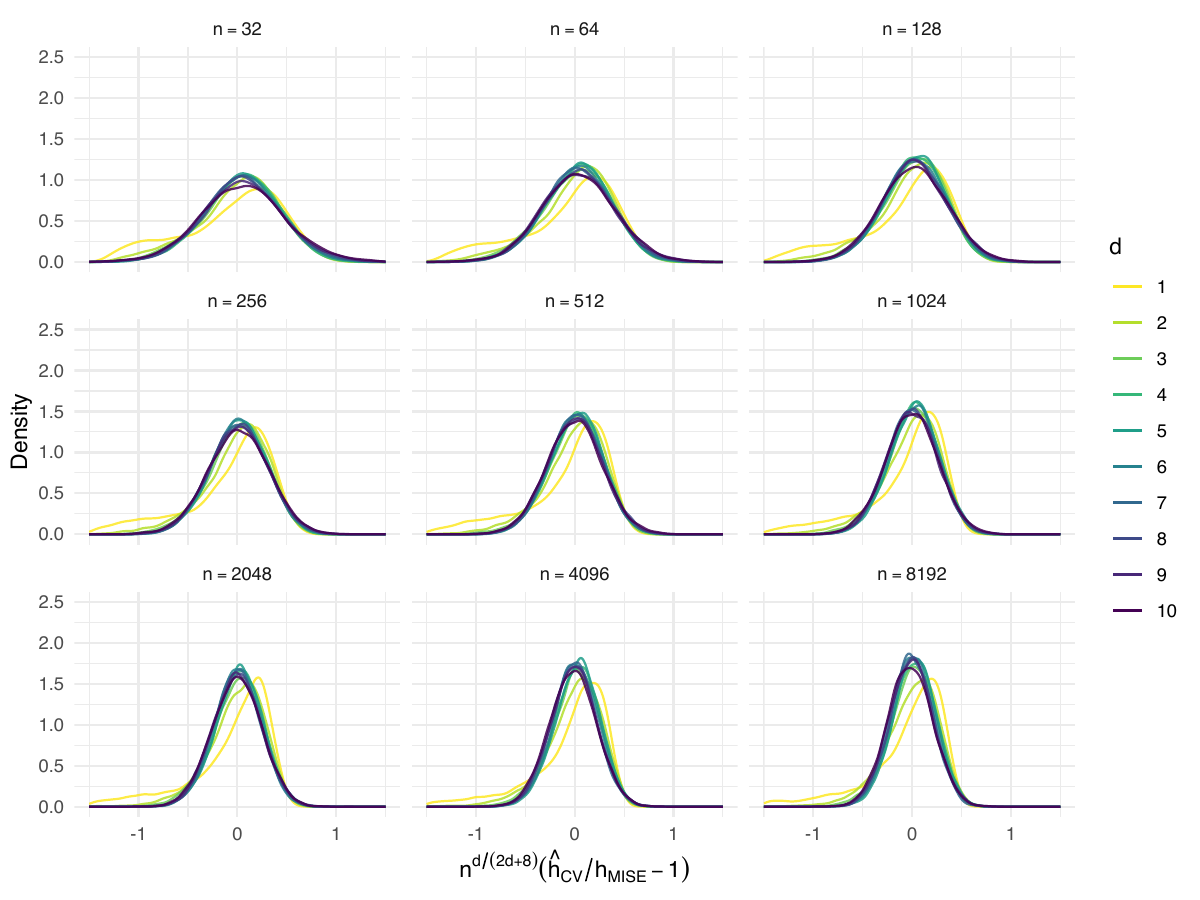}
  \caption{Evaluation of the asymptotic normality of $n^{-\beta_*(d)}R(n,d)$ for sample sizes $n=2^\ell$, $\ell=5,5.5,\ldots,13$ and dimensions $d=1,2,\ldots,10$. The curves show the kernel density estimates for the sample $\{n^{-\beta_*(d)}R^{(j)}(n,d)\}_{j=1}^M$ featuring normal scale bandwidths.\label{fig:dens}}
\end{figure}


\subsection{Mixtures of von Mises--Fisher distributions}
\label{sec:num:MvMF}

We assume in this section that the data generating process is a Mixture of four vMF densities (MvMF) with antipodal mean directions $\bmu_{0,1}=\be_1$, $\bmu_{0,2}=-\be_1$, $\bmu_{0,3}=\be_{d+1}$, and $\bmu_{0,4}=-\be_{d+1}$, with common concentrations $\kappa_{0,j} = 5$ and proportions $p_{0,j}=1/4$, for $j=1,\ldots,4$.

In contrast with the vMF model in Section \ref{sec:num:vMF}, the MvMF distribution is challenging for any bandwidth selector: the uniform density estimate obtained with $h=\infty$ is a competitive choice depending on the sample size and dimension. The ill-definedness of $h_{\mathrm{MISE}}$ for small-to-moderate sample sizes is depicted in Figure \ref{fig:rmse-mise:mise}: (\textit{i}) for $(n=32,64,d=1)$, the MISE curves seem to have a global minima located at $h=\infty$; (\textit{ii}) for $(n=64,d=1)$, a local finite minimum appears at $h\approx 0.40$, but the global minimum seems to be still located at $h=\infty$; (\textit{iii}) the global minimum below $h=0.40$ becomes better localized for $n\geq 256$; (\textit{iv}) for $d=10$, the gap between the global minima for $n=32,64$ and the error for $h=\infty$ narrows with respect to $d=1$. Due to the ill-definedness of $h_{\mathrm{MISE}}$, the curves $h\mapsto \mathrm{MISE}\{\hat{f}(\cdot;h)\}$ that $\mathrm{CV}(h)$ estimates and minimizes show highly variable patterns, with ill-localized global minima, unless $n$ is large.

\begin{figure}[htb!]
  \centering
  \begin{subfigure}{0.5\linewidth}
    \includegraphics[width=\linewidth]{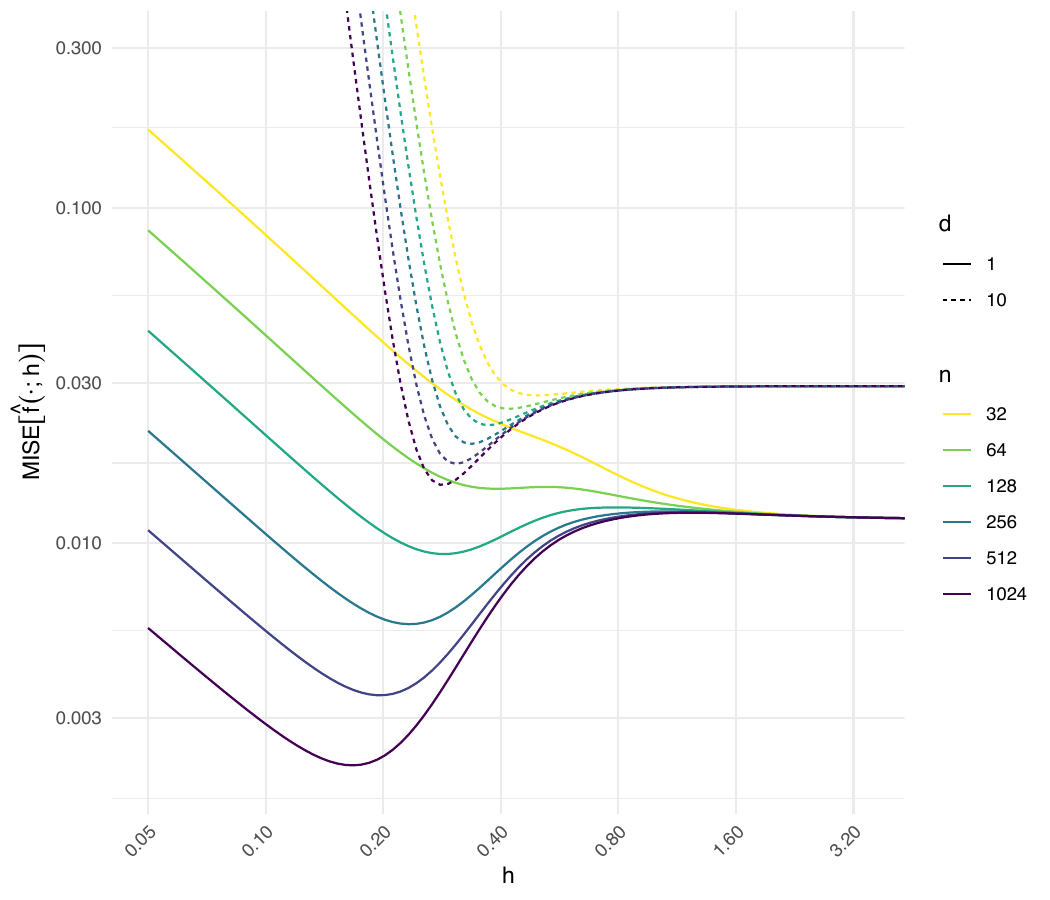}
    \caption{Errors $h\mapsto \mathrm{MISE}[\hat{f}(\cdot;h)]$.\label{fig:rmse-mise:mise}}
  \end{subfigure}%
  \begin{subfigure}{0.5\linewidth}
    \includegraphics[width=\linewidth]{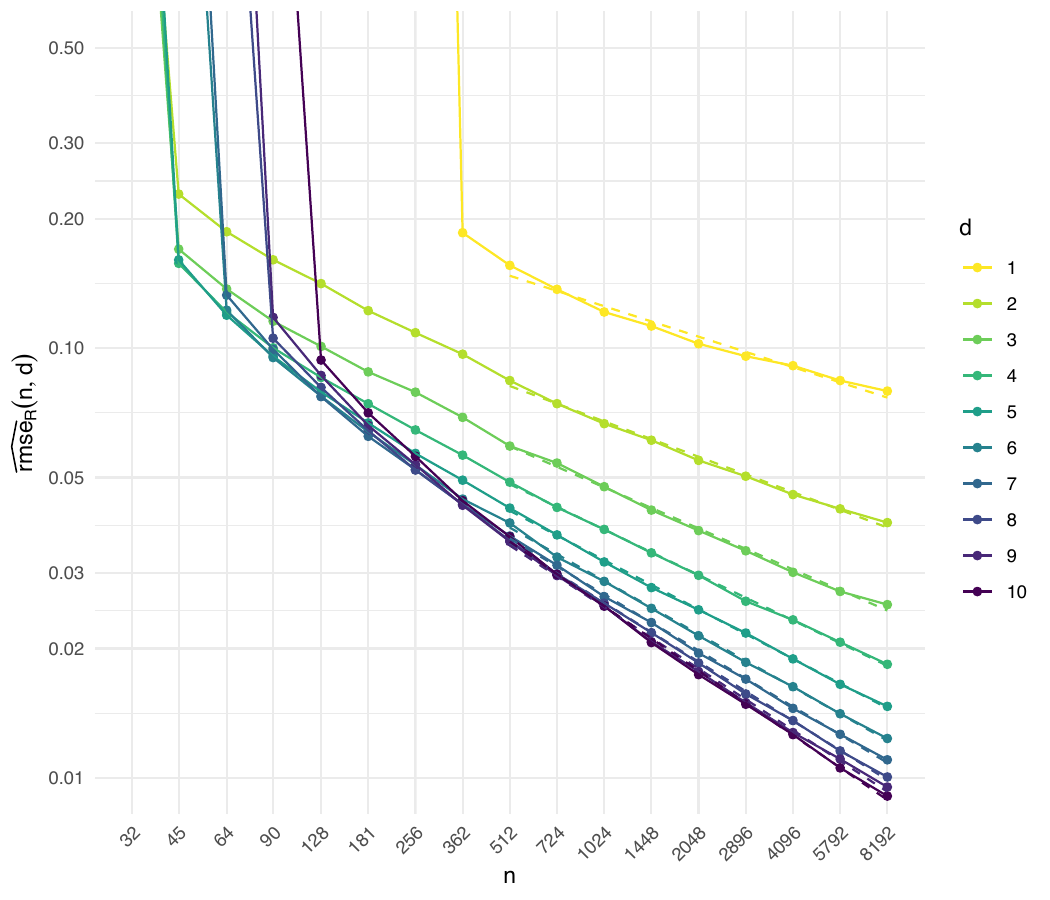}
    \caption{RMSEs $n\mapsto\widehat{\mathrm{rmse}}(n,d)$.\label{fig:rmse-mise:rmse}}
  \end{subfigure}%
  \caption{In Figure \ref{fig:rmse-mise:mise}, curves $h\mapsto \mathrm{MISE}\{\hat{f}(\cdot;h)\}$ for dimensions $d=1,10$ and sample sizes $n=32,64,\ldots,1024$ showing the varying degree of identifiability of the global minima in the MvMF distribution. In Figure \ref{fig:rmse-mise:rmse}, curves $n\mapsto\widehat{\mathrm{rmse}}(n,d)$ for dimensions $d=1,2,\ldots,10$. A $\log_2$-scale is used in both axes. Linear model fits for $\{(\log_2(\widehat{\mathrm{rmse}}(n,d)),\log_2(n))\}_{n\geq n_0}$ with $n_0=512$ are shown with dashed lines (see Table \ref{tab:linfit:MvMF} for their summaries).\label{fig:rmse-mise}}
\end{figure}

Figure \ref{fig:rmse-mise:rmse} shows the evolution of $\widehat{\mathrm{rmse}}(n,d)$. For large sample sizes, depending on the dimension, it decreases linearly (in $\log_2$-scale) in the sample size $n$, as seen in the linear fits. As before, the pattern clearly shows that, monotonically, the larger the dimension $d$, the steeper the negative slope capturing the reduction of the RMSE. The linear stabilization of the trends is inversely monotonic on $d$, except for $d=1$: for dimensions $d=2,3$ the stabilization starts happening already from $n=64$, while for $d=9,10$ it takes a larger sample size, and for $d=1$ the RMSEs are only within a comparable range with other dimensions for $n=362$. These phenomena are a direct consequence of the relative gaps between $\mathrm{MISE}\{\hat{f}(\cdot;h_{\mathrm{MISE}})\}$ and $\mathrm{MISE}\{\hat{f}(\cdot;\infty)\}$, which are maximized for $d=2$ and decrease monotonically until $d=10$ (uniformly in $n$), and are minimized for $d=1$ (for small-to-moderate sample sizes).

Table \ref{tab:linfit:MvMF} collects the outcomes of the linear fits $\log_2(\widehat{\mathrm{rmse}}(n,d))\approx\hat{\alpha}(d)+\hat{\beta}(d)\log_2(n)$ done for $n\geq 512$. In all of them, the coefficient of determination is almost one for all dimensions. As for Table \ref{tab:linfit:vMF}, the estimated slopes $\hat{\beta}(d)$ are significantly smaller than the theoretical rate $\beta_*(d)$ ($p$-values smaller than $2\times 10^{-6}$; omitted). The relative difference attains its maximum at $d=1$ and $d=2$, and then exhibits a decaying trend until $d=6$, with an increase toward $d=10$. This latter behavior might be a consequence of the slower stabilization of the RMSEs for larger dimensions. The results corroborate that, although $\hat{h}_\mathrm{CV}$ passes through a very unstable phase with potentially infinite variances, $R(n,d)$ eventually converges at least as fast as the prescribed rate as $n$ diverges to infinity.

\begin{table}[htb!]
\centering
\small
\setlength{\tabcolsep}{4pt}
\begin{tabular}{llrrrrrrrrrr}
\toprule
$\hat{h}$ & Metric & $d = 1$ & $d = 2$ & $d = 3$ & $d = 4$ & $d = 5$ & $d = 6$ & $d = 7$ & $d = 8$ & $d = 9$ & $d = 10$\\
\midrule
CV & $\beta_*(d)$ & $-0.10$ & $-0.17$ & $-0.21$ & $-0.25$ & $-0.28$ & $-0.30$ & $-0.32$ & $-0.33$ & $-0.35$ & $-0.36$\\
& $\hat{\beta}(d)$ & $-0.24$ & $-0.27$ & $-0.32$ & $-0.35$ & $-0.38$ & $-0.41$ & $-0.43$ & $-0.46$ & $-0.48$ & $-0.50$\\
& $\Delta\pct$ & $136\pct$ & $64\pct$ & $48\pct$ & $40\pct$ & $37\pct$ & $37\pct$ & $36\pct$ & $37\pct$ & $37\pct$ & $40\pct$\\
& $R^2(d)$ & $0.98$ & $1.00$ & $1.00$ & $1.00$ & $1.00$ & $1.00$ & $1.00$ & $1.00$ & $1.00$ & $1.00$\\\midrule
AMI & $\beta_*(d)$ & $-0.40$ & $-0.33$ & $-0.29$ & $-0.25$ & $-0.22$ & $-0.20$ & $-0.18$ & $-0.17$ & $-0.15$ & $-0.14$\\
& $\hat{\beta}(d)$ & $-0.55$ & $-0.45$ & $-0.38$ & $-0.33$ & $-0.29$ & $-0.27$ & $-0.25$ & $-0.24$ & $-0.23$ & $-0.23$\\
& $\Delta\pct$ & $38\pct$ & $34\pct$ & $31\pct$ & $31\pct$ & $32\pct$ & $35\pct$ & $39\pct$ & $45\pct$ & $52\pct$ & $59\pct$\\
& $R^2(d)$ & $1.00$ & $1.00$ & $1.00$ & $1.00$ & $1.00$ & $1.00$ & $1.00$ & $0.99$ & $0.99$ & $0.99$\\\midrule
EMI & $\beta_*(d)$ & $-0.50$ & $-0.50$ & $-0.50$ & $-0.50$ & $-0.50$ & $-0.50$ & $-0.50$ & $-0.50$ & $-0.50$ & $-0.50$\\
& $\hat{\beta}(d)$ & $-0.56$ & $-0.59$ & $-0.61$ & $-0.62$ & $-0.64$ & $-0.67$ & $-0.70$ & $-0.72$ & $-0.75$ & $-0.79$\\
& $\Delta\pct$ & $12\pct$ & $17\pct$ & $21\pct$ & $24\pct$ & $28\pct$ & $34\pct$ & $40\pct$ & $44\pct$ & $51\pct$ & $59\pct$\\
& $R^2(d)$ & $1.00$ & $1.00$ & $1.00$ & $1.00$ & $1.00$ & $1.00$ & $1.00$ & $1.00$ & $0.99$ & $0.99$\\
\bottomrule
\end{tabular}
\caption{Same table as Table \ref{tab:linfit:vMF}, but for the MvMF data-generating process and $n_0=512$. \label{tab:linfit:MvMF}}
\end{table}

\subsection{Comparison with plug-in bandwidths}
\label{sec:num:PI}

We compare now the convergence rates of $\hat{h}_\mathrm{CV}$ with plug-in bandwidth selectors that target the minimization of the AMISE and MISE. We do so under scenarios that benefit plug-in selectors to illustrate the striking degree of competitiveness that $\hat{h}_\mathrm{CV}$ achieves as the dimension grows.

The Asymptotic MIxtures (AMI) and Exact MIxtures (EMI) in \cite{Garcia-Portugues2013a} are plug-in bandwidth selectors that exploit mixtures of vMF densities \eqref{eq:vmfmix}. The AMI selector follows the traditional approach to plug-in bandwidth selection that replaces the curvature term $R\big(\nabla^2 \bar{f}\big)$ in the $h_\mathrm{AMISE}$ bandwidth \eqref{eq:hamise} with $R\big(\nabla^2 \bar{g}(\cdot;\hat{\btheta})\big)$, resulting in $\hat{h}_\mathrm{AMI}\defin\arg\min_{h>0}\mathrm{AMISE}_{\hat{\btheta}}\{\hat f(\cdot;h)\}$. The EMI selector uses a similar approach, but leveraging \eqref{eq:exactmise} to use the estimated non-asymptotic error: $\hat{h}_\mathrm{EMI}\defin\arg\min_{h>0}\mathrm{MISE}_{\hat{\btheta}}\{\hat f(\cdot;h)\}$. To estimate the mixtures, the Expectation--Maximization (EM) algorithm, as implemented in the R package \texttt{movMF} \citep{Hornik2014}, can be used. In the numerical experiments, we used \texttt{movMF()} from the latter package with five runs of the EM with $100$ maximum iterations and the default initialization. In the case of a single vMF distribution, maximum likelihood is direct and benefits from the discussion on the evaluation of the normalizing constant from Section \ref{sec:num:vMF}.

Assuming the data generating process is truly $g(\cdot;\btheta_0)$, then $\hat{h}_\mathrm{AMI}$ and $\hat{h}_\mathrm{EMI}$ are parametric estimates of $h_\mathrm{AMISE}$ and $h_\mathrm{MISE}$. In this case, if $\hat{\btheta}-\btheta_0=O_\mathbb{P}(n^{-1/2})$, then $R_\mathrm{EMI}(n,d)\defin(\hat{h}_\mathrm{EMI}-h_\mathrm{MISE})/h_\mathrm{MISE}=O_\mathbb{P}\lrp{n^{\beta_{*,\mathrm{EMI}}(d)}}$, with $\beta_{*,\mathrm{EMI}}(d)\defin -1/2$. However, even in this best-case scenario, the rate for $R_\mathrm{AMI}(n,d)\defin(\hat{h}_\mathrm{AMI}-h_\mathrm{MISE})/h_\mathrm{MISE}=O_\mathbb{P}\lrp{n^{\beta_{*,\mathrm{AMI}}(d)}}$ is not parametric, but rather $\beta_{*,\mathrm{AMI}}(d)\defin -4/(2d+8)$, as induced by the rate in \eqref{eq:Ohamise}.

Tables \ref{tab:linfit:vMF}--\ref{tab:linfit:MvMF} collect the outcomes of the linear fits $\log_2(\widehat{\mathrm{rmse}}(n,d))\approx\hat{\alpha}(d)+\hat{\beta}(d)\log_2(n)$ for the estimated RMSEs of $R_\mathrm{AMI}(n,d)$ and $R_\mathrm{EMI}(n,d)$. For both AMI and EMI, the coefficients of determination of the linear fits are almost one. As expected, the estimated slopes $\hat{\beta}_{\mathrm{EMI}}(d)$ are close to or faster than $-1/2$, across all the explored dimensions, evidencing the parametric rate. The estimated slopes for AMI follow a decaying trend in absolute value, paralleling the theoretical rates, and becoming surpassed by the slopes of CV for $d\geq 5$. Consistent with the situation for CV, the empirical rates for AMI and EMI are also faster than the theoretical ones, especially in the MvMF distribution.

\subsection{Density error rates}
\label{sec:num:dens}

We explore in this final section the convergence of the $L^2$-errors $\|\hat{f}(\cdot;\hat{h})-g(\cdot;\btheta_0)\|_2=[\mathrm{ISE}_{\btheta_0}\{\hat{f}(\cdot;\hat{h})\}]^{1/2}$ for the bandwidth selectors $\hat{h}_\mathrm{CV}$, $ \hat{h}_\mathrm{AMI}$, and $\hat{h}_\mathrm{EMI}$ under the vMF and MvMF densities $g(\cdot;\btheta_0)$ from Sections \ref{sec:num:vMF}--\ref{sec:num:MvMF}. As a benchmark, we compute the oracle errors $\|\hat{f}(\cdot;\hat{h}_\mathrm{ISE})-g(\cdot;\btheta_0)\|_2$ with $\hat{h}_\mathrm{ISE}=\min_{h>0}\mathrm{ISE}\{\hat f(\cdot;h)\}$ and the parametric errors $\|g(\cdot;\hat{\btheta})-g(\cdot;\btheta_0)\|_2$. Note that, within this setup, $\hat{f}(\cdot;\hat{h}_\mathrm{CV})$ is the only density estimator that is fully agnostic to the underlying density, as the parametric fit and the plug-in selectors rely on the parametric specification $g(\cdot;\btheta)$, while $\hat{h}_\mathrm{ISE}$ directly uses $g(\cdot;\btheta_0)$.

We explored two possibilities to compute the $\mathrm{ISE}\{\hat f(\cdot;h)\}$ for the vMF kernel and mixtures of vMF densities. The exact expression can be written as
\begin{align}
    \mathrm{ISE}_{\btheta}\{\hat f(\cdot;h)\}
    =&\;\int_{\Sd} \hat{f}(\bx;h)^2 \,\sigmad(\rd\bx)+\int_{\Sd} g(\bx;\btheta)^2 \,\sigmad(\rd\bx)-2\int_{\Sd} \hat{f}(\bx;h) g(\bx;\btheta) \,\sigmad(\rd\bx)\nonumber\\
    =&\;\frac{1}{n^2}\sum_{i,j=1}^n \Psi_0(h^{-2},h^{-2},\bX_i,\bX_j) + \sum_{i,j=1}^r p_ip_j\Psi_0(\kappa_i,\kappa_j,\bmu_i,\bmu_j)\nonumber\\
    &- \frac{2}{n}\sum_{i=1}^n\sum_{j=1}^r p_j\Psi_0(h^{-2},\kappa_j,\bX_i,\bmu_j),\label{eq:isevmf}
\end{align}
where 
\begin{align*}
    \Psi_{0}(\kappa_1,\kappa_2,\mu_1,\mu_2)\defin&\;\int_{\Sd} f_{\mathrm{vMF}}(\bx;\bmu_1,\kappa_1)f_{\mathrm{vMF}}(\bx;\bmu_2,\kappa_2)\,\sigmad(\rd\bx)
    \frac{c_d^{\mathrm{vMF}}(\kappa_1)c_d^{\mathrm{vMF}}(\kappa_2)}{c_d^{\mathrm{vMF}}(\|\kappa_1\bmu_1+\kappa_2\bmu_2\|)}.
\end{align*}
While exact, \eqref{eq:isevmf} scales poorly on the sample size and adds a substantial computational overhead. Therefore, we used importance sampling instead to be able to handle large sample sizes:
\begin{align*}
    \mathrm{ISE}_{\btheta}\{\hat f(\cdot;h)\}= \int_{\Sd} \frac{(\hat{f}(\bx;h)-g(\bx;\btheta))^2}{g(\bx;\btheta)}g(\bx;\btheta)\,\sigmad(\rd\bx) \approx \frac{1}{B}\sum_{b=1}^B \frac{(\hat{f}(\bY_b;h)-g(\bY_b;\btheta))^2}{g(\bY_b;\btheta)},
\end{align*}
where $\bY_1,\ldots,\bY_B$ is a random sample from $g(\cdot;\btheta)$. This approximation, with $B=10,\!000$, was used to compute all the $L^2$-errors.

\begin{figure}[t]
  \centering
  \includegraphics[width=0.9\linewidth]{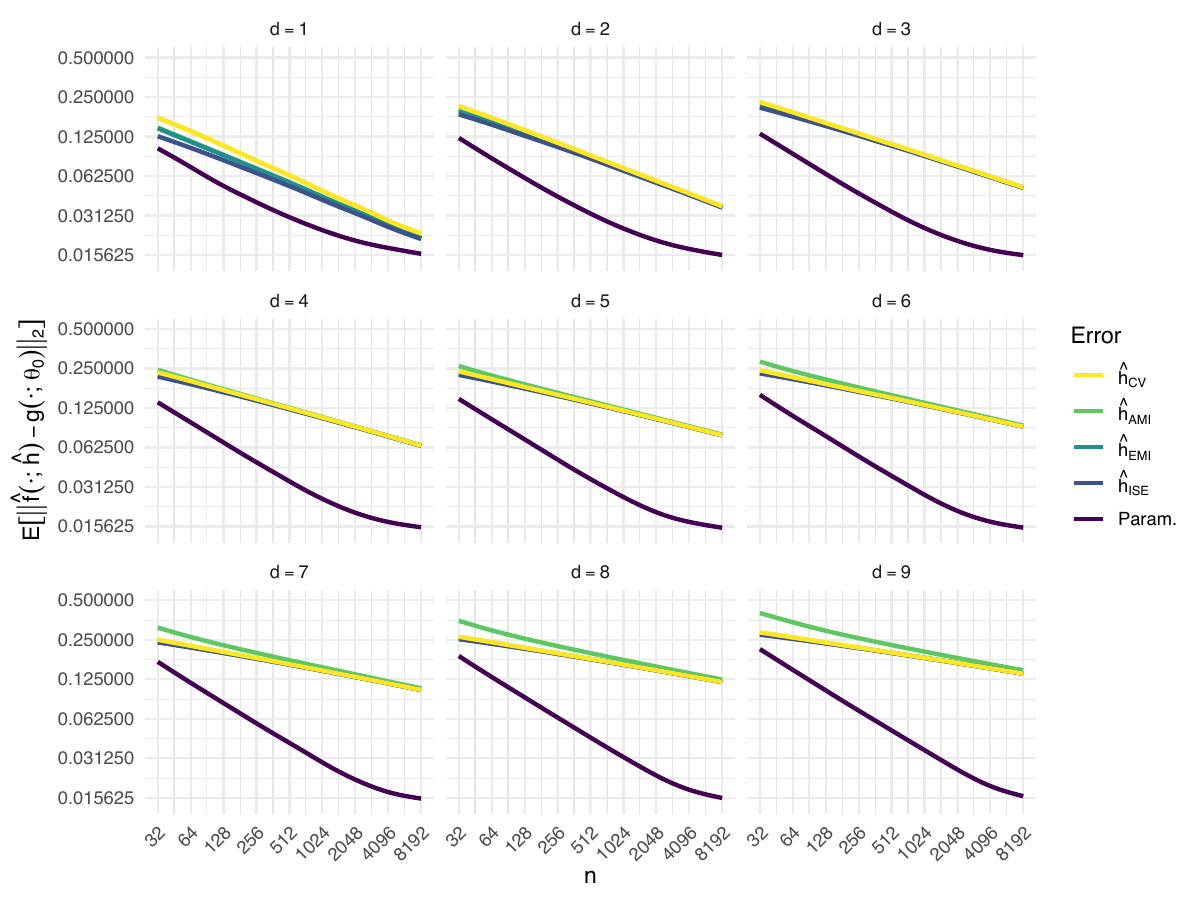}
  \caption{Error curves $n\mapsto \mathbb{E}\big[\|\hat{f}(\cdot;\hat{h})-g(\cdot;\btheta_0)\|_2\big]$ for the bandwidth selectors $\hat{h}\in\{\hat{h}_\mathrm{CV},\hat{h}_\mathrm{AMI},\hat{h}_\mathrm{EMI},\hat{h}_\mathrm{ISE}\}$. These curves are benchmarked with respect to the parametric error curve $n\mapsto \mathbb{E}\big[\|g(\cdot;\hat{\btheta})-g(\cdot;\btheta_0)\|_2\big]$. The density $g(\cdot;\btheta_0)$ is the vMF from Section \ref{sec:num:vMF} and $\hat{\btheta}$ is the maximum likelihood estimator.\label{fig:ise:vMF}}
\end{figure}

The averages of $L^2$-errors are reported in Figures \ref{fig:ise:vMF}--\ref{fig:ise:MvMF} as the functions $n\mapsto \mathbb{E}\big[\|\hat{f}(\cdot;\hat{h})-g(\cdot;\btheta_0)\|_2\big]$, $\hat{h}\in\{\hat{h}_\mathrm{CV},\hat{h}_\mathrm{AMI},\hat{h}_\mathrm{EMI},\hat{h}_\mathrm{ISE}\}$, and $n\mapsto \mathbb{E}\big[\|g(\cdot;\hat{\btheta})-g(\cdot;\btheta_0)\|_2\big]$. In both figures, it is clearly seen how the density errors for $\hat{h}_\mathrm{CV}$ converge to those of $\hat{h}_\mathrm{ISE}$ as $d$ grows, coinciding also with the errors of $\hat{h}_\mathrm{EMI}$. The asymptotic errors of $\hat{h}_\mathrm{AMI}$ are smaller than those of $\hat{h}_{\mathrm{CV}}$ for $d\leq 3$, equal to for $d=4$, and slower for $d>4$. The parametric estimate has the asymptotically smallest errors, yet for MvMF (Figure \ref{fig:ise:MvMF}), the errors are larger than nonparametric approaches for small-to-moderate sample sizes. This behavior is attributed to the challenging density form and the use of the EM algorithm (despite fixing the number of mixture components to four). The difficulty in parametric estimation also translates to the errors of $\hat{h}_\mathrm{AMI}$ and $\hat{h}_\mathrm{EMI}$. However, the speed of convergence of the errors of $\hat{h}_\mathrm{CV}$ to those of $\hat{h}_\mathrm{ISE}$ remains unaltered, evidencing the competitiveness of $\hat{f}(\cdot;\hat{h}_{\mathrm{CV}})$ as the only assumption-free density estimator.

\begin{figure}[htb!]
  \centering
  \includegraphics[width=\linewidth]{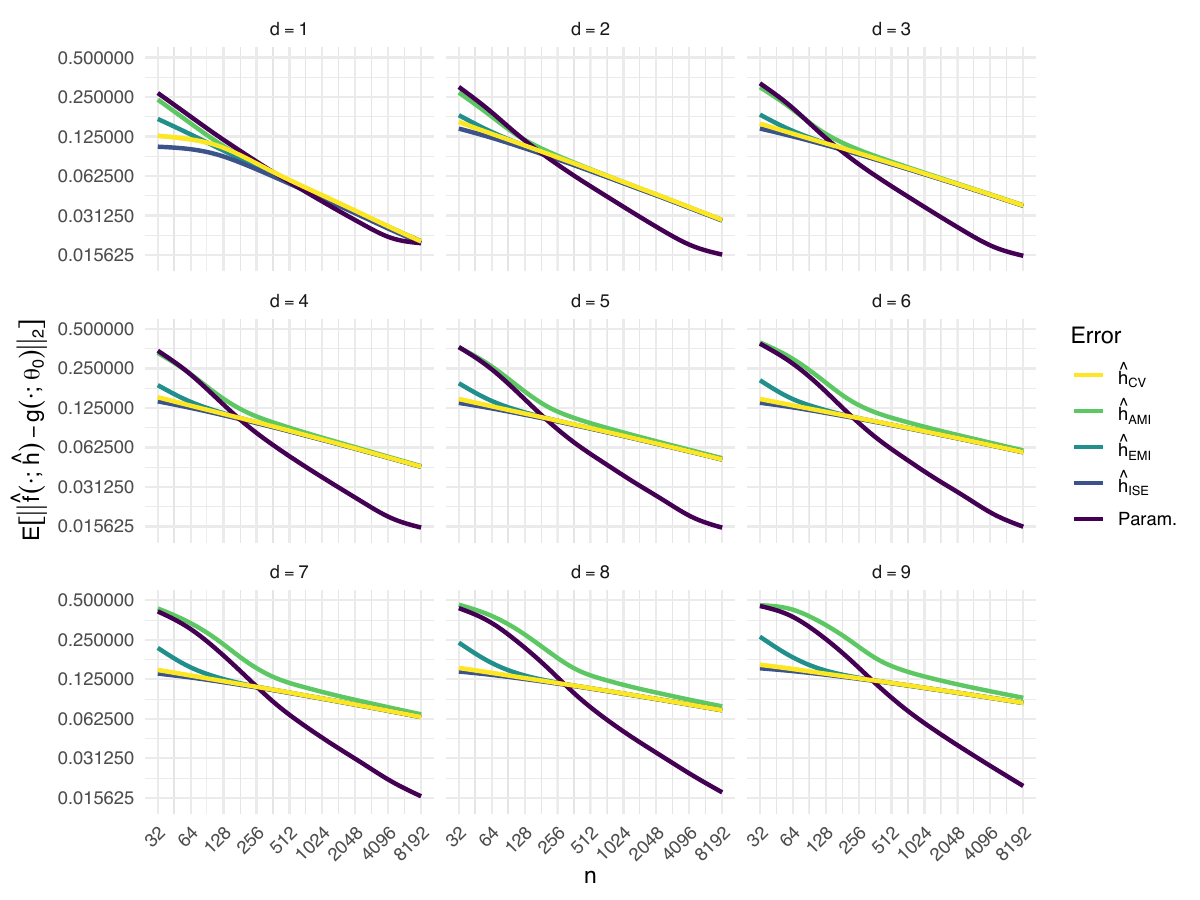}
  \caption{Same figure as Figure \ref{fig:ise:vMF} but with $g(\cdot;\btheta_0)$ the MvMF from Section \ref{sec:num:MvMF} and $\hat{\btheta}$ the EM estimator.\label{fig:ise:MvMF}}
\end{figure}

\section{Discussion}
\label{sec:disc}

In this work, we develop the asymptotic theory for least-squares cross-validated (CV) bandwidth selection in kernel density estimation on the hypersphere $\Sd$. Under mild non-uniformity conditions, we establish the existence of the MISE-optimal bandwidth and derive the exact relative rate of convergence of the CV selector, which is $n^{-d/(2d+8)}$. This rate reveals a genuine blessing of dimensionality, approaching the $n^{-1/2}$ parametric benchmark as $d$ grows. We also provide explicit expressions for the asymptotic variance of the CV bandwidth, allowing direct comparison with the Euclidean setting and highlighting how high-dimensional structures enhance the stability and reliability of cross-validation.

These theoretical results clarify that, in contrast to kernel density estimation itself, which deteriorates with increasing dimension, the selection of the smoothing parameter becomes more accurate as the dimension grows. Extensive numerical experiments further corroborate these conclusions. The performance of the CV selector was evaluated across a range of directional distributions, encompassing both von Mises--Fisher and mixtures of von Mises--Fisher densities, in order to capture a variety of distributional complexities. The results confirm that the relative convergence rate improves with dimension and that CV bandwidths yield density estimates with errors that closely mirror those of oracle and parametric benchmarks, even under challenging mixture settings.

Overall, the combination of theoretical derivations and numerical evidence establishes cross-validation as a robust and practically effective method for bandwidth selection in high-dimensional directional data. This approach yields accurate density estimates without reliance on parametric assumptions, and its performance improves in higher dimensions, providing a clear illustration of the blessing of dimensionality in the context of bandwidth selection.

\section*{Acknowledgments}
 The authors are supported by grant PID2021-124051NB-I00, funded by MCIN/\-AEI/\-10.13039/\-501100011033 and by ``ERDF A way of making Europe''. The authors acknowledge the computational resources from the Centro de Supercomputación de Galicia (CESGA).

\bibliographystyle{apalike}

\newpage

\title{Supplementary material for ``Blessing of dimensionality in cross-validated bandwidth selection on the sphere''}
\setlength{\droptitle}{-1cm}
\predate{}%
\postdate{}%
\date{}
\author{José E. Chacón$^{1,3}$, Eduardo Garc\'ia-Portugu\'es$^{2}$ and Andrea Meil\'an-Vila$^{2}$}
\footnotetext[1]{Department of Mathematics, Universidad de Extremadura (Spain).}
\footnotetext[2]{Department of Statistics, Universidad Carlos III de Madrid (Spain).}
\footnotetext[3]{Corresponding author. e-mail: \href{mailto:jechacon@unex.es}{jechacon@unex.es}.}
\maketitle

\begin{abstract}
    The supplementary material consists of five parts. Section \ref{sec:proofs} contains the proofs of the main results of the paper. Section \ref{sec:aux} provides auxiliary results along with their corresponding proofs. Section \ref{sec:euclidean} presents the proofs of the results from Section \ref{sec:comp} regarding the cross-validation bandwidth selector in the multivariate Euclidean case. Section \ref{sec:calc} provides the calculations for the von Mises--Fisher kernel, while Section \ref{sec:add} presents additional numerical experiments.
\end{abstract}

\appendix

\section{Proofs of the main results}
\label{sec:proofs}

\begin{proof}[Proof of Theorem \ref{thm:exist1}]
    We rely on Lemma \ref{lem:limh0} in Section \ref{sec:aux}, which shows that ${\rm MISE}2(\nu)$ is a continuous function of $\nu$ such that ${\rm MISE}2(\nu)\to\infty$ as $\nu\to\infty$. This implies that it is possible to find $\nu_0>0$ such that ${\rm MISE}2(\nu)>{\rm MISE}2(0)$, for all $\nu\geq\nu_0$. But since ${\rm MISE}2$ is continuous, it must have a minimizer $\nu_{\rm MISE}$ on $[0,\nu_0]$. In particular, this implies that ${\rm MISE}2(\nu_{\rm MISE})\leq{\rm MISE}2(0)$. The choice of $\nu_0$ ensures that also ${\rm MISE}2(\nu_{\rm MISE})<{\rm MISE}2(\nu)$ for all $\nu\geq\nu_0$. Hence, $\nu_{\rm MISE}$ is indeed the global minimizer of ${\rm MISE}2(\nu)$.
\end{proof}

\begin{proof}[Proof of Theorem \ref{thm:numise}]
    This proof relies on the asymptotic expansions provided in Corollary \ref{cor:MISEh_big} in Section \ref{sec:aux}.

    If $m_1(f)>0$, then part \ref{cor1:a} of Corollary \ref{cor:MISEh_big} implies that $\nu^{-2}\omega_d\{{\rm MISE}2(\nu)-{\rm MISE}2(0)\}\to 2L'(0)m_1(f)/L(0)$ as $\nu\to0$. Hence, the assumption that $L'(0)<0$ ensures that ${\rm MISE}2(\nu)<{\rm MISE}2(0)$ for small enough $\nu$, so that it must be $\nu_{\rm MISE}>0$.

    If $m_1(f)=0$ but $m_2(f)>z_2$ then part \ref{cor1:b} of Corollary \ref{cor:MISEh_big} gives $\nu^{-4}\omega_d\{{\rm MISE}(h)-{\rm MISE}2(0)\}\to n^{-1}L'(0)^2z_2/L(0)^2-L''(0)\{m_2(f)-z_2\}/L(0)$ as $\nu\to0$. This limit is negative as soon as $n>L'(0)^2z_2/\{L''(0)L(0)[m_2(f)-z_2]\}$, which yields the desired statement.

    Finally, when  $m_1=z_1=0$, $m_2=z_2$, $m_3=z_3=0$ and $L=L_{\rm vMF}$, part \ref{cor1:c} of Corollary \ref{cor:MISEh_big} ensures that
        $\omega_d\{{\rm MISE}2(\nu)-{\rm MISE}2(0)\}=\nu^{4}(A+B\nu^{4})+o(\nu^{8}).
        $ as $\nu\to0$,  where $A=n^{-1}z_2$ and $B=\tfrac{1}{12}\big\{n^{-1}(7z_4-15z_2^2)-[m_4(f)-z_4]\big\}.$
    Then, for every
    $n>(7z_4-15z_2^2)/\{m_4(f)-z_4\}$
    we have $B<0$, and we already know that $A>0$, so that $\nu_n\defin(-2A/B)^{1/4}>0$ is such that $\nu_n\sim\big(\frac{24z_2}{m_4-z_4}\big)^{1/4}n^{-1/4}\to0$ as $n\to\infty$. Hence,
    $$\omega_d\{{\rm MISE}2(\nu_n)-{\rm MISE}2(0)\}=2A^2/B+o(n^{-2})=-\nu_n^4A+o(n^{-2})\sim\,-\tfrac{24z_2^2}{m_4-z_2}n^{-2}.$$
    Therefore, $n^2\omega_d\{{\rm MISE}2(\nu_n)-{\rm MISE}2(0)\}\to-\tfrac{24z_2^2}{m_4-z_2}<0$ as $n\to\infty$ so that ${\rm MISE}2(\nu_n)<{\rm MISE}2(0)$ for large enough $n$, which implies that $\nu_{\rm MISE}>0$.
\end{proof}

\begin{proof}[Proof of Lemma \ref{lemma:varCV}]
	The cross-validation criterion, aside from a non-stochastic term, is a $U$-statistic of order 2 with the kernel $\varphi(\bX_1,\bX_2)=K_h(\bX_1,\bX_2)$. To prove the result, it suffices to use \eqref{eq:varU}, followed by the calculation of $\mathbb E\{K_h(\bX_1,\bX_2)\}$, $\mathbb E\{K_h(\bX_1,\bX_2)^2\}$ and $ \mathbb E\{K_h(\bX_1,\bX_2)K_h(\bX_1,\bX_3)\}$.

    From \eqref{eq:expker}, it follows that $\mathbb E\{K_h(\bX_1,\bX_2)\}=R_{K_h}(f)$. In addition,
    \begin{align*}
    \mathbb E\{K_h(\bX_1,\bX_2)^2\}&=\int_{\Sd} \int_{\Sd} K_h(\bx,\by)^2 f(\bx) f(\by)\,\sigma_d(\rd \bx)\,\sigma_d(\rd \by)=R_{(K_h)^2}(f).
    \end{align*}

    Finally,
    \begin{align*}
    \mathbb E\{&K_h(\bX_1,\bX_2)K_h(\bX_1,\bX_3)\}\\
    &=
   \int_{\Sd} \int_{\Sd}\int_{\Sd} K_h(\bx,\by)K_h(\bx,\bz) f(\bx) f(\by)f(\bz)\,\sigma_d(\rd \bx)\,\sigma_d(\rd \by)\,\sigma_d(\rd \bz)\\
   &=\int_{\Sd} (K_h*f)(\bx)^2 f(\bx)\,\sigma_d(\rd \bx)=S_{K_h}(f).
    \end{align*}
\end{proof}

\begin{proof}[Proof of Lemma \ref{lemma:varACV}]
    Let us denote $\bar Y=n^{-1}\sum_{i=1}^n f(\bX_i)$, so that ${\rm ACV}(h)={\rm CV}(h)+2\bar Y-2R(f)$. Then, $\mathbb{V}\mathrm{ar}\{{\rm ACV}(h)\}=\mathbb{V}\mathrm{ar}\{{\rm CV}(h)\}+4n^{-1}\mathbb{V}\mathrm{ar}\{f(\bX_1)\}+4\mathbb{C}\mathrm{ov}\{\bar Y,{\rm CV}(h)\}$. The last term involves the covariance between a $U$-statistic of order 1 (a sample average) and a $U$-statistic of order 2, for which \citet[Section 1.4]{lee1990u} has an exact formula. Precisely,
    \begin{align*}
        \mathbb{C}\mathrm{ov}\{\bar Y,{\rm CV}(h)\}&=2n^{-1}\mathbb{C}\mathrm{ov}\{f(\bX_1),K_h(\bX_1,\bX_2)\}\\
        &=2n^{-1}\lrb{\int_{\Sd} K_h*f(\bx)f(\bx)^2\,\sigma_d(\rd \bx)-R(f)R_{K_h}(f)}.
    \end{align*}
\end{proof}

\begin{proof}[Proof of Lemma \ref{lemma:varCVa}] To prove part $i)$ we use Lemma \ref{lemma:varCV} and compute the asymptotic expressions of $R_{K_h}(f)$,  $S_{K_h}(f)$, and $R_{(K_h)^2}(f)$.

First, using Corollary \ref{cor:ord4} twice and Corollary \ref{cor:RLhf}, as the bandwidth $h$ vanishes, it follows that
\begin{align}R_{K_h}(f)=&\int_{\Sd} \int_{\Sd} K_h(\bx,\by) f(\bx) f(\by)\,\sigma_d(\rd \bx)\,\sigma_d(\rd \by)\nonumber\\
=&\int_{\Sd} \int_{\Sd}\int_{\Sd} L_h(\bx,\bz) L_h(\by,\bz)f(\bx) f(\by)\,\sigma_d(\rd \bx)\,\sigma_d(\rd \by) \,\sigma_d(\rd \bz)\nonumber\\
&-2\int_{\Sd} \int_{\Sd} L_h(\bx,\by) f(\bx) f(\by)\,\sigma_d(\rd \bx)\,\sigma_d(\rd \by)\nonumber\\
\sim&\, R(f)-2R(f)=-R(f).\label{eq:RKf}\end{align}

On the other hand, using Corollary \ref{cor:ord4_2} twice, $S_{K_h}(f)\sim\int_{\Sd} f(\bx)^3\,\sigma_d(\rd \bx)$ as $h\to0$.

The most involved part concerns the asymptotic expression for $R_{(K_h)^2}(f)$. Next we will show that $R_{(K_h)^2}(f)\sim a_d(L)R(f)h_n^{-d}$. Start by writing
  \begin{align}
    R_{(K_h)^2}(f)&=\int_{\Sd} \int_{\Sd} (\tilde L_h-2L_h)(\bx,\by)^2 f(\bx) f(\by)\,\sigma_d(\rd \bx)\,\sigma_d(\rd \by)\nonumber\\
    &=\int_{\Sd} \int_{\Sd} \tilde L_h(\bx,\by)^2 f(\bx) f(\by)\,\sigma_d(\rd \bx)\,\sigma_d(\rd \by)\nonumber\\
    &\quad-4\int_{\Sd} \int_{\Sd} \tilde L_h(\bx,\by)L_h(\bx,\by) f(\bx) f(\by)\,\sigma_d(\rd \bx)\,\sigma_d(\rd \by)\nonumber\\
    &\quad+4\int_{\Sd} \int_{\Sd} L_h(\bx,\by)^2 f(\bx) f(\by)\,\sigma_d(\rd \bx)\,\sigma_d(\rd \by)\nonumber\\
    &=R_1-4R_2+4R_3.\label{eq:RK2dec}
    \end{align}
          Regarding the first term in \eqref{eq:RK2dec}, applying Proposition \ref{prop:AhBhChDh}, it follows that
    \begin{align*}
        R_1&=\int_{\Sd} \int_{\Sd} \tilde L_h(\bx,\by)^2 f(\bx) f(\by)\,\sigma_d(\rd \bx)\,\sigma_d(\rd \by)\nonumber\\
        &=\int_{\Sd} \int_{\Sd} (L_h*L_h)(\bx,\by) (L_h*L_h)(\bx,\by)f(\bx) f(\by)\,\sigma_d(\rd \bx)\,\sigma_d(\rd \by)\nonumber\\
        & \sim h^{3d}c^4_{d,L}(h)R(f)  \tilde{\gamma}_d  \int_0^{\infty}\Bigg\{ \int_0^{\infty} L(r) r^{d/2-1}\varphi_d(L, r, s)\, \rd r\Bigg\}^2s^{d/2-1}\, \rd s\nonumber\\
        & \sim h^{-d}\lambda^{-4}_d(L)R(f)  \tilde{\gamma}_d  \int_0^{\infty}\Bigg\{ \int_0^{\infty} L(r) r^{d/2-1}\varphi_d(L, r, s)\, \rd r\Bigg\}^2s^{d/2-1}\, \rd s.
    \end{align*}
    For the second term in \eqref{eq:RK2dec}, using Proposition \ref{prop:AhBhCh}, it follows that
    \begin{align*}
       R_2&= \int_{\Sd} \int_{\Sd} \tilde L_h(\bx,\by)L_h(\bx,\by) f(\bx) f(\by)\,\sigma_d(\rd \bx)\,\sigma_d(\rd \by)\nonumber\\
        &=\int_{\Sd} \int_{\Sd} (L_h*L_h)(\bx,\by)L_h(\bx,\by) f(\bx) f(\by)\,\sigma_d(\rd \bx)\,\sigma_d(\rd \by)\nonumber\\
        &\sim h^{2d} c^3_{d,L}(h)R(f)  \gamma_d \int_0^{\infty} L(r) r^{d/2-1}\left\{\int_0^{\infty}  L(s)s^{d/2-1} \varphi_d(L, r, s)\, \rd s\right\}\, \rd r\nonumber\\
        &\sim h^{-d} \lambda^{-3}_d(L)R(f)  \gamma_d \int_0^{\infty} L(r) r^{d/2-1}\left\{\int_0^{\infty}  L(s)s^{d/2-1} \varphi_d(L, r, s)\, \rd s\right\}\, \rd r.
    \end{align*}
   For the last term in \eqref{eq:RK2dec}, using Corollary \ref{cor:AhBh3}, we have
    \begin{align}
       R_3&=\int_{\Sd} \int_{\Sd} L_h(\bx,\by)^2 f(\bx) f(\by)\,\sigma_d(\rd \bx)\,\sigma_d(\rd \by)\nonumber\\
       &\sim h^{-d}\lambda_d(L^2)\lambda_d^{-2}(L)R(f).
        \label{eq:R3}
    \end{align}

Combining \eqref{eq:RK2dec}--\eqref{eq:R3}, it follows that  \begin{align*}
R_{(K_h)^2}(f)&\sim h^{-d}R(f)\lambda_d^{-2}(L)\Bigg\{\lambda_d^{-2}(L) \tilde{\gamma}_d  \int_0^{\infty}\left\{\int_0^{\infty} L(r) r^{d/2-1}\varphi_d(L, r, s)\, \rd r\right\}^2s^{d/2-1}\, \rd s\\
&-4\lambda_d^{-1}(L)\gamma_d \int_0^{\infty} L(r) r^{d/2-1}\left\{\int_0^{\infty}  L(s)s^{d/2-1} \varphi_d(L, r, s)\, \rd s\right\}\, \rd r+4\lambda_d(L^2)\Bigg\}\\
&= a_d(L)R(f)h_n^{-d}.\end{align*}

 Consequently, from Lemma \ref{lemma:varCV} we have
\begin{align}\label{eq:AVarCV}
	\mathbb{V}\mathrm{ar}\{{\rm CV}(h_n)\}&\sim 2n^{-2}h_n^{-d}a_d(L)R(f)
    +4n^{-1}\big\{\textstyle \int_{\Sd} f(\bx)^3\,\sigma_d(\rd \bx)-R(f)^2\big\}.
\end{align}
And, under the assumption that $nh_n\to\infty$, the first term in the previous expression is of smaller order than the second one, so part $i)$ follows.

Regarding part $ii)$, now use Lemma \ref{lemma:varACV} and \eqref{eq:AVarCV}, so that the only thing left is to find an asymptotic form for the term $\int_{\Sd} (K_h*f)(\bx)f(\bx)^2\,\sigma_d(\rd \bx)$ involved in the exact variance of ${\rm ACV}(h)$. But, proceeding as in the proof of \eqref{eq:RKf}, it can be checked that $\int_{\Sd} (K_h*f)(\bx)f(\bx)^2\,\sigma_d(\rd \bx)\sim-\int_{\Sd} f(\bx)^3 \,\sigma_d(\rd \bx)$, which entails
\begin{align}\label{eq:Khff2}
\int_{\Sd} (K_h*f)(\bx)f(\bx)^2\,\sigma_d(\rd \bx)-R(f)R_{K_h}(f)\sim-\left\{\int_{\Sd} f(\bx)^3\,\sigma_d(\rd \bx)-R(f)^2\right\}.
\end{align}
Combining the exact formula in Lemma \ref{lemma:varACV} with \eqref{eq:AVarCV} and \eqref{eq:Khff2} yields the desired asymptotic form for $\mathbb{V}\mathrm{ar}\{{\rm ACV}(h_n)\}$.
\end{proof}

\begin{proof}[Proof of Theorem \ref{th:consist}]
From \eqref{eq:ACV} and $\eqref{eq:amise}$, for any $c>0$ we have  $n^{4/(d+4)}{\rm ACV}(cn^{-1/(d+4)})\stackrel{\mathbb P}{\longrightarrow} A(c)-R(f)$, where $A(c)=c^{-d}v_d(L)+c^4b_d(L)^2R\lrpbig{\nabla^2 \bar{f}}$. Since $Z_n=n^{1/(d+4)}\hat h_{\rm CV}\in[\varepsilon,M]$, it follows that the sequence $Z_n$ is uniformly tight. By Prohorov's theorem \citep[see][Theorem 2.4]{vander1998cambridge}, there exists a random variable $Z$ such that $Z_n\stackrel{d}{\longrightarrow} Z$ along some subsequence, which we denote the same.

Using (\ref{eq:ACV}) again, along this subsequence we have
\begin{align*}
		n^{4/(d+4)}{\rm ACV}(\hat h_{\rm CV})&=n^{4/(d+4)}{\rm ACV}(Z_n n^{-1/(d+4)})\\
		&=Z_n^{-d}v_d(L)+Z_n^4b^{2}_d(L)R\lrpbig{\nabla^2 \bar{f}}-R(f)+o_\mathbb{P}(1)\\
        &\stackrel{d}{\longrightarrow}A(Z)-R(f),
\end{align*}
since $A(c)$ is a continuous function. By definition of $\hat h_{\rm CV}$, we also have ${\rm ACV}(\hat h_{\rm CV})\leq {\rm ACV}(cn^{-1/(d+4)})$ for all $c>0$, which implies $A(Z)\leq A(c)$ for all $c>0$. The function $A(c)$ has its unique minimum at $c=c_0$ given by \eqref{eq:hamise}, which implies $Z=c_0$.

Thus, a subsequence of  $Z_n=n^{1/(d+4)}\hat h_{\rm CV}$ converges in law to $c_0$. Since $c_0$ is a constant, this implies that $Z_n \stackrel{\mathbb P}{\longrightarrow} c_0$. Furthermore, any convergent subsequence of $Z_n$ must also converge to $c_0$, implying that $Z_n \stackrel{\mathbb P}{\longrightarrow} c_0$.

Therefore, this proves that $Z_n$ converges in probability, and that its limit is $c_0$. That is, $n^{1/(d+4)}\hat h_{\rm CV}\stackrel{\mathbb P}{\longrightarrow} c_0$ or, equivalently, $\hat h_{\rm CV}/h_0\stackrel{\mathbb P}{\longrightarrow}1$, which completes the proof since $h_0\sim h_M\sim h_{\rm MISE}$.
\end{proof}

\begin{proof}[Proof of Lemma \ref{lemma:varCVd}]
First, using formula \eqref{eq:varU}, the variance of ${\rm CV}'(h)$ is
\begin{align}
     \mathbb{V}\mathrm{ar}\{{\rm CV}'(h)\}&=\frac{2}{n(n-1)}\mathbb{V}\mathrm{ar}\{\nu_h(\bX_1,\bX_2)\}+\frac{4(n-2)}{n(n-1)}\mathbb{C}\mathrm{ov}\{\nu_h(\bX_1,\bX_2),\nu_h(\bX_1,\bX_3)\}\nonumber\\
     &=\frac{2}{n(n-1)}\{\E{\nu_h(\bX_1,\bX_2)^2}-\E{\nu_h(\bX_1,\bX_2)}^2\}\nonumber\\
     &\quad+\frac{4(n-2)}{n(n-1)}\{\E{\nu_h(\bX_1,\bX_2)\nu_h(\bX_1,\bX_3)}-\E{\nu_h(\bX_1,\bX_2)}^2\}\nonumber\\
     &\indef\frac{2}{n(n-1)}\{\mu_1(h)-\mu_2(h)^2\}+\frac{4(n-2)}{n(n-1)}\{\mu_3(h)-\mu_2(h)^2\}.\label{eq:CV1}
\end{align}

We now have to compute  $\mu_j(h)$, $j=1,2,3$, to obtain the exact order of $\mathbb{V}\mathrm{ar}\{{\rm CV}'(h)\}$. First, we have that \begin{align}
    \mu_2(h)=\mathbb E[\nu_h(\bX_1,\bX_2)]&=\int_{\Sd}\int_{\Sd} \nu_h(\bx,\by) f(\bx)f(\by)\,\sigmad(\rd \bx)\,\sigmad(\rd \by)\nonumber\\&=\int_{\Sd}(\nu_h*f)(\bx)f(\bx)\,\sigmad(\rd \bx)\nonumber\\
    &=O(h^3).\label{eq:mu2}
\end{align}
In \eqref{eq:mu2}, we use that
\begin{align}(\nu_h* f)(\bx)&=\int_{\Sd}\nu_h(\bx,\by)f(\by)\,\sigma_d(\rd \by)\nonumber\\
&=\frac{4c_{d,L}(h)}{h c_{d,G}(h)}\int_{\Sd}\big\{(L_h-G_h)*L_h-(L_h-G_h)\big\}(\bx,\by) f(\by)\,\sigma_d(\rd\by)\nonumber\\
&=\frac{4c_{d,L}(h)}{h c_{d,G}(h)}O(h^4)\nonumber\\
&=O(h^3),\label{eq:nuf}\end{align}
since using Corollaries \ref{cor:ord4} and \ref{cor:ord4_2} we have that
    \begin{align*}
        &\int_{\Sd}\big\{(L_h-G_h)*L_h-(L_h-G_h)\big\}(\bx,\by) \bar f(\by)\,\sigma_d(\rd\by)\\
        &=\int_{\Sd}\lrp{L_h*L_h-G_h*L_h-L_h+G_h}(\bx,\by)\bar f(\by)\,\sigma_d(\rd\by)\\
       &=\frac{h^2}{2d}\lrb{2\widetilde m_{d,2}(L,h)-[\widetilde m_{d,2}(G,h)+\widetilde m_{d,2}(L,h)]-\widetilde m_{d,2}(L,h)+\widetilde m_{d,2}(G,h)}{\rm tr}\,\bHcal\bar f(\bx)\\
       &\quad+O(h^4)\nonumber\\
       &=O(h^4).
    \end{align*}

Using \eqref{eq:nuf} twice, we have that
\begin{align}
    \mu_3(h)=&\; \E{\nu_h(\bX_1,\bX_2)\nu_h(\bX_1,\bX_3)}\nonumber\\
    =&\int_{\Sd}\int_{\Sd}\int_{\Sd} \nu_h(\bx,\by)  \nu_h(\bx,\bz) f(\bx) f(\by) f(\bz) \,\sigma_d(\rd \bx)\,\sigma_d(\rd \by)\,\sigma_d(\rd \bz)\nonumber\\
    =&\int_{\Sd} (\nu_h*f)(\bx)^2 f(\bx)\,\sigma_d(\rd \bx)\nonumber\\
    =&\; O(h^6). \label{eq:mu3}
\end{align}

Finally,
\begin{align}
  \mu_1(h)=&\;\E{\nu_h(\bX_1,\bX_2)^2}\nonumber\\
    =&\;\int_{\Sd}\int_{\Sd} \nu_h(\bx,\by)^2 f(\bx)f(\by)\,\sigmad(\rd \bx)\,\sigmad(\rd \by)\nonumber\\
    =&\;\frac{16c_{d,L}(h)^2}{h^2 c_{d,G}(h)^2}\nonumber\\
    &\times \int_{\Sd}\int_{\Sd}\{(L_h-G_h)*L_h-(L_h-G_h)\}(\bx,\by)^2f(\bx)f(\by)\,\sigmad(\rd \bx)\,\sigmad(\rd \by) \label{eq:mu11}\\
      =&\;\frac{16c_{d,L}(h)^2}{h^2 c_{d,G}(h)^2}O(h^{-d})\nonumber\\
       =&\;O(h^{-d-2}).\label{eq:mu1}
\end{align}

From \eqref{eq:CV1}, \eqref{eq:mu2}, \eqref{eq:mu3} and \eqref{eq:mu1}, it can be observed that $\mu_1(h_0)$ is the leading term of $\mathbb{V}\mathrm{ar}\{{\rm CV}'(h_0)\}$. Therefore, we will compute the asymptotic constant of this term, which will then be the asymptotic constant of $\mathbb{V}\mathrm{ar}\{{\rm CV}'(h_0)\}$.

For the calculation of \eqref{eq:mu11}, define
\begin{align*}
I_3(A,B,C)&\defin \int_0^{\infty}\int_0^\infty A(u)B(v)\varphi_d(C, u, v) (uv)^{d/2-1} \,\rd u \,\rd v,\\
I_4(A,B,C,D)&\defin \int_0^{\infty}\int_0^{\infty}\int_0^\infty A(u)B(v)\varphi_d(C, u, w)\varphi_d(D, v, w) (uvw)^{d/2-1} \,\rd u \,\rd v \,\rd w,
\end{align*}
where the function $\varphi_d$ is introduced in Proposition \ref{prop:AhBhCh}. Using Corollary \ref{cor:AhBh3} and Propositions \ref{prop:AhBhCh} and \ref{prop:AhBhChDh}, it follows that
\begin{align}
\int_{\Sd}\int_{\Sd}&\{(L_h-G_h)*L_h-(L_h-G_h)\}(\bx,\by)^2f(\bx)f(\by)\,\sigmad(\rd \bx)\,\sigmad(\rd \by)\nonumber\\
=&\int_{\Sd}\int_{\Sd}\{(L_h*L_h)-(G_h*L_h)-L_h+G_h\}(\bx,\by)^2f(\bx)f(\by)\,\sigmad(\rd \bx)\,\sigmad(\rd \by)\nonumber\\
=&\int_{\Sd}\int_{\Sd}(L_h*L_h)(\bx,\by)^2f(\bx)f(\by)\,\sigmad(\rd \bx)\,\sigmad(\rd \by)\nonumber\\
&+\int_{\Sd}\int_{\Sd}(G_h*L_h)(\bx,\by)^2f(\bx)f(\by)\,\sigmad(\rd \bx)\,\sigmad(\rd \by)\nonumber\\
&+\int_{\Sd}\int_{\Sd}L_h(\bx,\by)^2f(\bx)f(\by)\,\sigmad(\rd \bx)\,\sigmad(\rd \by)\nonumber\\
&+\int_{\Sd}\int_{\Sd}G_h(\bx,\by)^2f(\bx)f(\by)\,\sigmad(\rd \bx)\,\sigmad(\rd \by)\nonumber\\
&-2\int_{\Sd}\int_{\Sd}\{(L_h*L_h)(G_h*L_h)\}(\bx,\by)f(\bx)f(\by)\,\sigmad(\rd \bx)\,\sigmad(\rd \by)\nonumber\\
&-2\int_{\Sd}\int_{\Sd}\{(L_h*L_h)L_h\}(\bx,\by)f(\bx)f(\by)\,\sigmad(\rd \bx)\,\sigmad(\rd \by)\nonumber\\
&+2\int_{\Sd}\int_{\Sd}\{(L_h*L_h)G_h\}(\bx,\by)f(\bx)f(\by)\,\sigmad(\rd \bx)\,\sigmad(\rd \by)\nonumber\\
&+2\int_{\Sd}\int_{\Sd}\{(G_h*L_h)L_h\}(\bx,\by)f(\bx)f(\by)\,\sigmad(\rd \bx)\,\sigmad(\rd \by)\nonumber\\
&-2\int_{\Sd}\int_{\Sd}\{(G_h*L_h)G_h\}(\bx,\by)f(\bx)f(\by)\,\sigmad(\rd \bx)\,\sigmad(\rd \by)\nonumber\\
&-2\int_{\Sd}\int_{\Sd}L_h(\bx,\by)G_h(\bx,\by)f(\bx)f(\by)\,\sigmad(\rd \bx)\,\sigmad(\rd \by)\nonumber\\
\sim&\,\,\frac{R(f)\tilde{\gamma}_d}{\lambda_d(L)^4h^d}  I_4(L,L,L,L)+\frac{R(f)\tilde{\gamma}_d}{\lambda_d(L)^2\lambda_d(G)^2h^d}  I_4(L,L,G,G)+\frac{R(f)\lambda_d(L^2)}{\lambda_d(L)^2h^d}\nonumber\\
&+\frac{R(f)\lambda_d(G^2)}{\lambda_d(G)^2h^d}-\frac{2R(f)\tilde{\gamma}_d}{\lambda_d(L)^3\lambda_d(G)h^d}  I_4(L,L,L,G)-\frac{2R(f){\gamma}_d}{\lambda_d(L)^3h^d}  I_3(L,L,L)\nonumber\nonumber\\
&+\frac{4R(f){\gamma}_d}{\lambda_d(L)^2\lambda_d(G)h^d}  I_3(L,L,G)-\frac{2R(f){\gamma}_d}{\lambda_d(G)^2\lambda_d(L)h^d}  I_3(L,G,G)-\frac{2R(f)\lambda_d(LG)}{\lambda_d(L)\lambda_d(G)h^d},\label{eq:cons1}
\end{align}
where the constants $\gamma_d$ and $\tilde{\gamma}_d$ are defined in Propositions \ref{prop:AhBhCh} and \ref{prop:AhBhChDh}, respectively.

Notice that
\begin{align*}
    \lambda_{d}(G)=&\;2^{d/2-1}\om{d-1}\int_0^{\infty} L'(s) s^{d/2}\,\rd s
    =\;\begin{cases}
        \rd v=L'(s) \,\rd s,\quad v=L(s)\\
        u=s^{d/2},\quad \rd u=(d/2)s^{d/2-1}\,\rd s
    \end{cases}\\
    =&\;2^{d/2-1}\om{d-1}\lrc{L(s) s^{d/2}}_{0}^{\infty}-2^{d/2-1}\om{d-1}(d/2)\int_0^{\infty} L(s) s^{d/2-1}\,\rd s\\
    =&-(d/2)\lambda_d(L),\\
    \lambda_{d}(LG)=&\;2^{d/2-1}\om{d-1}\int_0^{\infty} L(s)L'(s) s^{d/2}\,\rd s =\;\begin{cases}
        \rd v=L(s)L'(s) \,\rd s,\quad v=L(s)^2/2\\
        u=s^{d/2},\quad \rd u=(d/2)s^{d/2-1}\,\rd s
    \end{cases}\\
     =&\;2^{d/2-2}\om{d-1}\lrc{L(s)^2 s^{d/2}}_{0}^{\infty}-2^{d/2-1}\om{d-1}(d/4)\int_0^{\infty} L(s)^2 s^{d/2-1}\,\rd s\\
     =&-(d/4)\lambda_d(L^2),
\end{align*}
and, therefore:
\begin{align}
\eqref{eq:cons1}=&\;\,\frac{R(f)\tilde{\gamma}_d}{\lambda_d(L)^4h^d}  I_4(L,L,L,L)+\frac{4R(f)\tilde{\gamma}_d}{d^2\lambda_d(L)^4h^d}  I_4(L,L,G,G)+\frac{4R(f)\lambda_d(G^2)}{d^2\lambda_d(L)^2h^d}\nonumber\\
&+\frac{4R(f)\tilde{\gamma}_d}{d\lambda_d(L)^4h^d}  I_4(L,L,L,G)-\frac{2R(f){\gamma}_d}{\lambda_d(L)^3h^d}  I_3(L,L,L)\nonumber\nonumber\\
&-\frac{8R(f){\gamma}_d}{d\lambda_d(L)^3h^d}  I_3(L,L,G)-\frac{8R(f){\gamma}_d}{d^2\lambda_d(L)^3h^d}  I_3(L,G,G)\nonumber\\
=&\;\,\frac{4R(f)\lambda_d(G^2)}{d^2\lambda_d(L)^2h^d}+\frac{R(f)\tilde{\gamma}_d}{\lambda_d(L)^4h^d}\{I_4(L,L,L,L)+4d^{-1}  I_4(L,L,L,G)+4d^{-2}  I_4(L,L,G,G)\}\nonumber\\
&-\frac{2R(f){\gamma}_d}{\lambda_d(L)^3h^d}\{I_3(L,L,L)+4d^{-1}I_3(L,L,G)+4d^{-2}I_3(L,G,G)\}. \label{eq:cons2}
\end{align}

In addition,
\begin{align*}
    I_3(G,L,L)&=\int_0^{\infty}\int_0^\infty G(u)L(v)\varphi_d(L, u, v) u^{d/2-1}v^{d/2-1}  \,\rd u \,\rd v\\
    &=\int_0^{\infty}\int_0^\infty L'(u)L(v)\varphi_d(L, u, v) u^{d/2}v^{d/2-1}  \,\rd u \,\rd v\\
    &=-\int_0^{\infty}\int_0^\infty L(u)L(v)\frac{\partial}{\partial u}\lrb{\varphi_d(L, u, v) u^{d/2}}v^{d/2-1}  \,\rd u \,\rd v\\
    &=-\frac d2\int_0^{\infty}\int_0^\infty L(u)L(v)\varphi_d(L, u, v) u^{d/2-1}v^{d/2-1}  \,\rd u \,\rd v\\
    &\quad-\int_0^{\infty}\int_0^\infty L(u)L(v)\lrb{\frac{\partial}{\partial u}\varphi_d(L, u, v)} u^{d/2}v^{d/2-1}  \,\rd u \,\rd v\\
    &=- \frac d2 I_3(L,L,L)-\int_0^{\infty}\int_0^\infty L(u)L(v)\lrb{\frac{\partial}{\partial u}\varphi_d(L, u, v)} u^{d/2}v^{d/2-1}  \,\rd u \,\rd v,
\end{align*}
\begin{align*}
     I_3(G,G,L)&=\int_0^{\infty}\int_0^\infty G(u)G(v)\varphi_d(L, u, v) u^{d/2-1}v^{d/2-1}  \,\rd u \,\rd v\\
    &=\int_0^{\infty}\int_0^\infty L'(u)L'(v)\varphi_d(L, u, v) u^{d/2}v^{d/2}  \,\rd u \,\rd v\\
    &=-\int_0^{\infty}\int_0^\infty L(u)L'(v)\frac{\partial}{\partial u}\lrb{\varphi_d(L, u, v) u^{d/2}}v^{d/2}  \,\rd u \,\rd v\\
    &=-\frac d2\int_0^{\infty}\int_0^\infty L(u)L'(v)\varphi_d(L, u, v) u^{d/2-1}v^{d/2}  \,\rd u \,\rd v\\
    &\quad-\int_0^{\infty}\int_0^\infty L(u)L'(v)\lrb{\frac{\partial}{\partial u}\varphi_d(L, u, v)} u^{d/2}v^{d/2}  \,\rd u \,\rd v\\
     &=\frac d2\int_0^{\infty}\int_0^\infty L(u)L(v)\lrb{\frac{\partial}{\partial v}\varphi_d(L, u, v)v^{d/2}}u^{d/2-1}  \,\rd u \,\rd v\\
    &\quad+\int_0^{\infty}\int_0^\infty L(u)L(v)\lrb{\frac{\partial}{\partial v}\lrb{\frac{\partial}{\partial u}\varphi_d(L, u, v)}v^{d/2} } u^{d/2} \,\rd u \,\rd v\\
    &=\frac{d^2}{4}\int_0^{\infty}\int_0^\infty L(u)L(v)\varphi_d(L, u, v)u^{d/2-1}v^{d/2-1}  \,\rd u \,\rd v\\
    &\quad +d\int_0^{\infty}\int_0^\infty L(u)L(v)\lrb{\frac{\partial}{\partial v}\varphi_d(L, u, v)}u^{d/2-1}v^{d/2}  \,\rd u \,\rd v\\
    &\quad+\int_0^{\infty}\int_0^\infty L(u)L(v)\lrb{\frac{\partial}{\partial v}\lrb{\frac{\partial}{\partial u}\varphi_d(L, u, v)} } u^{d/2}v^{d/2} \,\rd u \,\rd v\\
    &= \frac{d^2}{4} I_3(L,L,L) +d\int_0^{\infty}\int_0^\infty L(u)L(v)\lrb{\frac{\partial}{\partial v}\varphi_d(L, u, v)}u^{d/2-1}v^{d/2}  \,\rd u \,\rd v\\
    &\quad+\int_0^{\infty}\int_0^\infty L(u)L(v)\lrb{\frac{\partial}{\partial v}\lrb{\frac{\partial}{\partial u}\varphi_d(L, u, v)} } u^{d/2}v^{d/2} \,\rd u \,\rd v,
\end{align*}
\begin{align*}
    I_4(G,L,L,L)&=\int_0^{\infty}\int_0^{\infty}\int_0^\infty G(u)L(v)\varphi_d(L, u, w)\varphi_d(L, v, w) u^{d/2-1}v^{d/2-1}w^{d/2-1} \,\rd u \,\rd v \,\rd w\\
    &=\int_0^{\infty}\int_0^{\infty}\int_0^\infty L'(u)L(v)\varphi_d(L, u, w)\varphi_d(L, v, w) u^{d/2}v^{d/2-1}w^{d/2-1}  \,\rd u \,\rd v \,\rd w\\
    &=-\int_0^{\infty}\int_0^\infty \int_0^{\infty}L(u)L(v)\frac{\partial}{\partial u}\lrb{\varphi_d(L, u, w) u^{d/2}}\varphi_d(L, v, w)v^{d/2-1} w^{d/2-1} \,\rd u \,\rd v \,\rd w\\
    &=-\frac d2\int_0^{\infty}\int_0^{\infty}\int_0^\infty L(u)L(v)\varphi_d(L, u, w) \varphi_d(L, v, w)u^{d/2-1}v^{d/2-1}w^{d/2-1}  \,\rd u \,\rd v  \,\rd w\\
    &\quad-\int_0^{\infty}\int_0^{\infty}\int_0^\infty L(u)L(v)\lrb{\frac{\partial}{\partial u}\varphi_d(L, u, w)} \varphi_d(L, v, w)u^{d/2}v^{d/2-1} w^{d/2-1} \,\rd u \,\rd v \,\rd w\\
    &=- \frac d2 I_4(L,L,L,L)\\
    &\quad -\int_0^{\infty}\int_0^\infty \int_0^\infty L(u)L(v)\lrb{\frac{\partial}{\partial u}\varphi_d(L, u, w)}\varphi_d(L, v, w) u^{d/2}v^{d/2-1} w^{d/2-1}  \,\rd u \,\rd v,
\end{align*}
\begin{align*}
     I_4(G,G,L,L)&=\int_0^{\infty}\int_0^\infty G(u)G(v)\varphi_d(L, u, w) \varphi_d(L, v, w)u^{d/2-1}v^{d/2-1}w^{d/2-1}   \,\rd u \,\rd v\\
    &=\int_0^{\infty}\int_0^{\infty}\int_0^\infty L'(u)L'(v)\varphi_d(L, u, w) \varphi_d(L, v, w)u^{d/2}v^{d/2} w^{d/2-1}  \,\rd u \,\rd v \,\rd w\\
    &=-\int_0^{\infty}\int_0^{\infty}\int_0^\infty L(u)L'(v)\frac{\partial}{\partial u}\lrb{\varphi_d(L, u, w) u^{d/2}}\varphi_d(L, v, w)v^{d/2} w^{d/2-1}  \,\rd u \,\rd v \,\rd w\\
    &=-\frac d2\int_0^{\infty}\int_0^{\infty}\int_0^\infty L(u)L'(v)\varphi_d(L, u, w) \varphi_d(L, v, w)u^{d/2-1}v^{d/2} w^{d/2-1}  \,\rd u \,\rd v\,\rd w\\
    &\quad-\int_0^{\infty}\int_0^{\infty}\int_0^\infty L(u)L'(v)\lrb{\frac{\partial}{\partial u}\varphi_d(L, u, w)}\varphi_d(L, v, w) u^{d/2}v^{d/2}w^{d/2-1}   \,\rd u \,\rd v\,\rd w\\
     &=\frac d2\int_0^{\infty}\int_0^{\infty}\int_0^\infty L(u)L(v)\lrb{\frac{\partial}{\partial v}\varphi_d(L, v, w)v^{d/2}}\varphi_d(L, u, w)u^{d/2-1} w^{d/2-1}  \,\rd u \,\rd v\,\rd w\\
&\quad+\int_0^{\infty}\int_0^{\infty}\int_0^\infty L(u)L(v)\lrb{\frac{\partial}{\partial u}\varphi_d(L, u, w)}\lrb{\frac{\partial}{\partial v}\varphi_d(L, v, w)v^{d/2}  }u^{d/2}w^{d/2-1}  \,\rd u \,\rd v\,\rd w\\
    &=\frac{d^2}{4}\int_0^{\infty}\int_0^{\infty}\int_0^\infty L(u)L(v)\varphi_d(L, u, w)\varphi_d(L, v, w)u^{d/2-1}v^{d/2-1}w^{d/2-1}   \,\rd u \,\rd v\,\rd w\\
    &\quad +d\int_0^{\infty}\int_0^{\infty}\int_0^\infty L(u)L(v)\lrb{\frac{\partial}{\partial v}\varphi_d(L, v, w)}\varphi_d(L, u, w)u^{d/2-1}v^{d/2} w^{d/2-1}  \,\rd u \,\rd v\,\rd w\\
&\quad+\int_0^{\infty}\int_0^{\infty}\int_0^\infty L(u)L(v)\lrb{\frac{\partial}{\partial u}\varphi_d(L, u, w)}\lrb{\frac{\partial}{\partial v}\varphi_d(L, v, w)} u^{d/2}v^{d/2}w^{d/2-1}  \,\rd u \,\rd v\,\rd w\\
    &= \frac{d^2}{4} I_4(L,L,L,L)\\ &\quad +d\int_0^{\infty}\int_0^{\infty}\int_0^\infty L(u)L(v)\lrb{\frac{\partial}{\partial v}\varphi_d(L, v, w)}\varphi_d(L, u, w)u^{d/2-1}v^{d/2} w^{d/2-1}  \,\rd u \,\rd v\,\rd w\\
&\quad+\int_0^{\infty}\int_0^{\infty}\int_0^\infty L(u)L(v)\lrb{\frac{\partial}{\partial u}\varphi_d(L, u, w)}\lrb{\frac{\partial}{\partial v}\varphi_d(L, v, w)} u^{d/2}v^{d/2}w^{d/2-1}  \,\rd u \,\rd v\,\rd w.
\end{align*}

Therefore:
\begin{align*}
\eqref{eq:cons2}=&\;\,\frac{4R(f)\lambda_d(G^2)}{d^2\lambda_d(L)^2h^d}-\frac{8R(f){\gamma}_d}{d^2\lambda_d(L)^3h^d}\int_0^{\infty}\int_0^\infty L(u)L(v)\lrb{\frac{\partial}{\partial v}\lrb{\frac{\partial}{\partial u}\varphi_d(L, u, v)} } u^{d/2}v^{d/2} \,\rd u \,\rd v\nonumber\\
&+\frac{4R(f)\tilde{\gamma}_d}{d^2\lambda_d(L)^4h^d}\int_0^{\infty}\lrc{\int_0^\infty L(u)\lrb{\frac{\partial}{\partial u}\varphi_d(L, u, w)} u^{d/2}\, \rd u}^2w^{d/2-1}   \,\rd w.
\end{align*}

Consequently:
\begin{align*}
\eqref{eq:mu11}\sim & \frac{16 R(f)}{h^{d+2}}\Bigg\{\frac{\lambda_d(G^2)}{\lambda_d(L)^2}-\frac{2{\gamma}_d}{\lambda_d(L)^3}\int_0^{\infty}\int_0^\infty L(u)L(v) \varphi_{12,d}(L, u, v)u^{d/2}v^{d/2} \,\rd u \,\rd v\nonumber\\
&+\frac{\tilde{\gamma}_d}{\lambda_d(L)^4}\int_0^{\infty}\lrc{\int_0^\infty L(u) \varphi_{1,d}(L, u, w)u^{d/2}\, \rd u}^2w^{d/2-1}   \,\rd w\Bigg\},
\end{align*}
where $\varphi_{1,d}(L, u, v)=\frac{\partial}{\partial u}\varphi_d(L, u, v)$ and  $\varphi_{12,d}(L, u, v)=\frac{\partial}{\partial v}\frac{\partial}{\partial u}\varphi_d(L, u, v)$.\\

\color{black}

For $h_0=c_0n^{-1/(d+4)}$, given that $\mu_1(h_0)$ is the leading term of $\mathbb{V}\mathrm{ar}\{{\rm CV}'(h)\}$ in \eqref{eq:CV1}, and using \eqref{eq:mu11}, it follows that
\begin{align*}
     \mathbb{V}\mathrm{ar}\{{\rm CV}'(h_0)\}\sim 2\sigma^2_0(L) R(f)n^{-2}h_0^{-d-2},
\end{align*}
where
\begin{align*}
    \sigma^2_0(L)=&\;16\Bigg\{\frac{\lambda_d(G^2)}{\lambda_d(L)^2}-\frac{2{\gamma}_d}{\lambda_d(L)^3}\int_0^{\infty}\int_0^\infty L(u)L(v) \varphi_{12,d}(L, u, v)u^{d/2}v^{d/2} \,\rd u \,\rd v\nonumber\\
&+\frac{\tilde{\gamma}_d}{\lambda_d(L)^4}\int_0^{\infty}\lrc{\int_0^\infty L(u) \varphi_{1,d}(L, u, w)u^{d/2}\, \rd u}^2w^{d/2-1}   \,\rd w\Bigg\}.
\end{align*}

The asymptotic normality is a consequence of Theorem 2.1 in \cite{Jammalamadaka1986} for $U$-statistics of order 2 with varying kernel. The notations of that result correspond to setting $f_n(\bx,\by)=\textstyle{\binom{n}{2}}^{-1}\nu_{h_0}(\bx,\by)$ and $\sigma_n^2=\mathbb{V}\mathrm{ar}\{{\rm CV}'(h_0)\}$. Then, the asymptotic normality follows once we check the conditions (\textit{i}) $\sup_{\bx,\by}|f_n(\bx,\by)|=o(\sigma_n)$ and (\textit{ii}) $\sup_{\bx}\mathbb E|f_n(\bx,\bY)|=o(\sigma_n/n)$.

Regarding (\textit{i}), it is clear that for a bounded function $A$ we have $\sup_{\bx,\by}|A_h(\bx,\by)|=c_{d,A}(h)\sup_{t}|A(t)|$, which is of order $h^{-d}$. Analogously, it can be checked that for arbitrary bounded functions $A,B$ with finite $\lambda_d(A), \lambda_d(B)$ we have that $\sup_{\bx,\by}|(A_h*B_h)(\bx,\by)|$ is of order $h^{-d}$ as well. So, taking into account \eqref{eq:nu}, we have that $\sup_{\bx,\by}|f_n(\bx,\by)|$ is of order $n^{-2}h_0^{-d-1}$. Since the first part of this result shows that $\sigma_n$ is of order $n^{-1}h_0^{-(d+2)/2}$ and we have $h_0=c_0n^{-1/(d+4)}$, then (\textit{i}) immediately follows.
On the other hand, $\mathbb E|L_h(\bx,\bY)|=|L_h|*f(\bx)$ so, for a bounded $f$, it can be proved that $\sup_{\bx}\mathbb E|L_h(\bx,\bY)|$ is of order $O(1)$. Then, $\sup_{\bx}\mathbb E|f_n(\bx,\bY)|$ can be shown to be of order $n^{-2}h_0^{-1}$. Taking into account that $\sigma_n/n$ is of order $n^{-2}h_0^{-(d+2)/2}$, then condition (\textit{ii}) immediately holds.
\end{proof}

\begin{proof}[Proof of Theorem \ref{th:rate}]
The asymptotic distribution follows immediately from Lemma \ref{lemma:varCVd}, together with \eqref{eq:re2} and \eqref{eq:h0M2}. Moreover, by collecting the asymptotic representations in those results, the asymptotic variance can be written as
$$\sigma_d^2(L,f)=2c_0^{-d-2}c_1^{-2}\sigma^2_0(L) R(f).$$

Then, use that $c_0=\alpha(L)\delta(f)$, where
$
\alpha(L)=\lrb{dv_d(L)}^{1/(d+4)}\lrb{2b_d(L)}^{-2/(d+4)}$ and $\delta(f)=R(\nabla^2 \bar f)^{-1/(d+4)}
$,
and the fact that $c_1=\alpha(L)^{-d-1}\lrb{d(d+4)v_d(L)}\delta(f)^{-d-1}$,
to unpack the constant in the asymptotic variance as
\begin{align*}
    \sigma_d^2(L,f)=&\;2c_1^{-2}c_0^{-d-2}\sigma^2_0(L) R(f)\\
    =&\;2[\alpha(L)^{-d-1}\lrb{d(d+4)v_d(L)}\delta(f)^{-d-1}]^{-2}[\alpha(L)\delta(f)]^{-d-2}\sigma^2_0(L) R(f)\\
    =&\;\frac{2}{\lrb{d(d+4)}^2}\frac{\alpha(L)^{2(d+1)}\delta(f)^{2(d+1)}}{v_d(L)^2}\frac{1}{\alpha(L)^{d+2}\delta(f)^{d+2}}\sigma^2_0(L) R(f)\\
    =&\;\frac{2}{\lrb{d(d+4)}^2} \frac{\alpha(L)^{d}}{v_d(L)^2} \sigma^2_0(L) \delta(f)^{d}R(f)\\
    =&\;\frac{2}{\lrb{d(d+4)}^2} \frac{\big[\lrb{dv_d(L)}^{1/(d+4)}\lrb{2b_d(L)}^{-2/(d+4)}\big]^{d}}{v_d(L)^2} \sigma^2_0(L) R(\nabla^2 \bar f)^{-d/(d+4)}R(f)\\
    =&\; \frac{2}{\lrb{d(d+4)}^2} \frac{d^{d/(d+4)}v_d(L)^{d/(d+4)}}{v_d(L)^2 2^{2d/(d+4)}b_d(L)^{2d/(d+4)}} \sigma^2_0(L) \frac{R(f)}{R(\nabla^2 \bar f)^{d/(d+4)}}\\
    =&\;\frac{2d^{d/(d+4)}}{\lrb{d(d+4)}^2 2^{2d/(d+4)}} \frac{\sigma^2_0(L)}{v_d(L)^{(d+8)/(d+4)} b_d(L)^{2d/(d+4)}}  \frac{R(f)}{R(\nabla^2 \bar f)^{d/(d+4)}}\\
    =&\;\frac{1}{2^{(d-4)/(d+4)} d^{(d+8)/(d+4)}(d+4)^2} \frac{\sigma^2_0(L)}{[v_d(L)^{d+8} b_d(L)^{2d}]^{1/(d+4)}}  \frac{R(f)}{R(\nabla^2 \bar f)^{d/(d+4)}}\\
    =&\;\frac{1}{[2^{d-4} d^{d+8}]^{1/(d+4)}(d+4)^2} \times \frac{\sigma^2_0(L)}{[v_d(L)^{d+8} b_d(L)^{2d}]^{1/(d+4)}}  \times  \frac{R(f)}{R(\nabla^2 \bar f)^{d/(d+4)}}.
\end{align*}

\end{proof}

\begin{proof}[Proof of Corollary \ref{cor:vmf:kernel}]
To compute $\tau_d(L_\mathrm{vMF})$, we begin by calculating the term $\sigma^2_0(L_\mathrm{vMF})$ given in \eqref{eq:sigma2}:
\begin{align}   \sigma^2_0(L_\mathrm{vMF})=&\;16\Bigg\{\frac{\lambda_d(G_\mathrm{vMF}^2)}{\lambda_d(L_\mathrm{vMF})^2}-\frac{2{\gamma}_d}{\lambda_d(L_\mathrm{vMF})^3}\int_0^{\infty}\int_0^\infty L_\mathrm{vMF}(u)L_\mathrm{vMF}(v) \varphi_{12,d}(L_\mathrm{vMF}, u, v)u^{d/2}v^{d/2} \,\rd u \,\rd v\nonumber\\
&+\frac{\tilde{\gamma}_d}{\lambda_d(L_\mathrm{vMF})^4}\int_0^{\infty}\lrc{\int_0^\infty L_\mathrm{vMF}(u) \varphi_{1,d}(L_\mathrm{vMF}, u, w)u^{d/2}\,\rd u}^2w^{d/2-1} \,  \rd w\Bigg\},\label{eq:sigma0vMF}
\end{align}
where $G_\mathrm{vMF}(t)\defin L_\mathrm{vMF}'(t)t= -te^{-t}$.

For the calculation of the integrals involved in \eqref{eq:sigma0vMF}, we consider the cases $d=1$ and $d\ge 2$ separately. To begin with, let us suppose $d=1$. First, using \eqref{eq:phid}, we have that
\begin{align*}
    \varphi_{d}(L_\mathrm{vMF}, u, w)=&\; e^{-\left(u+w-2(u w)^{1/2}\right)}+e^{-\left(u+w+2(u w)^{1/2}\right)}=2e^{-\left(u+w\right)}\cosh(2(u w)^{1/2}),
\end{align*}
and, therefore
\begin{align*}
\varphi_{1,d}(L_\mathrm{vMF}, u, w)=&\;2u^{-1}e^{-\left(u+w\right)}[-u\cosh(2(u w)^{1/2})+(uw)^{1/2}\sinh(2(u w)^{1/2})],\\
\varphi_{12,d}(L_\mathrm{vMF}, u, w)=&\;(uw)^{-1/2}e^{-\left(u+w\right)}[4(u w)^{1/2}\cosh(2(u w)^{1/2})+(1-2u-2w)\sinh(2(u w)^{1/2})].
\end{align*}

Then, using the above expressions, we have that
\begin{align*}
\int_0^{\infty} L_\mathrm{vMF}(u) &  \varphi_{1,d}(L, u, w)u^{1/2} \,\rd u \\
=&\;\int_0^{\infty} 2u^{-1/2}e^{-\left(2u+w\right)}[-u\cosh(2(u w)^{1/2})+(uw)^{1/2}\sinh(2(u w)^{1/2})]\,\rd u\\
=&-2e^{-w}\int_0^{\infty} u^{1/2}e^{-2u}\cosh(2(u w)^{1/2})\,\rd u\\
&+2e^{-w}w^{1/2}\int_0^{\infty} e^{-2u}\sinh(2(u w)^{1/2})\,\rd u\\
=&-2e^{-w}\int_0^{\infty} x^{2}(e^{-2x^2+2x w^{1/2}}+e^{-2x^2-2x w^{1/2}})\,\rd x\\
&+e^{-w}w^{1/2}\int_0^{\infty} e^{-2u}(e^{2(u w)^{1/2}}-e^{-2(u w)^{1/2}})\,\rd u\\
=&-2e^{-w}\frac{\sqrt{\pi}}{4\sqrt{2}}(1+w)e^{w/2}+e^{-w}w^{1/2}\frac{\sqrt{\pi w}}{\sqrt{2}}e^{w/2}\\
=&\;(\pi/8)^{1/2}e^{-w/2}(w-1),
\end{align*}
and, consequently,
\begin{align}
\int_0^{\infty}\bigg[\int_0^\infty L_\mathrm{vMF}(u) &\varphi_{1,d}(L, u, w)u^{1/2}\, \rd u\bigg]^2w^{-1/2} \,  \rd w\nonumber\\
=&\;(\pi/8)\int_0^{\infty}e^{-w}(w-1)^2 w^{-1/2}  \, \rd w\nonumber\\
=&\;(\pi/8)\int_0^\infty e^{-w}w^{3/2} \,\rd w
-(\pi/4)\int_0^\infty e^{-w} w^{1/2} \,\rd w\nonumber
\\&+\, (\pi/8)\int_0^\infty e^{-w} w^{-1/2} \,\rd w\nonumber\\
=&\;(\pi/8)\Gamma(5/2)
-(\pi/4)\Gamma(3/2)+ (\pi/8)\Gamma(1/2)\nonumber\\
=&\;(\pi/8)3/4\pi^{1/2}-(\pi/4)\pi^{1/2}/2+ (\pi/8)\pi^{1/2}\nonumber\\
=&\;(3/32)\pi^{3/2}.\label{eq:LvMFphi1d1}
\end{align}

In addition,
\begin{align}
\int_0^{\infty}\int_0^\infty& L_\mathrm{vMF}(u)L_\mathrm{vMF}(v) \varphi_{12,d}(L, u, v)u^{1/2}v^{1/2} \,\rd u \,\rd v\nonumber
\\=&\;\int_0^{\infty}\lrc{\int_0^\infty e^{-2u}4(u v)^{1/2}\cosh(2(u v)^{1/2})\,\rd u} e^{-2v} \,\rd v\nonumber\\
&+\int_0^{\infty}\lrc{\int_0^\infty e^{-2u}(1-2u-2v)\sinh(2(u v)^{1/2})\,\rd u} e^{-2v} \,\rd v\nonumber\\
=&\;8\int_0^{\infty}\lrc{\int_0^\infty e^{-2t^2/v} t^2 \cosh(2t)  \, \rd t } \frac{1}{v} e^{-2v} \,\rd v\nonumber\\
&+2\int_0^{\infty}\lrc{\int_0^\infty e^{-\frac{2t^2}{v}}t \sinh(2t) \left(1-\frac{2t^2}{v}-2v\right) \, \rd t
} \frac{1}{v}e^{-2v} \,\rd v\nonumber\\
=&\;\sqrt{\frac{\pi}{2}}\int_0^{\infty}\lrc{v^{3/2} (1 + v) e^{v/2}} \frac{1}{v} e^{-2v} \,\rd v+\,-\sqrt{\frac{\pi}{2^5}}\int_0^{\infty}\lrc{v^{3/2} (1+5v)e^{v/2}
} \frac{1}{v}e^{-2v} \,\rd v\nonumber\\
=&\;\sqrt{\frac{\pi}{2}}\int_0^{\infty}v^{1/2}(1+v)e^{-3v/2}  \,\rd v -\sqrt{\frac{\pi}{2^5}}\int_0^{\infty} v^{1/2}(1+5v)e^{-3v/2}\,\rd v
\nonumber\\
=&\;\frac{2\pi}{3\sqrt{3}}-\frac{\pi}{2\sqrt{3}}\nonumber\\
=&\;\pi/(6\sqrt{3}).\label{eq:LvMFphi12d1}
\end{align}

Then, using \eqref{eq:LvMFphi1d1} and \eqref{eq:LvMFphi12d1} in \eqref{eq:sigma0vMF}, it follows that
\begin{align}\label{eq:sigma0vMFd1}
   \sigma^2_0(L_\mathrm{vMF})=\pi^{-1/2}(-4/9\sqrt{6}+3/2+3/16
   \sqrt{2}),
\end{align}
since from Lemma \ref{lemma:lambdaLGvmf} we have that $\lambda_1(L_\mathrm{vMF})=(2\pi)^{1/2}$ and $\lambda_1(G_\mathrm{vMF}^2)=\frac{3}{16}\pi^{1/2}$. Moreover, $\gamma_1=1$ and $\tilde{\gamma}_1=2^{-1/2}$.

For the calculation of the integrals in  \eqref{eq:sigma0vMF} in the case $d\ge 2$, similar steps can be followed.  First, using \eqref{eq:phid}, we have that
\begin{align*}
	\varphi_d(L_\mathrm{vMF}, u, w)&=e^{-\left(u+w\right)}
		\int_{-1}^1\left(1-\theta^2\right)^{(d-3)/2} e^{2 \theta(uw)^{1/2}} \,\rd \theta\nonumber\\
        &=e^{-\left(u+w\right)}\pi^{1/2}\Gamma(1/2(d-1))\left(uw \right)^{-(d-2)/4}\mathcal{I}_{(d-2)/2}\left(2(uw)^{1/2}\right)\nonumber\\
        &=:e^{-\left(u+w\right)}\pi^{1/2}\frac{\Gamma(1/2(d-1))}{\Gamma(d/2)}{}_0F_1(;d/2; uw)\nonumber\\
        &=e^{-\left(u+w\right)}2^{-(d-2)}\pi\frac{\Gamma(d-1)}{\Gamma(d/2)^2}{}_0F_1(;d/2; uw),
        \end{align*}
where ${}_0F_1(; b; z)$ denotes the confluent hypergeometric limit function. Therefore,
\begin{align}
\varphi_{1,d}(L_\mathrm{vMF},u,w)=&\;2^{-(d-2)}\pi\frac{\Gamma(d-1)}{\Gamma(d/2)^2}\nonumber\\
&\times\bigg[{}_0F_1(;d/2;uw)\frac{\partial}{\partial u}e^{-(u+w)}+e^{-(u+w)}\frac{\partial}{\partial u}{}_0F_1(;d/2;uw)\bigg]\nonumber\\
=&e^{-(u+w)}2^{-(d-2)}\pi\frac{\Gamma(d-1)}{\Gamma(d/2)^2}\nonumber\\
&\times\bigg[
-{}_0F_1(;d/2;uw)+\frac{2w}{d}{}_0F_1(;d/2+1;uw)
\bigg].\label{eq:phi1dvMF}
\end{align}

In addition,
\begin{align*}
\varphi_{12,d}(L_\mathrm{vMF},u,w)=&\;2^{-(d-2)}\pi\frac{\Gamma(d-1)}{\Gamma(d/2)^2}\nonumber\\
&\times\Bigg[
\left(-{}_0F_1(;d/2;uw)+\frac{2w}{d}{}_0F_1(;d/2+1;uw)\right)\frac{\partial}{\partial w}e^{-(u+w)}\nonumber\\
&+e^{-(u+w)}\frac{\partial}{\partial w}\left(-{}_0F_1(;d/2;uw)+\frac{2w}{d}{}_0F_1(;d/2+1;uw)\right)
\Bigg]\nonumber\\
=&\;e^{-(u+w)}2^{-(d-2)}\pi\frac{\Gamma(d-1)}{\Gamma(d/2)^2}\nonumber\\
&\times \Bigg[
-\left(-{}_0F_1(;d/2;uw)+\frac{2w}{d}{}_0F_1(;d/2+1;uw)\right)\nonumber\\
&+\bigg(
-\frac{u}{d/2}{}_0F_1(;d/2+1;uw)
+\frac{2}{d}{}_0F_1(;d/2+1;uw)\nonumber\\
&+\frac{2uw}{d(d/2+1)}{}_0F_1(;d/2+2;uw)
\bigg)\Bigg]\nonumber\\
=&\;e^{-(u+w)}2^{-(d-2)}\pi\frac{\Gamma(d-1)}{\Gamma(d/2)^2}\Bigg[
{}_0F_1(;d/2;uw)\nonumber\\
&-\frac{2u+2w-2}{d}{}_0F_1(;d/2+1;uw)
+\frac{2uw}{d(d/2+1)}{}_0F_1(;d/2+2;uw)
\Bigg]\nonumber\\
=&\;e^{-(u+w)}2^{-(d-2)}\pi\frac{\Gamma(d-1)}{\Gamma(d/2)^2}\nonumber\\
&\times\Bigg[
{}_0F_1(;d/2;uw)
-\frac{2u+2w-2}{d}{}_0F_1(;d/2+1;uw)\nonumber\\
&+[{}_0F_1(;d/2;uw)-{}_0F_1(;d/2+1;uw)]
\Bigg]\nonumber\\
=&\;e^{-(u+w)}2^{-(d-2)}\pi\frac{\Gamma(d-1)}{\Gamma(d/2)^2}\nonumber\\
&\times\Bigg[
2{}_0F_1(;d/2;uw)
-\frac{2u+2w-2+d}{d}{}_0F_1(;d/2+1;uw)\Bigg],
\end{align*}
since ${}_0F_1(;d/2+2;uw)=d(d+2)/(4uw)[{}_0F_1(;d/2;uw)-{}_0F_1(;d/2+1;uw)].$

Then, using \eqref{eq:phi1dvMF}, we have that
\begin{align*}
\int_0^{\infty} L_\mathrm{vMF}(u) \varphi_{1,d}(L, u, w)u^{d/2} \,\rd u =&\;2^{-(d-2)}\pi\frac{\Gamma(d-1)}{\Gamma(d/2)^2}e^{-w}\int_0^{\infty}e^{-2u}\nonumber\\
&\times\bigg[
-{}_0F_1(;d/2;uw)+\frac{2w}{d}{}_0F_1(;d/2+1;uw)
\bigg]u^{d/2} \,\rd u\nonumber\\
=&\;2^{-(d-2)}\pi\frac{\Gamma(d-1)}{\Gamma(d/2)^2}e^{-w}\big[-\Gamma(d/2)2^{-(d/2+2)}e^{w/2}(d+w)\nonumber\\
&+(2w/d)\Gamma(d/2+1)2^{-(d/2+1)}e^{w/2}\big]\nonumber\\
=&\;2^{-3d/2}\pi\frac{\Gamma(d-1)}{\Gamma(d/2)}e^{-w/2}(w-d).
\end{align*}
In the previous derivation, the following Laplace transforms evaluated at $s=2$ were applied:
\begin{align}
\mathcal{L}\{u^{d/2}{}_0F_1(;d/2; uw)\}(s)=&(1/2)\Gamma(d/2)s^{-(d/2+2)} (ds+2w)e^{w/s},\label{eq:lt1}\\
\mathcal{L}\{u^{d/2}{}_0F_1(; d/2+1; uw)\}(s)
=&\Gamma(d/2+1) s^{-(d/2+1)}e^{w/s}\label{eq:lt2}.
\end{align}

Consequently:
\begin{align}\int_0^{\infty}&\lrc{\int_0^\infty L_\mathrm{vMF}(u) \varphi_{1,d}(L, u, w)u^{d/2}\, \rd u}^2w^{d/2-1}  \, \rd w\nonumber\\
=&\;2^{-3d}\pi^2\frac{\Gamma(d-1)^2}{\Gamma(d/2)^2}\int_0^{\infty}e^{-w}(w^2+d^2-2dw)w^{d/2-1}   \,\rd w\nonumber\\
=&\;2^{-3d}\pi^2\frac{\Gamma(d-1)^2}{\Gamma(d/2)^2}[\Gamma(d/2+2)+d^2\Gamma(d/2)-2d\Gamma(d/2+1)]\nonumber\\
=&\;2^{-3d}\pi^2\frac{\Gamma(d-1)^2\Gamma(d/2+2)}{\Gamma(d/2)^2}.\label{eq:LvMFphi1d}
\end{align}

In addition,
\begin{align}
\int_0^{\infty}&\int_0^\infty L_\mathrm{vMF}(u)L_\mathrm{vMF}(v) \varphi_{12,d}(L, u, v)u^{d/2}v^{d/2} \,\rd u \,\rd v\nonumber\\
=&\;2^{-(d-2)}\pi\frac{\Gamma(d-1)}{\Gamma(d/2)^2}\nonumber\\
&\times \int_0^{\infty}\lrc{\int_0^\infty e^{-2u}\Bigg(
2{}_0F_1(;d/2;uv)
-\frac{2u+2v-2+d}{d}{}_0F_1(;d/2+1;uv)\Bigg)u^{d/2}\,\rd u}e^{-2v}v^{d/2}  \,\rd v\nonumber\\
=&\;2^{-(d-2)}\pi\frac{\Gamma(d-1)}{\Gamma(d/2)^2}\nonumber\\
&\times \int_0^{\infty}\Big[\Gamma(d/2)2^{-(d/2+1)}e^{v/2}(d+v)-\frac{2}{d}\Gamma(d/2+1) s^{-(d/2+3)}(d+v+2)e^{v/s}\nonumber\\
&-\frac{2v-2+d}{d}\Gamma(d/2+1)2^{-(d/2+1)}e^{v/2}\Big]e^{-2v}v^{d/2}  \,\rd v\nonumber\\
=&\;2^{-(d-2)}\pi\frac{\Gamma(d-1)}{\Gamma(d/2)^2}\nonumber\\
&\times \int_0^{\infty}\lrc{\Gamma(d/2)2^{-(d/2+1)}e^{v/2}(d+v)-\Gamma(d/2)2^{-(d/2+3)}e^{v/2}(3d+5v-2)}e^{-2v}v^{d/2}  \,\rd v\nonumber\\
=&\;2^{-(3d/2+1)}\pi\frac{\Gamma(d-1)}{\Gamma(d/2)} \int_0^{\infty}(d-v+2)e^{-3v/2}v^{d/2} \,\rd v\nonumber\\
=&\;2^{-(3d/2+1)}\pi\frac{\Gamma(d-1)}{\Gamma(d/2)}\lrc{(d + 2) \frac{\Gamma(d/2+1)}{(3/2)^{d/2+1}}-\frac{\Gamma(d/2+2)}{(3/2)^{d/2+2}}}\nonumber\\
=&\;2^{-d}3^{-(d/2+2)}d(d+2)\pi\Gamma(d-1).\label{eq:LvMFphi12d}
\end{align}

In \eqref{eq:LvMFphi12d}, we used \eqref{eq:lt1} and \eqref{eq:lt2}, and the following Laplace transform evaluated at $s=2$:
\begin{align*}
\mathcal{L}\{u^{d/2+1}{}_0F_1(; d/2+1; uw)\}(s)
=&(1/2)\Gamma(d/2+1) s^{-(d/2+3)}((2+d)s+2w)e^{w/s}.
\end{align*}

Then, using \eqref{eq:LvMFphi1d} and \eqref{eq:LvMFphi12d} in \eqref{eq:sigma0vMF}, it follows that
\begin{align}
   \sigma^2_0(L_\mathrm{vMF})=&\;\frac{d(d+2)\pi^{d/2}}{(2\pi)^d}-\frac{16\om{d-1}\om{d-2}2^{d-1}}{(2\pi)^{3d/2}}2^{-d}3^{-(d/2+2)}d(d+2)\pi\Gamma(d-1)\nonumber\\
&+\frac{16\om{d-1}\om{d-2}^22^{(3d-6)/2}}{(2\pi)^{2d}}2^{-3d}\pi^2\frac{\Gamma(d-1)^2\Gamma(d/2+2)}{\Gamma(d/2)^2}\nonumber\\
   =&\;\pi^{-d/2}d(d+2)(2^{-d}-2^{-(d/2-3)}3^{-(d/2+2)}+2^{-(3d/2+2)})\nonumber\\
   =&\;\pi^{-d/2}d(d+2)(2^{-d}-(8/9)2^{-d/2}3^{-d/2}+(1/4)2^{-3d/2)})\nonumber\\
   =&\;(2\pi)^{-d/2}d(d+2)(2^{-d/2}-(8/9)3^{-d/2}+(1/4)2^{-d})\nonumber\\
   =&\;2^{-d}\pi^{-d/2}d(d+2)(1+2^{-(d/2+2)}-2(3/2)^{-(d/2+2)}), \label{eq:sigma0vMFd}
\end{align}
since from Lemma \ref{lemma:lambdaLGvmf} we have that $\lambda_d(L_\mathrm{vMF})=(2\pi)^{d/2}$ and $\lambda_d(G_\mathrm{vMF}^2)=\frac{d(d+2)}{16}\pi^{d/2}$. Moreover, $\gamma_d=\om{d-1} \om{d-2} 2^{d-2}$ and $\tilde{\gamma}_d=\om{d-1}\om{d-2}^22^{(3d-6)/2}$, for $d \geq 2$.

Notice that although \eqref{eq:sigma0vMFd} was obtained for $d \geq 2$, if one sets $d = 1$ in that expression, it yields the same result as the one obtained in~\eqref{eq:sigma0vMFd1}.

Now, for the calculation of $\tau_d(L_\mathrm{vMF})$, using \eqref{eq:sigma0vMFd} and that
$b_d(L_{\mathrm{vMF}})=1/2$ and $v_d(L_{\mathrm{vMF}})=(2\pi^{1/2})^{-d}$, we have
\begin{align*}    \tau_d(L_\mathrm{vMF})=&\;\frac{1}{[2^{d-4} d^{d+8}]^{1/(d+4)}(d+4)^2} \frac{\sigma^2_0(L_\mathrm{vMF})}{[v_d(L_\mathrm{vMF})^{d+8} b_d(L_\mathrm{vMF})^{2d}]^{1/(d+4)}} \\
    =&\;\frac{2^{2d/(d+4)} (2\pi^{1/2})^{d(d+8)/(d+4)}}{[2^{d-4} d^{d+8}]^{1/(d+4)}(d+4)^2} \sigma^2_0(L_\mathrm{vMF})  \\
    =&\;\frac{2^{(5d+4)/(d+4)} \pi^{2d/(d+4)}}{ d^{(d+8)/(d+4)}} \lrb{1+2^{-(d/2+2)}-2(3/2)^{-(d/2+2)}} \frac{d(d+2)}{(d+4)^2}.
\end{align*}
\end{proof}

\begin{proof}[Proof of Corollary \ref{cor:vmf:dens}]
On the one hand, from \eqref{eq:vmf},
\begin{align*}
    R(f_{\mathrm{vMF}}(\cdot;\bmu, \kappa))=\frac{c_d^{\mathrm{vMF}}(\kappa)^2}{c_d^{\mathrm{vMF}}(2\kappa)}
    =\frac{\kappa^{(d-1)/2}\mathcal{I}_{(d-1)/2}(2\kappa)}{2^{d}\pi^{(d+1)/2}\mathcal{I}_{(d-1)/2}(\kappa)^2}.
\end{align*}
On the other hand, from Proposition 2 in \cite{Garcia-Portugues2013a},
\begin{align*}
    R(\nabla^2 \bar{f}_{\mathrm{vMF}}(\cdot;\bmu, \kappa))=\frac{d\kappa^{(d+1)/2}}{2^{d+2} \pi^{(d+1)/2} \mathcal{I}_{(d-1)/2}(\kappa)^2}\left[2 d \mathcal{I}_{(d+1)/2}(2 \kappa)+(d+2) \kappa \mathcal{I}_{(d+3)/2}(2 \kappa)\right].
\end{align*}
Combining both results, we have
\begin{align*}
    \rho_d(\kappa)&=\frac{R(f_{\mathrm{vMF}}(\cdot;\bmu, \kappa))}{R(\nabla^2 \bar{f}_{\mathrm{vMF}}(\cdot;\bmu, \kappa))^{d/(d+4)}}\\
    &=\frac{\kappa^{(d-1)/2}\mathcal{I}_{(d-1)/2}(2\kappa)}{2^{d}\pi^{(d+1)/2}\mathcal{I}_{(d-1)/2}(\kappa)^2}\\
    &\quad\times\lrb{\frac{2^{d+2} \pi^{(d+1)/2} \mathcal{I}_{(d-1)/2}(\kappa)^2}{d\kappa^{(d+1)/2}\left[2 d \mathcal{I}_{(d+1)/2}(2 \kappa)+(d+2) \kappa \mathcal{I}_{(d+3)/2}(2 \kappa)\right]}}^{d/(d+4)}\\
    &=\mathcal{I}_{(d-1)/2}(2\kappa)\\
    &\quad\times \frac{2^{-2d/(d+4)} \pi^{-2(d+1)/(d+4)} (\mathcal{I}_{(d-1)/2}(\kappa))^{-8/(d+4)}}{d^{d/(d+4)}\kappa^{(2 - d)/(d + 4)} \left[2 d \mathcal{I}_{(d+1)/2}(2 \kappa)+(d+2) \kappa \mathcal{I}_{(d+3)/2}(2 \kappa)\right]^{d/(d+4)}}\\
    &=\mathcal{I}_{(d-1)/2}(2\kappa)\\
    &\quad\times \lrb{4 \pi^{2(1+1/d)} d (\mathcal{I}_{(d-1)/2}(\kappa))^{8/d}  \kappa^{2/d - 1} \left[2 d \mathcal{I}_{(d+1)/2}(2 \kappa)+(d+2) \kappa \mathcal{I}_{(d+3)/2}(2 \kappa)\right]}^{-d/(d+4)}.
\end{align*}

We now turn our attention to \ref{cor:vmf:dens:kappa} and \ref{cor:vmf:dens:d}. For that, we use Equations 10.40.1 and 10.41.1 from \cite{NIST:DLMF}:
\begin{align*}
    \mathcal{I}_{\nu}(\kappa)\sim \frac{e^\kappa}{\sqrt{2\pi \kappa}} \text{ as }k\to\infty;\quad \mathcal{I}_{\nu}(\kappa)\sim \frac{e^\kappa}{\sqrt{2\pi \nu}} \lrp{\frac{e \kappa}{2\nu}}^{\nu}\text{ as }\nu\to\infty.
\end{align*}
First, as $\kappa\to\infty$,
\begin{align*}
    \rho_d(\kappa)&=\mathcal{I}_{(d-1)/2}(2\kappa)\\
    &\quad\times \lrb{4\pi^{2(1+1/d)} d (\mathcal{I}_{(d-1)/2}(\kappa))^{8/d}  \kappa^{2/d - 1} \left[2 d \mathcal{I}_{(d+1)/2}(2 \kappa)+(d+2) \kappa \mathcal{I}_{(d+3)/2}(2\kappa)\right]}^{-d/(d+4)}\\
    &\sim \frac{e^{2\kappa}}{\sqrt{4\pi \kappa}} \lrb{4 \pi^{2(1+1/d)} d \frac{e^{(8/d)\kappa}}{(2\pi \kappa)^{4/d}} \kappa^{2/d - 1} \left[2 d \frac{e^{2\kappa}}{\sqrt{4\pi \kappa}}+(d+2) \kappa \frac{e^{2\kappa}}{\sqrt{4\pi \kappa}}\right]}^{-d/(d+4)}\\
    &=\frac{1}{\sqrt{4\pi}} \lrb{\frac{ d(2d+(d+2)\kappa)}{2^{4/d-1}\pi^{2/d-3/2}}}^{-d/(d+4)} \kappa^{d/(d+4)}\\
    &=\lrb{4\pi^{2} d(2d+(d+2)\kappa)}^{-d/(d+4)}  \kappa^{d/(d+4)}\\
    &\sim\lrb{4\pi^{2} d(d+2)}^{-d/(d+4)}.
\end{align*}

Second, as $d\to\infty$,
\begin{align}
    \rho_d(\kappa)&=\mathcal{I}_{(d-1)/2}(2\kappa)\nonumber\\
    &\quad\times\lrb{4\pi^{2(1+1/d)} d (\mathcal{I}_{(d-1)/2}(\kappa))^{8/d}  \kappa^{2/d - 1} \left[2 d \mathcal{I}_{(d+1)/2}(2 \kappa)+(d+2) \kappa \mathcal{I}_{(d+3)/2}(2\kappa)\right]}^{-d/(d+4)}\nonumber\\
    &\sim \frac{e^{2\kappa}}{\sqrt{\pi (d-1)}} \lrp{\frac{2\kappa e}{d-1}}^{(d-1)/2} \nonumber\\
    &\quad\times\Bigg\{4\pi^{2(1+1/d)} d \lrc{\frac{e^{\kappa}}{\sqrt{\pi (d-1)}} \lrp{\frac{\kappa e}{d-1}}^{(d-1)/2}}^{8/d}  \kappa^{2/d - 1} \nonumber\\
    &\quad\times \left[2 d \frac{e^{2\kappa}}{\sqrt{\pi (d+1)}} \lrp{\frac{2\kappa e}{d+1}}^{(d+1)/2}+(d+2) \frac{e^{2\kappa}}{\sqrt{\pi (d+3)}} \lrp{\frac{2\kappa e}{d+3}}^{(d+3)/2}\right]\Bigg\}^{-d/(d+4)}\nonumber\\
    &\indef \frac{e^{2\kappa}}{\sqrt{\pi (d-1)}} \lrp{\frac{2\kappa e}{d-1}}^{(d-1)/2}\lrb{A\times B}^{-d/(d+4)}.\label{eq:rhod}
\end{align}
We compute first the leading term for $A\times B$. We start with $A$:
\begin{align*}
A&=2 d \frac{e^{2\kappa}}{\sqrt{\pi (d+1)}} \lrp{\frac{2\kappa e}{d+1}}^{(d+1)/2}+(d+2) \frac{e^{2\kappa}}{\sqrt{\pi (d+3)}} \lrp{\frac{2\kappa e}{d+3}}^{(d+3)/2}\\
&=d \frac{e^{2\kappa}}{\sqrt{\pi (d+1)}} \lrp{\frac{2\kappa e}{d+1}}^{(d+1)/2}\lrb{2 +\frac{d+2}{d} \sqrt{\frac{d+1}{d+3}} \lrp{\frac{d+1}{d+3}}^{(d+1)/2} \frac{2\kappa e}{d+3}}\\
&\sim d \frac{e^{2\kappa}}{\sqrt{\pi (d+1)}} \lrp{\frac{2\kappa e}{d+1}}^{(d+1)/2}\lrb{2 + \frac{2\kappa}{d+3}}\\
&\sim 2d \frac{e^{2\kappa}}{\sqrt{\pi (d+1)}} \lrp{\frac{2\kappa e}{d+1}}^{(d+1)/2}
\end{align*}
using
$\lrp{\frac{d+1}{d+3}}^{(d+1)/2}
\sim e^{-1}$. Similarly,
\begin{align*}
B&=4\pi^{2(1+1/d)} d \lrc{\frac{e^{\kappa}}{\sqrt{\pi (d-1)}} \lrp{\frac{\kappa e}{d-1}}^{(d-1)/2}}^{8/d}  \kappa^{2/d - 1}\\
&=4\pi^{2(1+1/d)} d \frac{e^{(8/d)\kappa}}{(\pi (d-1))^{4/d}} \lrp{\frac{\kappa e}{d-1}}^{4(d-1)/d} \kappa^{2/d - 1}\\
&\sim 4\pi^2 d \lrp{\frac{\kappa e}{d-1}}^{4} \kappa^{-1}\\
&\sim \frac{4\pi^2 \kappa^3 e^4}{d^3}.
\end{align*}
Consequently,
\begin{align}
    \lrb{A\times B}^{-d/(d+4)}&\sim \lrb{\frac{4\pi^2 \kappa^3 e^4}{d^3} 2d \frac{e^{2\kappa}}{\sqrt{\pi (d+1)}} \lrp{\frac{2\kappa e}{d+1}}^{(d+1)/2}}^{-d/(d+4)}\nonumber\\
    &\sim\lrb{\frac{8\pi^{3/2} \kappa^3 e^{2\kappa+4}}{d^{5/2}}\lrp{\frac{2\kappa e}{d+1}}^{(d+1)/2}}^{-d/(d+4)}. \label{eq:AB}
\end{align}
Replacing \eqref{eq:AB} in \eqref{eq:rhod} yields
\begin{align*}
    \rho_d(\kappa)&\sim \frac{e^{2\kappa}}{\sqrt{\pi (d-1)}} \lrp{\frac{2\kappa e}{d-1}}^{(d-1)/2} \lrb{\frac{8\pi^{3/2} \kappa^3 e^{2\kappa+4}}{d^{5/2}}  \lrp{\frac{2\kappa e}{d+1}}^{(d+1)/2}}^{-d/(d+4)}\\
    &=\frac{e^{2\kappa(1-d/(d+4))-4d/(d+4)}}{\sqrt{\pi}} \frac{\lrp{d+1}^{d(d+1)/(2d+8)}}{\lrp{d-1}^{d/2}}  \lrp{2\kappa e}^{(d-1)/2-d(d+1)/(2d+8)}\\
    &\quad\times \lrb{\frac{8\pi^{3/2} \kappa^3}{d^{5/2}}}^{-d/(d+4)}\\
    &\sim \frac{e^{-4}}{\sqrt{\pi}}
     \frac{\lrp{d+1}^{(d/2)(d+1)/(d+4)}}{\lrp{d-1}^{d/2}}\lrp{2\kappa e}
    \frac{1}{8\pi^{3/2} \kappa^3}d^{5d/(2d+8)}\\
    &=\frac{1}{4\pi^2e^3\kappa^2}
    d^{5d/(2d+8)} \lrp{\frac{d+1}{d-1}}^{d/2} \lrp{d+1}^{-3d/(2d+8)}\\
    &\sim \frac{1}{4\pi^2e^3 \kappa^2} d^{5/3}e d^{-3/2}\\
    &=\frac{d}{4(\pi e  \kappa)^2},
\end{align*}
concluding the proof.
\end{proof}

\section{Auxiliary results}
\label{sec:aux}

\begin{lemma}[Moments of the uniform distribution] \label{lem:moments}
Let $\bXi\sim\mathrm{Unif}(\mathbb{S}^d)$ and $r\geq1$. Then the vector of moments $\boldsymbol\zeta_{d,r}\defin\mathbb E[\bXi^{\otimes r}]\in\mathbb{R}^{(d+1)^r}$ is $\boldsymbol\zeta_{d,r}=c_{d,r}\,\boldsymbol{\mathcal{S}}_{d+1,r}(\vect\bI_{d+1})^{\otimes r/2}$ with
\begin{align*}
    c_{d,r}=\frac{\Gamma((d+1)/2)(r-1)!!}{2^{r/2}\Gamma((r+d+1)/2)}=\frac{(r-1)!!}{\prod_{k=0}^{r/2-1} (d+1+2k)}
\end{align*}
if $r$ is even and $\boldsymbol\zeta_{d,r}=\zero$ if $r$ is odd, where $(r-1)!!=(r-1)(r-3)\cdots5\cdot3\cdot1$. In particular,
\begin{align*}
    \boldsymbol\zeta_{d,2}=\frac{1}{d+1}\vect\bI_{d+1},\quad \boldsymbol\zeta_{d-1,4}=\frac{3}{(d+3)(d+1)}\,\boldsymbol{\mathcal{S}}_{d,4}(\vect\bI_d)^{\otimes2}.
\end{align*}
Here, $\boldsymbol{\mathcal{S}}_{d,r}$ denotes the $d^r\times d^r$ symmetrizer matrix \citep[see][p. 95]{Chacon2018}.
\end{lemma}

\begin{proof}[Proof of Lemma \ref{lem:moments}]
\citet[Equation (9.6.3)]{Mardia1999a} gives the connection between the mixed moments of the uniform distribution on $\Sd$ and those of a normal distribution $\mathcal{N}_{d+1}(\zero,\bI_{d+1})$:
\begin{align*}
    \mathbb{E}[\bXi^{\otimes r}]=\frac{\Gamma((d+1)/2)}{2^{r/2}\Gamma((r+d+1)/2)}\mathbb{E}[\bZ^{\otimes r}],
\end{align*}
where $\bZ\sim \mathcal{N}_{d+1}(\zero,\bI_{d+1})$. The result then follows from the formula for the multivariate normal central moments of arbitrary order given in \citet{Holmquist1988}: if $\bX\sim \mathcal{N}_{d+1}(\bmu,\bSigma)$, then $\mathbb E[(\bX-\bmu)^{\otimes r}]=(r-1)!!\,\boldsymbol{\mathcal{S}}_{d+1,r}(\vect\bSigma)^{\otimes r/2}$ if $r$ is even, and zero otherwise.
\end{proof}

\begin{lemma}\label{lem:limh0}
   Assume that \ref{A2} holds.
    \begin{enumerate}[label=\textit{(\alph{*})}, ref=\textit{(\alph{*})}]
        \item If $f$ is continuous, then $(L_h*f)(\bx)\to f(\bx)$ as $h\to0$ for every fixed $\bx\in\Sd$. \label{lem1:a}
        \item If $f$ is square integrable, then $h\mapsto R_{L_h}(f)$ and $h\mapsto R_{\tilde L_h}(f)$ are continuous functions such that $R_{L_h}(f)\to R(f)$ and $R_{\tilde L_h}(f)\to R(f)$ as $h\to0$.\label{lem1:b}
        \item If $L$ is continuous, then the function $h\mapsto v_{h,d}(L)$ is continuous and such that $v_{h,d}(L)\to \lambda_{d}(L^2)/\lambda_{d}(L)^2$ as $h\to0$.\label{lem1:c}
        \item If $L$ is continuous at 0, with $L(0)\neq0$, then $c_{d,L}(h)\to\{\omega_dL(0)\}^{-1}$ as $h\to\infty$. It also holds that $L_h(\bx,\by)\to\omega_d^{-1}$ and $(L_h*f)(\bx)\to\omega_d^{-1}$ as $h\to\infty$, uniformly in $\bx,\by\in\mathbb{S}^d$.\label{lem2:a}
         \item If $L$ is continuous at 0, with $L(0)\neq0$, then $R_{L_h}(f)\to\omega_d^{-1}$ and $R_{\tilde L_h}(f)\to\omega_d^{-1}$ as $h\to\infty$.\label{lem2:c}
        \item If $f$ is square integrable and $L$ is continuous, with $L(0)\neq0$, then ${\rm MISE}2(\nu)$ is a continuous function of $\nu$ such that ${\rm MISE}2(\nu)\to\infty$ as $\nu\to\infty$. \label{lem1:d}
    \end{enumerate}
\end{lemma}

\begin{proof}[Proof of Lemma \ref{lem:limh0}]
    First, we derive a convenient expression for the convolution $(L_h*f)(\bx)$. To begin with, employ the change of variables
    $$\by=t \bx+\left(1-t^2\right)^{1/2} \bB_{\bx} \bxi, \quad \sigma_d(\rd \by)=\left(1-t^2\right)^{d/2-1}\,\sigma_{d-1}(\rd \bxi)\,\rd t,$$
    where $\bx\in\Sd$, $\bxi\in\mathbb{S}^{d-1}$ and $\bB_{\bx}=\left(\mathbf{b}_1, \ldots, \mathbf{b}_d\right)_{(d+1) \times d}$ is the semi-orthonormal matrix resulting from the completion of $\bx$ to the orthonormal basis $\left\{\bx, \mathbf{b}_1, \ldots, \mathbf{b}_d\right\}$ of $\mathbb{R}^{d+1}$.  Then, continue with the change of variables $(1-t)/h^2=s$, $\rd t=-h^2\rd s$, to get
    \begin{align}
        (L_h*f)(\bx)&=c_{d,L}(h)\int_{\Sd}L\lp{\frac{1-\bx^\top\by}{h^2}}\rp f(\by)\,\sigmad(\rd \by)\nonumber\\
    &=c_{d,L}(h)\int_{-1}^1\int_{\mathbb{S}^{d-1}}L\lp{\frac{1-t}{h^2}}\rp {f}\lrp{t\bx+\left(1-t^2\right)^{1/2} \bB_{\bx} \bxi}\left(1-t^2\right)^{d/2-1}\,\sigma_{d-1}(\rd \bxi)\,\rd t\nonumber\\
    &=h^dc_{d,L}(h)\int_0^{2h^{-2}}\int_{\mathbb{S}^{d-1}}L(s)s^{d/2-1}(2-sh^2)^{d/2-1}{f}(\bx+\ba_{\bx,\bxi})\,\sigma_{d-1}(\rd\bxi)\,\rd s,\nonumber\\
    &=\lambda_{h,d}(L)^{-1}\int_0^{2h^{-2}}\int_{\mathbb{S}^{d-1}}L(s)s^{d/2-1}(2-sh^2)^{d/2-1}{f}(\bx+\ba_{\bx,\bxi})\,\sigma_{d-1}(\rd\bxi)\,\rd s,\label{eq:Lhf_int}
    \end{align}
    where $\ba_{\bx,\bxi}=-sh^2\bx+h[s(2-sh^2)]^{1/2}\bB_{\bx} \bxi$.
    Let $M_f$ be an upper bound for $f$ (since $f$ is continuous on $\Sd$, it is bounded). If $d\geq2$, the integrand appearing in \eqref{eq:Lhf_int} can be bounded by
    $$|L(s)s^{d/2-1}(2-sh^2)^{d/2-1}{f}(\bx+\ba_{\bx,\bxi})I_{[0,2h^{-2}]}(s)|\leq M_fL(s)(2s)^{d/2-1}I_{[0,\infty)}(s),$$
    which is an integrable function by assumption. So, taking into account that $f$ is continuous, \eqref{eq:equiv} and $\ba_{\bx,\bxi}\to\mathbf 0$ as $h\to0$, we can apply the DCT to conclude that
    $$\lim_{h\to0}(L_h*f)(\bx)=\lambda_d(L)^{-1}f(\bx)2^{d/2-1}\omega_{d-1}\int_0^\infty L(s) s^{d/2-1}\,\rd s=f(\bx),$$
    thus showing part \ref{lem1:a}.

    If $d=1$, proceeding as in the proof of Lemma 1 in \cite{Garcia-Portugues2013b}, split the integral in \eqref{eq:Lhf_int} into the sum of the integrals over $[0,h^{-2}]$ and $[h^{-2},2h^{-2})$. The first one has a bounded integrand
    $$|L(s)s^{-1/2}(2-sh^2)^{-1/2}{f}(\bx+\ba_{\bx,\xi})I_{[0,h^{-2})}(s)|\leq M_fL(s)(2s)^{-1/2}I_{[0,\infty)}(s)$$
    and hence, the same application of the DCT as for $d\geq 2$ follows. We show now that the second term,
    \begin{align}
    \sum_{\xi\in\{-1,1\}}\int_{h^{-2}}^{2h^{-2}} L(s)s^{-1/2}(2-sh^2)^{-1/2}{f}(\bx+\ba_{\bx,\xi})\,\rd s,\label{eq:secterm}
    \end{align}
    converges to zero as $h\to0$. The integrability assumption $\int_{0}^\infty L(s)s^{-1/2}\,\rd s<\infty$ implies that $L(s)=O\big(s^{-(1/2+\varepsilon)}\big)$ with $\varepsilon>0$ as $s\to\infty$, and hence $L(s)\leq C s^{-(1/2+\varepsilon)}$ for $C>0$ and $s$ large enough. Therefore,
    $$|L(s)s^{-1/2}(2-sh^2)^{-1/2}{f}(\bx+\ba_{\bx,\xi})I_{[h^{-2},2h^{-2})}(s)|\leq C M_f s^{-(1+\varepsilon)} (2-sh^2)^{-1/2} I_{[h^{-2},2h^{-2})}(s)$$
    for $h$ small enough. Since $\int_{h^{-2}}^{2h^{-2}} s^{-(1+\varepsilon)} (2-sh^2)^{-1/2} \,\rd s=O(h^{2\varepsilon})$, the DCT shows that \eqref{eq:secterm} converges to zero as $h\to0$.

    The key to showing part \ref{lem1:b} is to prove the desired properties in the case where $f$ is continuous (and bounded), and then extend the results for any square integrable density $f$.
    First, note that under the conditions of part \ref{lem1:a} it also follows that $|(L_h*f)(\bx)|\leq M_f$ for all $h>0$ and $\bx\in\Sd$. Therefore, we can apply the DCT again to obtain
    \begin{align*}
        \lim_{h\to0}R_{L_h}(f)&=\lim_{h\to0}\int_{\Sd}(L_h*f)(\bx)f(\bx)\,\sigma_d(\rd\bx)=\int_{\Sd}f(\bx)^2\,\sigma_d(\rd\bx)=R(f).
    \end{align*}
     Next, for an arbitrary square integrable $f$ and a given $\varepsilon>0$, take a continuous $g$ such that $\|g-f\|_2<\varepsilon$. Then,
     \begin{align}|R_{L_h}(f)-R(f)|\leq|R_{L_h}(f)-R_{L_h}(g)|+|R_{L_h}(g)-R(g)|+|R(g)-R(f)|.\label{eq:RLhineq}\end{align}
     The third term in \eqref{eq:RLhineq} can be made arbitrarily small since $|R(g)-R(f)|=|\|g\|_2^2-\|f\|_2^2|\leq \|g+f\|_2\|g-f\|_2$ and we showed above that the same happens to the second term as $h\to0$. Finally, for the first term in \eqref{eq:RLhineq}, note that
     $$f(\bx)f(\by)-g(\bx)g(\by)=(f-g)(\bx)f(\by)+(f-g)(\by)g(\bx).$$
    Hence,
     \begin{align}
        |R_{L_h}(f)-R_{L_h}(g)|
        &=\bigg|\int_{\Sd}\int_{\Sd}L_h(\bx,\by)[f(\bx)f(\by)-g(\bx) g(\by)]\,\sigma_d(\rd\bx)\,\sigma_d(\rd\by)\bigg|\nonumber\\
        &\leq\bigg|\int_{\Sd}\int_{\Sd}L_h(\bx,\by)(f-g)(\bx)f(\by)\,\sigma_d(\rd\bx)\,\sigma_d(\rd\by)\bigg|\nonumber\\
        &\quad+\bigg|\int_{\Sd}\int_{\Sd}L_h(\bx,\by)(f-g)(\by)g(\bx)\,\sigma_d(\rd\bx)\,\sigma_d(\rd\by)\bigg|.\label{eq:RLhf_RLhg}
     \end{align}
    Then, since $L_h$ is positive, we apply the Cauchy--Schwarz inequality to the functions $L_h(\bx,\by)^{1/2}(f-g)(\bx)$ and $L_h(\bx,\by)^{1/2}f(\by)$ to obtain
    \begin{align*}
        \Big|\int_{\Sd}\int_{\Sd}&L_h(\bx,\by)(f-g)(\bx)f(\by)\,\sigma_d(\rd\bx)\,\sigma_d(\rd\by)\Big|\\&\leq\Big\{\int_{\Sd}\int_{\Sd}L_h(\bx,\by)(f-g)(\bx)^2\,\sigma_d(\rd\bx)\,\sigma_d(\rd\by)\Big\}^{1/2}\\
        &\quad\times\Big\{\int_{\Sd}\int_{\Sd}L_h(\bx,\by)f(\by)^2\,\sigma_d(\rd\bx)\,\sigma_d(\rd\by)\Big\}^{1/2}\\
        &=\Big\{\int_{\Sd}(f-g)(\bx)^2\,\sigma_d(\rd\bx)\Big\}^{1/2}
        \Big\{\int_{\Sd}f(\by)^2\,\sigma_d(\rd\by)\Big\}^{1/2}\\
        &=\|f-g\|_2\|f\|_2,
    \end{align*}
    where we used the fact that $L_h$ is a density in either $\bx$ or $\by$ on $\Sd$. Analogously, the second integral in \eqref{eq:RLhf_RLhg} can be bounded by $\|f-g\|_2\|g\|_2$. Hence, $ |R_{L_h}(f)-R_{L_h}(g)|\leq (\|f\|_2+\|g\|_2)\|f-g\|_2$, which can also be made arbitrary small as $\varepsilon\to0$. This finally shows that $R_{L_h}(f)\to R(f)$ as $h\to0$ for any arbitrary square integrable $f$.

    The continuity can be shown similarly. Specifically, for any $h_0>0$, we have to prove that $R_{L_h}(f)\to R_{L_{h_0}}(f)$ as $h\to h_0$. The steps are analogous to those used in the previous reasoning, which corresponds to the limit case $h_0=0$.

    The proof for $R_{\tilde L_h}(f)$ follows a similar argument. On the one hand, since $R_{\tilde L_h}(f)=\int_{\Sd} (L_h*f)(\bx)^2\,\sigma_d(\rd \bx)$ and $|(L_h*f)(\bx)|\leq M_f$ when $f$ is bounded by $M_f$, part \ref{lem1:a} readily implies that $R_{\tilde L_h}(f)\to R(f)$ as $h\to0$ when $f$ is continuous (hence, bounded). On the other hand, applying the Cauchy--Schwarz inequality as before, we have
    \begin{align*}
        (L_h*f)(\bx)^2&=\Big\{\int_{\mathbb{S}^d}L_h(\bx,\by)f(\by)\,\sigma_d(\rd\by)\Big\}^2\\
        &\leq\Big\{\int_{\mathbb{S}^d}L_h(\bx,\by)\,\sigma_d(\rd\by)\Big\}\Big\{\int_{\mathbb{S}^d}L_h(\bx,\by)f(\by)^2\,\sigma_d(\rd\by)\Big\}\\
        &=\int_{\mathbb{S}^d}L_h(\bx,\by)f(\by)^2\,\sigma_d(\rd\by),
    \end{align*}
    which implies that $R(L_h*f)\leq R(f)$ for all $h>0$. Then, reasoning as before, we have
    \begin{align*}
        |R_{\tilde L_h}(f)-R_{\tilde L_h}(g)|&=|R(L_h*f)-R(L_h*g)|\\
        &\leq\|L_h*(f+g)\|_2\|L_h*(f-g)\|_2\\
        &\leq \|f+g\|_2\|f-g\|_2,
    \end{align*}
    thus showing that $R_{\tilde L_h}(f)-R_{\tilde L_h}(g)$ can be made arbitrarily small by choosing $g$ sufficiently close to $f$ in $L_2$ norm.

    Regarding part \ref{lem1:c}, notice that $v_{h,d}(L)=\lambda_{h,d}(L^2)\lambda_{h,d}(L)^{-2}$, so the assumptions on $L$ yield the desired limit. For the continuity, it suffices to show that the function $h\mapsto c_{d,L}(h)^{-1}$ is continuous. Fix an arbitrary $h_0>0$. With the same changes of variable as before, we can write
    \begin{align*}
        c_{d,L}(h)^{-1}&=\int_{\mathbb{S}^d}L\lp{\frac{1-\bx^\top\by}{h^2}}\rp\,\sigmad(\rd \bx)=\omega_{d-1}\int_{-1}^1L\lp{\frac{1-t}{h^2}}\rp \left(1-t^2\right)^{d/2-1}\,\rd t.
    \end{align*}
    For a fixed $t\in[-1,1]$, the continuity of $L$ ensures that the previous integrand is continuous as a function of $h$. Besides, we can find $E>0$ such that $0\leq(1-t)/h^2\leq E$ for all $h$ in a neighborhood $N_0$ of $h_0$ and all $t\in[-1,1]$ and, since $L$ is continuous on $[0,E]$, there is $M_L>0$ such that $|L(u)|\leq M_L$ for all $u\in[0,E]$. So, for all $(h,t)\in N_0\times[-1,1]$ we can bound
    $$\Big|L\Big(\frac{1-t}{h^2}\Big) \left(1-t^2\right)^{d/2-1}\Big|\leq M_L\left(1-t^2\right)^{d/2-1}.$$
    Since this bound does not depend on $h$ and is integrable on $[-1,1]$, the continuity follows from the DCT.

    To show part \ref{lem2:a}, note that $0\leq 1-\bx^\top\by\leq 2$ so that $(1-\bx^\top\by)/h^2\to0$ as $h\to\infty$, uniformly in $\bx,\by\in\Sd$, and therefore $L((1-\bx^\top\by)/h^2)\to L(0)$ as $h\to\infty$ by continuity. The uniform convergence easily implies that
    $$c_{d,L}(h)^{-1}=\int_{\mathbb{S}^d}L\Big(\frac{1-\bx^\top\by}{h^2}\Big)\,\sigma_d(\rd\bx)\to\int_{\mathbb{S}^d}L(0)\,\sigma_d(\rd\bx)=\omega_d L(0),\quad\text{as }h\to\infty.$$
    Hence, $L_h(\bx,\by)=c_{d,L}(h)L((1-\bx^\top\by)/h^2)\to\omega_d^{-1}$ as $h\to\infty$, uniformly in $\bx,\by\in\Sd$.
    The limit of $(L_h*f)(\bx)$ follows immediately from the previous uniform convergence and the inequality
    $$|(L_h*f)(\bx)-\omega_d^{-1}|=\Big|\int_{\mathbb{S}^d}\{L_h(\bx,\by)-\omega_d^{-1}\}f(\by)\,\sigma_d(\rd \by)\Big|\leq\sup_{\by}|L_h(\bx,\by)-\omega_d^{-1}|.$$

    To show part \ref{lem2:c}, use the uniform convergence in part \ref{lem2:a} together with the fact that we can write $R_{L_h}(f)=\int_{\mathbb{S}^d}(L_h*f)(\bx)f(\bx)\,\sigma_d(\rd\bx)$ and $R_{\tilde L_h}(f)=\int_{\mathbb{S}^d}(L_h*f)(\bx)^2\,\sigma_d(\rd\bx)$.

    The limit in part \ref{lem1:d} follows from the parts \ref{lem1:b}--\ref{lem1:c} and \eqref{eq:mise2}. Regarding the continuity, parts \ref{lem1:b}--\ref{lem1:c} guarantee the continuity of ${\rm MISE}2(\nu)$ at all $\nu\neq0$. It only remains to prove the continuity at $\nu=0$; that is, we need to show that ${\rm MISE}2(\nu)\to{\rm MISE}2(0)=\int_{\Sd}\{\omega_d^{-1}-f(\bx)\}^2\,\sigma_d(\rd \bx)=R(f)-\omega_d^{-1}$ as $\nu\to0$. But that limit immediately follows from parts \ref{lem2:a}--\ref{lem2:c} by taking into account that  $v_{h,d}(L)=h^{d}c_{d,L^2}(h)^{-1}c_{d,L}(h)^2$.
\end{proof}

The following result shows higher-order expansions of the terms of ${\rm MISE}2(\nu)$ as $\nu\to0$. Or, equivalently, of the terms of ${\rm MISE}(h)$ as $h\to\infty$. We will need to assume higher-order differentiability of $L$ at 0, and we will denote $L_j=L^{(j)}(0)/j!$ and $r_j=L_j/L_0$ for $j=0,1,2,\ldots$. It will also be useful to introduce the notation
 $$p_j=\sum_{\ell=0}^j(-1)^\ell\binom{j}{\ell}z_\ell=\sum_{\ell=0}^{\lfloor j/2\rfloor}\binom{j}{2\ell}z_{2\ell},$$
    with the last equality due to the fact that $z_\ell=0$ for odd $\ell$. For example, $p_0=p_1=1$, $p_2=1+z_2$, $p_3=1+3z_2$ and $p_4=1+6z_2+z_4$.

\begin{lemma}\label{lem:highorder}
Assume that $L$ is $4$-times continuously differentiable at 0, with $L(0)\neq0$. We have the following expansions as $h\to\infty$:
\begin{enumerate}[label=\textit{(\alph{*})}, ref=\textit{(\alph{*})}]
    \item $c_{d,L}(h)^{-1}=\omega_dL_0\sum_{j=0}^4b_jh^{-2j}+o(h^{-8}),$
    where $b_j=r_jp_j$.
    \item $c_{d,L}(h)=\omega_d^{-1}L_0^{-1}\sum_{j=0}^4 c_j h^{-2j}+o(h^{-8})$, where the coefficients $c_j$ are defined recursively as
    $c_0=1$, $c_j=-\sum_{k=1}^jb_kc_{j-k}$ for $j\geq1$; namely,
     $c_1=-b_1$, $c_2=b_1^2-b_2$, $c_3=-b_1^3+2b_1b_2-b_3$ and $c_4=b_1^4 - 3 b_1^2 b_2 + b_2^2 + 2 b_1 b_3 - b_4$.
    \item $L_h(\bx,\by)=\omega_d^{-1}\sum_{j=0}^4D_j(\bx^\top\by)h^{-2j}+o(h^{-8})$, where
    $D_j(t)=\sum_{\ell=0}^j\gamma_{j,\ell}t^\ell$,
    with $\gamma_{j,\ell}=(-1)^\ell \sum_{k=0}^{j-\ell}c_k r_{j-k}\binom{j-k}\ell$.  Similarly,  $(L_h*f)(\bx)=\omega_d^{-1}\sum_{j=0}^4\tilde D_j(\bx)h^{-2j}+o(h^{-8})$, where $\tilde D_j(\bx)=\sum_{\ell=0}^j\gamma_{j,\ell}\bmu_\ell(f)^\top\bx^{\otimes \ell}$. The remainders are uniform on $\bx,\by\in\mathbb{S}^d$.
\item $R_{L_h}(f)=\omega_d^{-1}\sum_{j=0}^4a_jh^{-2j}+o(h^{-8}),$ where $a_0=1$ and
$a_j=\sum_{\ell=1}^j\gamma_{j,\ell}(m_\ell-z_\ell)$ for all $j\geq 1$. Here and henceforth, $m_\ell$ is short for $m_\ell(f)$.
\item The expansion for $R_{\tilde L_h}(f)$ is
\begin{align*}
    R_{\tilde L_h}(f)&=\omega_d^{-1}\Big\{1+r_1^2m_1z_2h^{-4}+2r_1(2r_2-r_1^2)m_1z_2h^{-6}\\
    &\quad+\Big[3 m_1 r_1^4 z_2
-2 m_1 r_1^2 r_2 z_2^2
-10 m_1 r_1^2 r_2 z_2 + 2 m_1 r_1 r_3\big(3z_2+z_4\big) \\
&\quad + r_2^2\big(4 m_1 z_2 - z_2^2 + \tfrac{1}{3}(1+2m_2)z_4\big)\Big]h^{-8}\Big\}+o(h^{-8})
\end{align*}
If $m_1=z_1=0$ and $m_2=z_2$ then the coefficients of $h^{-4}$, $h^{-6}$ and $h^{-8}$ all vanish.
\item The expansion for $h^{-d}v_{h,d}(L)$ is
\begin{align*}
    h^{-d}v_{h,d}(L)
&=\omega_d^{-1}\Big\{1+r_1^2z_2h^{-4}+2r_1(2r_2-r_1^2)z_2h^{-6}+\\
&\quad+\big[3r_1^4z_2-r_1^2r_2(10z_2+2z_2^2)+2r_1r_3(3z_2+z_4)+r_2^2(4z_2+z_4-z_2^2)\big]h^{-8}\Big\}\\
&\quad+o(h^{-8}).
\end{align*}
\end{enumerate}
\end{lemma}

\begin{proof}[Proof of Lemma \ref{lem:highorder}]
    Employing the same changes of variables as in the proof of Lemma \ref{lem:limh0}, we can write
    \begin{align*}
        c_{d,L}(h)^{-1}&=\int_{\mathbb{S}^d}L\Big(\frac{1-\bx^\top\by}{h^2}\Big)\,\sigma_d(\rd\bx)=\omega_{d-1}\int_{-1}^1L\Big(\frac{1-t}{h^2}\Big)(1-t^2)^{d/2-1}\,\rd t\\
        &=\omega_{d-1}\sum_{j=0}^4L_jh^{-2j}\int_{-1}^1(1-t)^j(1-t^2)^{d/2-1}\,\rd t+o(h^{-8})\\
        &=\omega_{d-1}\sum_{j=0}^4L_jI_{j,d}h^{-2j}+o(h^{-8}),
    \end{align*}
    where $I_{j,d}=\int_{-1}^1(1-t)^j(1-t^2)^{d/2-1}\,\rd t$. It is not hard to obtain an explicit formula for $I_{j,d}$; however, it will be more convenient to relate that integral to the uniform moment norms $z_j$. We have
    \begin{align*}
        z_j=\boldsymbol \zeta_j^\top\boldsymbol \zeta_j&=\omega_d^{-2}\int_{\mathbb{S}^d}\int_{\mathbb{S}^d}(\bx^{\otimes j})^\top\by^{\otimes j}\,\sigma_d{(\rd\bx)}\,\sigma_d{(\rd\by)}\\
        &=\omega_d^{-2}\int_{\mathbb{S}^d}\int_{\mathbb{S}^d}(\bx^\top\by)^{j}\,\sigma_d{(\rd\bx)}\,\sigma_d{(\rd\by)}\\
        &=\omega_d^{-2}\omega_{d-1}\int_{\mathbb{S}^d}\int_{-1}^1t^{j}(1-t^2)^{d/2-1}\,\rd t\,\sigma_d{(\rd\by)}\\
         &=\omega_d^{-1}\omega_{d-1}\int_{-1}^1t^{j}(1-t^2)^{d/2-1}\,\rd t.
    \end{align*}
    Therefore,
    $$I_{j,d}=\sum_{\ell=0}^j(-1)^\ell\binom{j}{\ell}\int_{-1}^1t^{\ell}(1-t^2)^{d/2-1}\,\rd t=(\omega_d/\omega_{d-1})\sum_{\ell=0}^j(-1)^\ell\binom{j}{\ell}z_\ell=(\omega_d/\omega_{d-1})p_j.$$
    Hence, we can write $c_{d,L}(h)^{-1}=\omega_dL_0\sum_{j=0}^4r_jp_jh^{-2j}+o(h^{-8})$.

    The expression for $c_{d,L}(h)$ with the recursive definition of the coefficients $c_j$ follows the standard description for the multiplicative inverse of power series, after noting that $b_0=1$.

    Next, consider $L_h(\bx,\by)=c_{d,L}(h)L\big((1-t)/h^2\big)$, where $t=\bx^\top\by$. Expand
    $
        L\big((1-t)/h^2\big)=\sum_{j=0}^4L_j(1-t)^jh^{-2j}+o(h^{-8})
    $
    and multiply by $c_{d,L}(h)=\omega_d^{-1}L_0^{-1}\sum_{j=0}^4 c_j h^{-2j}+o(h^{-8})$ to obtain
    $L_h(\bx,\by)=\omega_d^{-1}\sum_{j=0}^4D_j(t)h^{-2j}+o(h^{-8}),$ where
       \begin{align*}
        D_j(t)&=L_0^{-1}\sum_{k=0}^jc_k L_{j-k}(1-t)^{j-k}
        =\sum_{k=0}^jc_k r_{j-k}\sum_{\ell=0}^{j-k}\binom{j-k}\ell(-1)^\ell t^\ell\\
        &=\sum_{\ell=0}^j\sum_{k=0}^{j-\ell}c_k r_{j-k}\binom{j-k}\ell(-1)^\ell t^\ell
        =\sum_{\ell=0}^jt^\ell(-1)^\ell \sum_{k=0}^{j-\ell}c_k r_{j-k}\binom{j-k}\ell\\&=\sum_{\ell=0}^j\gamma_{j,\ell}t^\ell,
    \end{align*}
    with $\gamma_{j,\ell}=(-1)^\ell \sum_{k=0}^{j-\ell}c_k r_{j-k}\binom{j-k}\ell$, so that $D_j(t)$ is a polynomial in $t$ of degree $\leq j$.
   Multiplying the previous expansion of $L_h(\bx,\by)$ by $f(\by)$ and integrating with respect to $\sigma_d(\rd\by)$, a direct consequence is that $(L_h*f)(\bx)=\omega_d^{-1}\sum_{j=0}^4\tilde D_j(\bx)h^{-2j}+o(h^{-8})$, where now $\tilde D_j(\bx)=\sum_{\ell=0}^j\gamma_{j,\ell}\bmu_\ell^\top\bx^{\otimes \ell}$, with the abbreviation $\bmu_\ell=\bmu_\ell(f)$ along this proof.

    Similarly, multiplying the previous expansion of $(L_h*f)(\bx)$ by $f(\bx)$ and integrating with respect to $\sigma_d(\rd\bx)$ we obtain $R_{L_h}(f)=\omega_d^{-1}\sum_{j=0}^4a_jh^{-2j}+o(h^{-8})$, where
$a_j=\sum_{\ell=0}^j\gamma_{j,\ell}m_\ell$. It immediately follows that $a_0=1$. Moreover, for $j\geq1$ we can write $$a_j=\sum_{\ell=0}^j\gamma_{j,\ell}m_\ell=\sum_{\ell=0}^j\gamma_{j,\ell}(m_\ell-z_\ell)+\sum_{\ell=0}^j\gamma_{j,\ell}z_\ell=\sum_{\ell=1}^j\gamma_{j,\ell}(m_\ell-z_\ell)+\sum_{\ell=0}^j\gamma_{j,\ell}z_\ell,$$
the last equality is due to the fact that $m_0=z_0=1$. Hence, to obtain the expression that is stated for $a_j$ for $j\geq1$, it suffices to prove that $\sum_{\ell=0}^j\gamma_{j,\ell}z_\ell=0$ for all $j\geq1$. But notice that for the uniform density $u(\bx)=\omega_{d}^{-1}$ we have $R_{L_h}(u)=\omega_d^{-1}$ for all $h$ and, at the same time, for any $q\in\mathbb N$ we have
$$R_{L_h}(u)=\omega_d^{-1}\sum_{j=0}^q\tilde a_jh^{-2j}+o(h^{-2q}),$$ where
$\tilde a_j=\sum_{\ell=0}^j\gamma_{j,\ell}z_\ell$. So it must be $\tilde a_j=0$ for all $j\geq1$, as desired.

Regarding $R_{\tilde L_h}(f)=\int_{\mathbb{S}^d}(L_h*f)(\bx)^2\,\sigma_d(\rd\bx)$, squaring the expansion of $(L_h*f)(\bx)$ we have
\begin{align*}
    R_{\tilde L_h}(f)
    &=\omega_d^{-2}\sum_{j=0}^4h^{-2j}\sum_{k=0}^j\int_{\mathbb{S}^d}\tilde D_k(\bx)\tilde D_{j-k}(\bx)\,\sigma_d(\rd\bx)+o(h^{-8})\\
    &=\omega_d^{-2}\Big\{\int_{\mathbb{S}^d}\tilde D_0(\bx)^2\,\sigma_d(\rd\bx)+h^{-2}
    \int_{\mathbb{S}^d}\big[2\tilde D_0(\bx)\tilde D_1(\bx)\big]\,\sigma_d(\rd\bx)\\
    &\quad+h^{-4}\int_{\mathbb{S}^d}\big[2\tilde D_0(\bx)\tilde D_2(\bx)+\tilde D_1(\bx)^2\big]\,\sigma_d(\rd\bx)\\
    &\quad+h^{-6}\int_{\mathbb{S}^d}\big[2\tilde D_0(\bx)\tilde D_3(\bx)+2\tilde D_1(\bx)\tilde D_2(\bx)\big]\,\sigma_d(\rd\bx)\\
    &\quad+h^{-8}\int_{\mathbb{S}^d}\big[2\tilde D_0(\bx)\tilde D_4(\bx)+2\tilde D_1(\bx)\tilde D_3(\bx)+\tilde D_2(\bx)^2\big]\,\sigma_d(\rd\bx)\\
    &\quad+o(h^{-8}).
\end{align*}
Let us compute each term explicitly:

Since $\gamma_{0,0}=1$ we have $\tilde D_0(\bx)=1$, and the coefficient of $h^0$ is $\omega_d^{-1}$.

We have $\gamma_{1,0}=0$, $\gamma_{1,1}=-r_1$ so that $\tilde D_1(\bx)=-r_1\bmu_1^\top\bx$, so the coefficient of $h^{-2}$ is $-2r_1\bmu_1^\top\boldsymbol\zeta_1=0$, because $\boldsymbol\zeta_1={\mathbf 0}$.

We have $\gamma_{2,0}=-r_2z_2$, $\gamma_{2,1}=r_1^2-2r_2$, $\gamma_{2,2}=r_2$. Hence, $\tilde D_2(\bx)=-r_2z_2+(r_1^2-2r_2)\bmu_1^\top\bx+r_2\bmu_2^\top\bx^{\otimes 2}$. Therefore, $\int_{\mathbb{S}^d}\tilde D_2(\bx)\,\sigma_d(\rd\bx)=\omega_d\{-r_2z_2+(r_1^2-2r_2)\bmu_1^\top\boldsymbol\zeta_1+r_2\bmu_2^\top\boldsymbol\zeta_2\}$. Again, simplify $\boldsymbol\zeta_1={\mathbf 0}$ and note that
\begin{align}\label{eq:mu2zeta2}
    \bmu_2^\top\boldsymbol\zeta_2=\frac{1}{d+1}\bmu_2^\top\vect\bI_{d+1}=\frac{1}{d+1}\int_{\mathbb{S}^d}\bx^\top\bx f(\bx)\,\sigma_d(\rd\bx)=\frac{1}{d+1}=z_2,
\end{align}
which leads to $\int_{\mathbb{S}^d}\tilde D_2(\bx)\,\sigma_d(\rd\bx)=0$. On the other hand,
\begin{align*}
\int_{\mathbb{S}^d}\tilde D_1(\bx)^2\,\sigma_d(\rd\bx)&=\omega_dr_1^2(\bmu_1^{\top})^{\otimes 2}\boldsymbol\zeta_2=\frac1{d+1}\omega_dr_1^2(\bmu_1^{\top})^{\otimes 2}\vect\bI_{d+1}\\
&=\frac1{d+1}\omega_dr_1^2\,\bmu_1^{\top}\bmu_1=\omega_dr_1^2m_1z_2,
\end{align*}
so eventually, the coefficient of $h^{-4}$ is $\omega_d^{-1}r_1^2m_1z_2$.

For the coefficient of $h^{-6}$ notice that we have $\int_{\mathbb{S}^d}\tilde D_3(\bx)\,\sigma_d(\rd\bx)=\omega_d(\gamma_{3,0}+\gamma_{3,2}z_2)$, because the terms involving $\gamma_{3,1}$ and $\gamma_{3,3}$ vanish, since they include $\boldsymbol\zeta_1=\mathbf 0$ and $\boldsymbol\zeta_3=\mathbf 0$. We also made use of $\bmu_2^\top\boldsymbol\zeta_2=z_2$, as shown in \eqref{eq:mu2zeta2}. Then, using the formulas for $c_i$, $i=0,1,2,3$ in terms of the $r_i$ and the $z_i$ (through the relationship between $p_i$ and $z_i$) shows that $\gamma_{3,0}+\gamma_{3,2}z_2=0$. Next, we need to compute
\begin{align*}
    \omega_d^{-1}\int_{\mathbb{S}^d}\tilde D_1(\bx)\tilde D_2(\bx)\,\sigma_d(\rd\bx)
    &=-r_1\omega_d^{-1}\sum_{\ell=0}^2\gamma_{2,\ell}\int_{\mathbb{S}^d}(\bmu_1^\top\bx)(\bmu_\ell^\top\bx^{\otimes\ell})\,\sigma_d(\rd\bx)\\
    &=-r_1\sum_{\ell=0}^2\gamma_{2,\ell}(\bmu_1^\top\otimes\bmu_\ell^\top)\boldsymbol\zeta_{\ell+1}\\
    &=-r_1\sum_{\ell=0}^2\gamma_{2,\ell}(\bmu_1^\top\otimes\bmu_\ell^\top)\boldsymbol\zeta_{\ell+1}\\
    &=-r_1\gamma_{2,1}(\bmu_1^\top)^{\otimes 2}\boldsymbol\zeta_2\\
    &=-r_1(r_1^2-2r_2)z_2(\bmu_1^\top)^{\otimes 2}\vect\bI_{d+1}\\
    &=-r_1(r_1^2-2r_2)z_2m_1,
\end{align*}
which implies that the coefficient of $h^{-6}$ is $\omega_d^{-1}2r_1(2r_2-r_1^2)z_2m_1$.

Regarding the coefficient of $h^{-8}$, proceeding as before it is possible to show that $$\omega_d^{-1}\int_{\mathbb{S}^d}\tilde D_4(\bx)\,\sigma_d(\rd\bx)=r_4(\bmu_4^\top\boldsymbol\zeta_4-z_4).$$ But using Lemma \ref{lem:moments} we have $\boldsymbol\zeta_{4}=z_4\boldsymbol{\mathcal{S}}_{d+1,4}(\vect\bI_{d+1})^{\otimes2}$, which implies $\bmu_4^\top\boldsymbol\zeta_4=z_4$ by reasoning as in \eqref{eq:mu2zeta2} (because $\bx^{\top\otimes 4}(\vect\bI_{d+1})^{\otimes 2}=(\bx^{\top\otimes 2}\vect\bI_{d+1})^2=(\bx^\top\bx)^2$), so that the integral of $\tilde D_4(\bx)$ vanishes. On the other hand, since $\tilde D_1(\bx)\tilde D_3(\bx)=-r_1\sum_{\ell=0}^3\gamma_{3,\ell}(\bmu_1^\top\otimes\bmu_\ell^\top)\bx^{\ell+1}$, then only the terms corresponding to $\ell=1$ and $\ell=3$ are nonnull after integration, so
\begin{align*}
     \omega_d^{-1}\int_{\mathbb{S}^d}\tilde D_1(\bx)\tilde D_3(\bx)\,\sigma_d(\rd\bx)
    &=-r_1\big\{\gamma_{3,1}(\bmu_1^\top)^{\otimes2}\boldsymbol\zeta_2+\gamma_{3,3}(\bmu_1^\top\otimes\bmu_3^\top)\boldsymbol\zeta_4\big\}.
\end{align*}
We have $\gamma_{3,1}=-r_1^3+r_1r_2z_2+3r_1r_2-3r_3$ and $\gamma_{3,3}=-r_3$, and we already know $(\bmu_1^\top)^{\otimes2}\boldsymbol\zeta_2=m_1z_2$, so taking into account that it can be shown that
\begin{align*}
(\bmu_1^\top\otimes\bmu_3^\top)\boldsymbol\zeta_4&=z_4(\bmu_1^\top\otimes\bmu_3^\top)\boldsymbol{\mathcal{S}}_{d+1,4}(\vect\bI_{d+1})^{\otimes2}=z_4m_1
\end{align*}
we obtain
\begin{align*}
    \omega_d^{-1}\int_{\mathbb{S}^d}\tilde D_1(\bx)\tilde D_3(\bx)\,\sigma_d(\rd\bx)&=r_1m_1\{(r_1^3-r_1r_2z_2-3r_1r_2+3r_3)z_2+r_3z_4\}.
\end{align*}
Finally,
\begin{align*}
    \omega_d^{-1}\int_{\mathbb{S}^d}\tilde D_2(\bx)^2\,\sigma_d(\rd\bx)&=\sum_{a,b=0}^2\gamma_{2,a}\gamma_{2,b}(\bmu_a^\top\otimes\bmu_b^\top)\boldsymbol\zeta_{a+b}\\
    &= \gamma_{2,0}^2+2\gamma_{2,0}\gamma_{2,2}\bmu_2^\top\boldsymbol\zeta_{2}+\gamma_{2,1}^2(\bmu_1^\top)^{\otimes2}\boldsymbol\zeta_{2}+ \gamma_{2,2}^2(\bmu_2^\top)^{\otimes2}\boldsymbol\zeta_{4}.
\end{align*}
Reasoning as before we further find $(\bmu_2^\top)^{\otimes2}\boldsymbol\zeta_{4}=\frac13(1+2m_2)z_4$, so combining this with previous calculation and simplifying we get
$$\omega_d^{-1}\int_{\mathbb{S}^d}\tilde D_2(\bx)^2\,\sigma_d(\rd\bx)=-r_2^2z_2^2+(r_1^2-2r_2)^2z_2m_1+\tfrac13r_2^2(1+2m_2)z_4.$$
Overall, grouping all the integrals, we eventually obtain
\begin{align*}
\omega_d^{-1}\int_{\mathbb{S}^d}\big[2\tilde D_0(\bx)\tilde D_4(\bx)+&2\tilde D_1(\bx)\tilde D_3(\bx)+\tilde D_2(\bx)^2\big]\,\sigma_d(\rd\bx)\\
&= 3 m_1 r_1^4 z_2
-2 m_1 r_1^2 r_2 z_2^2
-10 m_1 r_1^2 r_2 z_2 + 2 m_1 r_1 r_3\big(3z_2+z_4\big) \\
&\quad + r_2^2\big[4 m_1 z_2 - z_2^2 + \tfrac{1}{3}(1+2m_2)z_4\big].
\end{align*}

For the asymptotic expansion of $h^{-d}v_{h,d}(L)=c_{d,L^2}(h)^{-1}c_{d,L}(h)^2$ we already have the expansion of $c_{d,L}(h)$, with coefficients $c_j$. The expansion for $c_{d,L}(h)^{-1}$ can be adapted to expand $c_{d,L^2}(h)^{-1}$ by simply changing the coefficients $r_j$ for $r_j'=L_j'/L_0^2$, where $L_j'=(L^2)^{(j)}(0)/j!$; that is,
$r_j'=\sum_{k=0}^jr_kr_{j-k}$. Then, the result follows by carefully computing the coefficients in the product of the two expansions.
\end{proof}

By combining all the earlier higher-order expansions, we obtain the behavior of ${\rm MISE}(h)$ for large $h$.

\begin{corollary}\label{cor:MISEh_big}
    Assume that $L$ is $4$-times continuously differentiable at zero, with $L(0)\neq0$. We have the following expansions as $h\to\infty$:
    \begin{enumerate}[label=\textit{(\alph{*})}, ref=\textit{(\alph{*})}]
        \item If $m_1>0$, then
        $\omega_d\{{\rm MISE}(h)-{\rm MISE}2(0)\}=2r_1m_1h^{-2}+o(h^{-2}).$\label{cor1:a}
        \item If $m_1=0$, then
        $\omega_d\{{\rm MISE}(h)-{\rm MISE}2(0)\}=\big\{n^{-1}r_1^2z_2-2r_2(m_2-z_2)\big\}h^{-4}+o(h^{-4}).$\label{cor1:b}
        \item If $m_1=z_1=0$, $m_2=z_2$, $m_3=z_3=0$ and $L=L_{\rm vMF}$, then
        $
        \omega_d\{{\rm MISE}(h)-{\rm MISE}2(0)\}=n^{-1}z_2h^{-4}+\tfrac{1}{12}\big\{n^{-1}(7z_4-15z_2^2)-(m_4-z_4)\big\}h^{-8}+o(h^{-8}).
        $
        \label{cor1:c}
    \end{enumerate}
\end{corollary}

\begin{proof}
    Recall that ${\rm MISE}(h)=n^{-1}h^{-d}v_{h,d}(L)+(1-n^{-1})R_{\tilde L_h}(f)-2R_{L_h}(f)+R(f)$ and ${\rm MISE}2(0)=R(f)-\omega_d^{-1}$, so that
    \begin{align}\label{eq:MhM0}
        \omega_d\{{\rm MISE}(h)-{\rm MISE}2(0)\}=1+\omega_d\big\{n^{-1}h^{-d}v_{h,d}(L)+(1-n^{-1})R_{\tilde L_h}(f)-2R_{L_h}(f)\big\}.
    \end{align}
    When $m_1>0$, Lemma \ref{lem:highorder} gives $\omega_dh^{-d}v_{h,d}=1+O(h^{-4})$, $\omega_dR_{\tilde L_h}(f)=1+O(h^{-4})$ and $\omega_dR_{L_h}(f)=1+a_1h^{-2}+O(h^{-4})$, leading to $\omega_d\{{\rm MISE}(h)-{\rm MISE}2(0)\}=-2a_1h^{-2}+o(h^{-2})$. Then, part \ref{cor1:a} follows by noting that $a_1=\gamma_{1,1}m_1=-r_1m_1$.

    When $m_1=0$, the coefficient of $h^{-2}$ in the previous expansion vanishes, so we must obtain the coefficient of $h^{-4}$. Lemma \ref{lem:highorder} gives $\omega_dh^{-d}v_{h,d}=1+r_1^2z_2h^{-4}+O(h^{-6})$, $\omega_dR_{\tilde L_h}(f)=1+O(h^{-8})$ and $\omega_dR_{L_h}(f)=1+a_2h^{-4}+O(h^{-6})$, with $a_2=\gamma_{2,1}(m_1-z_1)+\gamma_{2,2}(m_2-z_2)=r_2(m_2-z_2)$ since $m_1=z_1=0$ and $\gamma_{2,2}=r_2$. Substitution in \eqref{eq:MhM0} yields part \ref{cor1:b}.

    Finally, if $m_j=z_j$ for $j=1,2,3$ and $L(t)=e^{-t}$ we have $\omega_dh^{-d}v_{h_d}(L)=1+r_1^2z_2h^{-4}+\frac1{12}(7z_4-15z_2^2)h^{-8}+o(h^{-8})$, $\omega_dR_{\tilde L_h}(f)=1+o(h^{-8})$ and $\omega_dR_{L_h}(f)=1+\frac{1}{4!}(m_4-z_4)h^{-8}+o(h^{-8})$, so again the desired formula in part \ref{cor1:c} is obtained by using \eqref{eq:MhM0}.
\end{proof}

In what follows, we use the notation $A_h(\bx,\by)=c_{d,A}(h)A\lrp{\frac{1-\bx^\top\by}{h^2}}$ to denote the normalization of an arbitrary function $A:\R_{\geq 0}\to\R$ not necessarily nonnegative. Also, $\lceil a\rceil$ represents
the ceiling function of a real number $a$; that is, the smallest integer that is greater than or equal to $a$.

\begin{proposition}[Taylor expansion of order $q$ of the convolution] \label{prp:tay}
Let $A:\R_{\geq 0}\to\R$ be an arbitrary function. Assume that $\bar f$ is $q$-times continuously differentiable. Then, for $\bx\in\mathbb{S}^d$,
\begin{align*}
    \int_{\Sd} A\lrp{\frac{1-\bx^\top\by}{h^2}} \bar{f}(\by)\,\sigmad(\rd \by)&=h^d\sum_{\ell=0}^{q} h^\ell\sum_{j=\lceil \ell / 2 \rceil}^\ell
    \frac{(-1)^{\ell-j}}{(\ell-j)!(2j-\ell)!}m_{d,\ell,j}(A,h)\\
    &\quad\times
    \Do{j} \bar{f}(\bx)^\top
    \left\{\bx^{\otimes \ell-j}\otimes \left(\bB_{\bx}^{\otimes 2j-\ell}\bzeta_{d-1,2j-\ell}\right)\right\} +o(h^{d+q})
\end{align*}
as $h\to 0$, where
\begin{align*}
    m_{d,\ell,j}(A,h)\defin\om{d-1}\int_0^{2h^{-2}}A(s)s^{\ell/2+d/2-1}(2-sh^2)^{j-\ell/2+d/2-1}\,\rd s.
\end{align*}
\end{proposition}

\begin{proof}[Proof of Proposition \ref{prp:tay}]
    Define $\ba_{\bx,\bxi}=\ba_{\bx,\bxi,1}+\ba_{\bx,\bxi,2}$, with $\ba_{\bx,\bxi,1}=-sh^2\bx$ and $\ba_{\bx,\bxi,2}=h\{s(2-sh^2)\}^{1/2}\bB_{\bx} \bxi$, where $\bx\in\mathbb{S}^d$, $\bxi\in\mathbb{S}^{d-1}$ and $\bB_{\bx}\in\mathcal M_{(d+1)\times d}$ is such that $\bB_{\bx}\bB_{\bx}^\top=\bI_{d+1}-\bx\bx^\top$ and $\bB_{\bx}^\top\bB_{\bx}=\bI_d$. Note that $\ba_{\bx,\bxi,1}=O(h^2)$ and $\ba_{\bx,\bxi,2}=O(h)$, so that $\ba_{\bx,\bxi}=O(h)$.

    Using the Kronecker binomial expansion $({\bf a}+{\bf b})^{\otimes j}=\boldsymbol{\mathcal{S}}_{d+1,j}\sum_{k=0}^j\binom{j}{k}{\bf a}^{\otimes k}\otimes {\bf b}^{\otimes j-k}$ for two vectors ${\bf a,b}\in\mathbb R^{d+1}$, where $\boldsymbol{\mathcal{S}}_{d+1,j}\in\mathcal{M}_{(d+1)^j\times (d+1)^j}$ is the symmetrizer matrix \citep[see][]{Holmquist1996}, and the fact that $\boldsymbol{\mathcal{S}}_{d+1,j}\Do{j}=\Do{j}$ then, for a $q$-times continuously differentiable function $\bar{f}$, Taylor's theorem leads to
\begin{align*}
    \bar{f}(\bx + \ba_{\bx, \bxi})
    &= \sum_{j=0}^q \frac{1}{j!} \Do{j} \bar{f}(\bx)^\top \ba_{\bx, \bxi}^{\otimes j}
    + o\left(\|\ba_{\bx, \bxi}\|^q\right) \\
    &= \sum_{j=0}^q \sum_{k=0}^j \frac{1}{j!} \binom{j}{k}
    \Do{j} \bar{f}(\bx)^\top
    \left(\ba_{\bx, \bxi, 1}^{\otimes k} \otimes \ba_{\bx, \bxi, 2}^{\otimes j-k}\right)
    + o\left(\|\ba_{\bx, \bxi}\|^q\right).
\end{align*}
    The product of the two $\ba$'s is of order $O(h^{k+j})$. We want to arrange the previous expression according to increasing powers of $h$, so we make a change of summation indices from $(j,k)$ to $(j,\ell)$, where $\ell=j+k$. This transforms $\sum_{j=0}^q\sum_{k=0}^j$ into $\sum_{j=0}^q\sum_{\ell=j}^{2j}$, because $0\leq k\leq j$ is equivalent to $j\leq\ell\leq 2j$. Then, we swap the summation order by noting the lower and upper limits of $j$ for a given $\ell$, resulting in
\begin{align*}
    \bar{f}(\bx + \ba_{\bx, \bxi})
    &= \sum_{\ell=0}^{2q} \sum_{j=\lceil \ell / 2 \rceil}^\ell
    \frac{1}{j!} \binom{j}{\ell - j}
    \Do{j} \bar{f}(\bx)^\top
    \left(\ba_{\bx, \bxi, 1}^{\otimes \ell - j} \otimes \ba_{\bx, \bxi, 2}^{\otimes 2j - \ell}\right)
    + o\left(\|\ba_{\bx, \bxi}\|^q\right)\\
    &=\sum_{\ell=0}^{q} h^\ell\sum_{j=\lceil \ell / 2 \rceil}^\ell
    \frac{1}{(\ell-j)!(2j-\ell)!}
    \Do{j} \bar{f}(\bx)^\top\\
    &\quad\times
    \left[\left\{(-s)^{\ell-j}\bx^{\otimes \ell-j}\right\}\otimes \left\{ s^{j-\ell/2}(2-sh^2)^{j-\ell/2}\bB_{\bx}^{\otimes 2j-\ell}\bxi^{\otimes 2j - \ell}\right\}\right]
    + o(h^q)\\
    &=\sum_{\ell=0}^{q} h^\ell s^{\ell/2}\sum_{j=\lceil \ell / 2 \rceil}^\ell
    \frac{(-1)^{\ell-j}}{(\ell-j)!(2j-\ell)!}
    \Do{j} \bar{f}(\bx)^\top\\
    &\quad\times
    (2-sh^2)^{j-\ell/2}\left\{\bx^{\otimes \ell-j}\otimes \left(\bB_{\bx}^{\otimes 2j-\ell}\bxi^{\otimes 2j - \ell}\right)\right\}
    + o(h^q),
\end{align*}
where we have truncated the series to terms of order bigger than or equal to $h^q$, leaving the remainder as $o(h^q)$.

Once we have the expansion for $f$, employ the same changes of variables as those leading to representation \eqref{eq:Lhf_int} in the proof of Lemma \ref{lem:limh0}
to get
\begin{align*}
    \int_{\mathbb{S}^d}&A\lp{\frac{1-\bx^\top\by}{h^2}}\rp \bar{f}(\by)\,\sigmad(\rd \by)\\
    &=\int_{-1}^1\int_{\mathbb{S}^{d-1}}A\lp{\frac{1-t}{h^2}}\rp \bar{f}\lrp{t\bx+\left(1-t^2\right)^{1/2} \bB_{\bx} \bxi}\left(1-t^2\right)^{d/2-1}\,\sigma_{d-1}(\rd \bxi)\rd t\\
    &=\int_0^{2h^{-2}}\int_{\mathbb{S}^{d-1}}A(s)h^ds^{d/2-1}(2-sh^2)^{d/2-1}\bar{f}(\bx+\ba_{\bx,\bxi})\,\sigma_{d-1}(\rd\bxi)\,\rd s\\
    &=\int_0^{2h^{-2}}\int_{\mathbb{S}^{d-1}}A(s)h^ds^{d/2-1}(2-sh^2)^{d/2-1}\times \sum_{\ell=0}^{q} h^\ell s^{\ell/2}\sum_{j=\lceil \ell / 2 \rceil}^\ell
    \frac{(-1)^{\ell-j}}{(\ell-j)!(2j-\ell)!}
    \Do{j} \bar{f}(\bx)^\top\\
    &\quad\times
    (2-sh^2)^{j-\ell/2}\left\{\bx^{\otimes \ell-j}\otimes \left(\bB_{\bx}^{\otimes 2j-\ell}\bxi^{\otimes 2j - \ell}\right)\right\} \,\sigma_{d-1}(\rd\bxi)\,\rd s+o(h^{d+q})\\
    &=h^d\sum_{\ell=0}^{q} h^\ell\sum_{j=\lceil \ell / 2 \rceil}^\ell
    \frac{(-1)^{\ell-j}}{(\ell-j)!(2j-\ell)!}\om{d-1}\int_0^{2h^{-2}}A(s)s^{\ell/2+d/2-1}(2-sh^2)^{j-\ell/2+d/2-1}\,\rd s\\
    &\quad\times\Do{j} \bar{f}(\bx)^\top
    \left\{\bx^{\otimes \ell-j}\otimes \left(\bB_{\bx}^{\otimes 2j-\ell}\bzeta_{d-1,2j-\ell}\right)\right\} +o(h^{d+q}),
\end{align*}
where in the last expression we use the notation from Lemma \ref{lem:moments}.
\end{proof}

\color{black}

\begin{corollary}[Taylor expansion of order two of the convolution]\label{cor:ord4}
Assume that $\bar f$ satisfies condition \ref{A1b} and let $A:\R_{\geq 0}\to\R$ be an arbitrary function satisfying condition \ref{A2b}. Then, for $\bx\in\mathbb{S}^d$,
\begin{align*}
    (A_h*\bar f)(\bx)=\int_{\Sd} A_h(\bx,\by)\bar f(\by)\,\sigmad(\rd \by)&=\bar f(\bx)+\frac{1}{2d}\widetilde m_{d,2}(A,h)\, {\rm tr}\,\bHcal\bar{f}(\bx)h^2 +O(h^4)
\end{align*}
as $h\to0$, where $\widetilde m_{d,2\ell}(A,h)=m_{d,2\ell,2\ell}(A,h)/\lambda_{h,d}(A)$ and $\bHcal\bar{f}$ is the Hessian matrix of $\bar f$.
\end{corollary}

\begin{proof}[Proof of Corollary \ref{cor:ord4}]
Note that $c_{d,A}(h)=h^{-d}\lambda_{h,d}(A)^{-1}$, so that Proposition \ref{prp:tay} gives
\begin{align*}
    \int_{\Sd} A_h(\bx,\by)\bar f(\by)\,\sigmad(\rd \by)&=c_{d,A}(h)\int_{\Sd} A\lrp{\frac{1-\bx^\top\by}{h^2}} \bar{f}(\by)\,\sigmad(\rd \by)\\
    &=\lambda_{h,d}(A)^{-1}\sum_{\ell=0}^{4} h^\ell\sum_{j=\lceil \ell / 2 \rceil}^\ell
    \frac{(-1)^{\ell-j}}{(\ell-j)!(2j-\ell)!} m_{d,\ell,j}(A,h)\\
    &\quad\times
    \Do{j} \bar{f}(\bx)^\top
    \left\{\bx^{\otimes \ell-j}\otimes \left(\bB_{\bx}^{\otimes 2j-\ell}\bzeta_{d-1,2j-\ell}\right)\right\} +o(h^4).
\end{align*}

For $\ell=0$ the coefficient associated to $h^0$ is
\begin{align*}
 \lambda_{h,d}(A)^{-1}m_{d,0,0}(A,h) \bar{f}(\bx)=\bar{f}(\bx).
\end{align*}

For $\ell=1$ the coefficient associated to $h^1$ is
\begin{align*}
\lambda_{h,d}(A)^{-1}m_{d,1,1}(A,h)\Do{1} \bar{f}(\bx)^\top
\bB_{\bx}\bzeta_{d-1,1}=0.
\end{align*}

For $\ell=2$ the coefficient associated to $h^2$ is
\begin{align}
\lambda_{h,d}&(A)^{-1}\sum_{j=1}^2
    \frac{(-1)^{2-j}}{(2-j)!(2j-2)!} m_{d,2,j}(A,h) \Do{j} \bar{f}(\bx)^\top
    \left\{\bx^{\otimes 2-j}\otimes \left(\bB_{\bx}^{\otimes 2j-2}\bzeta_{d-1,2j-2}\right)\right\}\nonumber\\
   &=-\lambda_{h,d}(A)^{-1}m_{d,2,1}(A,h)\Do{1} \bar{f}(\bx)^\top\bx +\frac{1}{2} \lambda_{h,d}(A)^{-1}m_{d,2,2}(A,h)\Do{2} \bar{f}(\bx)^\top
   \bB_{\bx}^{\otimes 2}\bzeta_{d-1,2}\nonumber\\
   &=\frac{1}{2d}\widetilde m_{d,2}(A,h)\Do{2} \bar{f}(\bx)^\top \{\vect\bI_{d+1}- \bx^{\otimes 2}\}\label{eq:coefl21}\\
   &=\frac{1}{2d}\widetilde m_{d,2}(A,h) \tr{\bHcal\bar{f}(\bx)}\label{eq:coefl22}.
\end{align}
In \eqref{eq:coefl21}, we use the fact that $\bx^\top\D\bar f(\bx)=0$ (due to Euler's homogeneous function theorem), Lemma \ref{lem:moments}, and that
\begin{align*}
   \bB_{\bx}^{\otimes 2}\bzeta_{d-1,2}=\frac{1}{d}(\bB_{\bx}\otimes\bB_{\bx})\vect\bI_{d}=\frac{1}{d}\vect(\bB_{\bx}\bB_{\bx}^\top)=\frac{1}{d} \{\vect\bI_{d+1}- \bx^{\otimes 2}\}.
\end{align*}
In \eqref{eq:coefl22}, we use $\Do{2} \bar{f}(\bx)^\top \vect\bI_{d+1}=\tr{\bHcal\bar{f}(\bx)}$ and $\Do{2} \bar{f}(\bx)^\top \bx^{\otimes 2}=\bx^\top \bHcal \bar{f}(\bx) \bx=0$ (the last equality, again, due to Euler's homogeneous function theorem).

For $\ell=3$ the coefficient associated to $h^3$ is
\begin{align*}
 \lambda_{h,d}(A)^{-1}\sum_{j=2}^3&
    \frac{(-1)^{3-j}}{(3-j)!(2j-3)!}m_{d,3,j}(A,h)\Do{j} \bar{f}(\bx)^\top
    \left\{\bx^{\otimes 3-j}\otimes \left(\bB_{\bx}^{\otimes 2j-3}\bzeta_{d-1,2j-3}\right)\right\}\\
    &=-\lambda_{h,d}(A)^{-1}m_{d,3,2}(A,h)\Do{2} \bar{f}(\bx)^\top
    \left\{\bx\otimes \left(\bB_{\bx}\bzeta_{d-1,1}\right)\right\}\\
    &\quad+\frac{1}{6}\lambda_{h,d}(A)^{-1}m_{d,3,3}(A,h)\Do{3} \bar{f}(\bx)^\top
     \left(\bB_{\bx}^{\otimes 3}\bzeta_{d-1,3}\right)\\
     &=0,
\end{align*}
because $\bzeta_{d-1,1}=\bzeta_{d-1,3}=\zero$ by Lemma \ref{lem:moments}.

Finally, for $\ell=4$ the coefficient associated to $h^4$ is
\begin{align*}
\lambda_{h,d}(A)^{-1}\sum_{j=2}^4&
    \frac{(-1)^{4-j}}{(4-j)!(2j-4)!}m_{d,4,j}(A,h) \Do{j} \bar{f}(\bx)^\top
    \left\{\bx^{\otimes 4-j}\otimes \left(\bB_{\bx}^{\otimes 2j-4}\bzeta_{d-1,2j-4}\right)\right\}\\
    &=\frac{1}{2}\lambda_{h,d}(A)^{-1}m_{d,4,2}(A,h)\Do{2} \bar{f}(\bx)^\top\bx^{\otimes 2}\\
    &\quad- \frac{1}{2}\lambda_{h,d}(A)^{-1}m_{d,4,3}(A,h)\Do{3} \bar{f}(\bx)^\top
 ( \bx\otimes \left(\bB_{\bx}^{\otimes 2}\bzeta_{d-1,2})\right)\\
 &\quad+ \frac{1}{24}\lambda_{h,d}(A)^{-1}m_{d,4,4}(A,h)\Do{4
   } \bar{f}(\bx)^\top
     \left(\bB_{\bx}^{\otimes 4}\bzeta_{d-1,4}\right)\\
     &=- \frac{1}{2}\lambda_{h,d}(A)^{-1}m_{d,4,3}(A,h)\Do{3} \bar{f}(\bx)^\top
 ( \bx\otimes \left(\bB_{\bx}^{\otimes 2}\bzeta_{d-1,2})\right)\\
 &\quad+ \frac{1}{24}\lambda_{h,d}(A)^{-1}m_{d,4,4}(A,h)\Do{4
   } \bar{f}(\bx)^\top
     \left(\bB_{\bx}^{\otimes 4}\bzeta_{d-1,4}\right).
\end{align*}
This last term is included in the corollary statement as $O(h^4)$.
\end{proof}

Define $(A_h*B_h)(\bx,\by)\defin\{A_h*B_h(\by,\cdot)\}(\bx)=\int_{\Sd}A_h(\bx,\bz)B_h(\by,\bz)\,\sigma_d(\rd\bz)$ for arbitrary kernels $A_h$ and $B_h$.

\begin{corollary}[Taylor expansion for the convolution of kernels]\label{cor:ord4_2}
Assume that $\bar f$ satisfies condition \ref{A1b} and let $A,B:\R_{\geq 0}\to\R$ be arbitrary functions satisfying condition \ref{A2b}. Then, for $\bx\in\mathbb{S}^d$,
\begin{align*}
    \int_{\Sd} (A_h*B_h)(\bx,\by)\bar f(\by)\,\sigmad(\rd \by)&=\bar f(\bx)+\frac{h^2}{2d}\lrb{\widetilde m_{d,2}(A,h)+\widetilde m_{d,2}(B,h)}{\rm tr}\,{\bHcal\bar{f}(\bx)}+O(h^4).
\end{align*}
\end{corollary}

\begin{proof}[Proof of Corollary \ref{cor:ord4_2}]
    Using Corollary \ref{cor:ord4} first with respect to $\sigmad(\rd \by)$ and then with respect to $\sigmad(\rd \bz)$ we have
    \begin{align*}
        \int_{\Sd} (A_h*B_h)&(\bx,\by)\bar f(\by)\,\sigmad(\rd \by)\\
        &=\int_{\Sd}\int_{\Sd} A_h(\bx,\bz)B_h(\by,\bz)\bar f(\by)\,\sigmad(\rd \by)\,\sigmad(\rd \bz)\\
        &=\int_{\Sd} A_h(\bx,\bz)\lrb{\bar f(\bz)+\frac{h^2}{2d}\widetilde m_{d,2}(B,h){\rm tr}\,\bHcal\bar f(\bz)+O(h^4)}\,\sigmad(\rd \bz)\\
        &=\bar f(\bx)+\frac{h^2}{2d}\widetilde m_{d,2}(A,h){\rm tr}\,\bHcal\bar f(\bx)+O(h^4)+\frac{h^2}{2d}\widetilde m_{d,2}(B,h){\rm tr}\,\bHcal\bar f(\bx)+O(h^4).
    \end{align*}
\end{proof}
Notice that the above result is completely analogous to Corollary \ref{cor:ord4}, with $\widetilde m_{d,2}(L,h)$ replaced by $\widetilde m_{d,2}(A,h)+\widetilde m_{d,2}(B,h)$.

\begin{corollary}\label{cor:RLhf}
Assume that $\bar f$ satisfies condition \ref{A1b} and let $A:\R_{\geq 0}\to\R$ be an arbitrary function satisfying condition \ref{A2b}. Then,
\begin{align*}
    R_{A_h}(\bar f)=R(f)+\frac{1}{2d}\widetilde m_{d,2}(A,h) \psi_1h^2 +O(h^4)
\end{align*}
as $h\to0$, where $\psi_1=\int_{\mathbb{S}^d}\bar f(\bx){\rm tr}\,\bHcal\bar f(\bx)\,\sigma_d(\rd\bx).$
\end{corollary}

\begin{proof}[Proof of Corollary \ref{cor:RLhf}]
    Apply Corollary \ref{cor:ord4} to $R_{A_h}(\bar f)=\int_{\mathbb{S}^d}(A_h*\bar f)(\bx)\bar f(\bx)\,\sigma_d(\rd\bx)$.
\end{proof}

\begin{corollary}\label{cor:AhBh3}
Let $A,B:\R_{\geq 0}\to\R$ be continuous functions satisfying condition \ref{A2}.  If $f$ is square integrable, then
\begin{align*}      \int_{\Sd}\int_{\Sd}A_h(\bx,\by)B_h(\bx,\by)\bar f(\bx)\bar f(\by)\,\sigma_d(\rd \bx)\,\sigma_d(\rd \by)
&=\frac{\lambda_d(AB)}{\lambda_d(A)\lambda_d(B)}R(f) h^{-d}\{1+o(1)\}
\end{align*}
as $h\to0$.
\end{corollary}

\begin{proof}[Proof of Corollary \ref{cor:AhBh3}]
    Note that
    \begin{align*}
       A_h(\bx,\by)B_h(\bx,\by)&=c_{d,A}(h)c_{d,B}(h)(AB)\lrp{\frac{1-\bx^\top \by}{h^2}}=\frac{c_{d,A}(h)c_{d,B}(h)}{c_{d,AB}(h)}(AB)_h(\bx,\by).
    \end{align*}
    Then, use Lemma \ref{lem:limh0}.
\end{proof}

\begin{proposition}
\label{prop:AhBhCh}
Let $\bar{f}:\mathbb{R}^{d+1}\backslash\{\zero\}\to\R$ be continuous and $A,B,C:\R_{\geq 0}\to\R$ be arbitrary functions satisfying condition \ref{A2}. Under \ref{A3}, it follows that
\begin{align*}
\int_{\Sd}\int_{\Sd}&(A_h*B_h)(\bx,\by)C_h(\bx,\by)\bar f(\bx)\bar f(\by)\,\sigma_d(\rd \bx)\,\sigma_d(\rd \by)\\
\sim&\,h^{2d} c_{d,A}(h)c_{d,B}(h)c_{d,C}(h)R(f)  \gamma_d \int_0^{\infty} A(r) r^{d/2-1}\left\{\int_0^{\infty} \rho^{d/2-1} C(\rho) \varphi_d(B, r, \rho)\,\rd \rho\right\}\,\rd r, \\
\varphi_d(L, r, \rho)= & \begin{cases}L\left(r+\rho-2(r \rho)^{1/2}\right)+L\left(r+\rho+2(r \rho)^{1/2}\right), & d=1, \\
\int_{-1}^1\left(1-\theta^2\right)^{(d-3)/2} L\left(r+\rho-2 \theta(r \rho)^{1/2}\right)\,\rd \theta, & d \geq 2,\end{cases} \\
\gamma_d= & \begin{cases}1, & d=1, \\
\om{d-1} \om{d-2} 2^{d-2}, & d \geq 2 .\end{cases}
\end{align*}
\end{proposition}

\begin{proof}[Proof of Proposition \ref{prop:AhBhCh}]
To begin, the integral can be expressed as:
   \begin{align}
\int_{\Sd}\int_{\Sd}&(A_h*B_h)(\bx,\by)C_h(\bx,\by)\bar f(\bx)\bar f(\by)\,\sigma_d(\rd \bx)\,\sigma_d(\rd \by)\nonumber\\
=&\;\int_{\Sd}\int_{\Sd}A_h(\bx,\bz)\left\{\int_{\Sd}C_h(\by,\bx)B_h(\by,\bz)\bar f(\by)\,\sigma_d(\rd \by)\right\}\bar f(\bx)\,\sigma_d(\rd \bx)\,\sigma_d(\rd \bz)\nonumber\\
=&\;c_{d,A}(h)c_{d,B}(h)c_{d,C}(h)\nonumber\\
&\times \int_{\Sd}\int_{\Sd}A\lrp{\frac{1-\bx^\top\bz}{h^2}}\left\{\int_{\Sd}C\lrp{\frac{1-\by^\top\bx}{h^2}}B\lrp{\frac{1-\by^\top\bz}{h^2}}\bar f(\by)\,\sigma_d(\rd \by)\right\}\nonumber\\
&\times \bar f(\bx)\,\sigma_d(\rd \bx)\,\sigma_d(\rd \bz).\label{eq:AhBhCh1}
\end{align}

The computation of \eqref{eq:AhBhCh1} will be divided into the cases $d \geq 2$ and $d=1$. To begin with, suppose that $d \geq 2$. Consider the change of variables:
$$
\bz=s \bx+\left(1-s^2\right)^{1/2} \bB_{\bx} \bxi, \quad \sigma_d(\rd \bz)=\left(1-s^2\right)^{d/2-1} \,\rd s \,\sigma_{d-1}(\rd \bxi)
$$
where $s \in(-1,1), \bxi \in \Sdm$ and $\bB_{\bx}=\left(\mathbf{b}_1, \ldots, \mathbf{b}_d\right)_{(d+1) \times d}$ is the semi-orthonormal matrix $\left(\bB_{\bx}^\top \bB_{\bx}=\bI_d\right.$ and $\left.\bB_{\bx} \bB_{\bx}^\top=\bI_{d+1}-\bx \bx^\top\right)$ resulting from the completion of $\bx$ to the orthonormal basis $\left\{\bx, \mathbf{b}_1, \ldots, \mathbf{b}_d\right\}$ of $\mathbb{R}^{d+1}$. Here $\bI_d$ represents the identity matrix with dimension $d$. Consider also the change of variable
$$
\by=t \bx+\tau \bB_{\bx} \bxi+\left(1-t^2-\tau^2\right)^{1/2} \bA_{\xib} \etab, \quad \sigma_d(\rd \by)=\left(1-t^2-\tau^2\right)^{\frac{d-3}{2}} \,\rd t \,\rd  \tau \,\sigma_{d-2}(\rd \etab)
$$
where $t, \tau \in(-1,1), t^2+\tau^2<1, \etab \in \Sdmm$ and $\bA_{\bxi}=\left(\mathbf{a}_1, \ldots, \mathbf{a}_d\right)_{(d+1) \times(d-1)}$ is the semi-orthonormal matrix $\left(\bA_{\bxi}^\top \bA_{\bxi}=\bI_d\right.$ and $\left.\bA_{\bxi} \bA_{\bxi}^\top=\bI_{d+1}-\bx \bx^\top-\bB_{\bx} \bxi \bxi^\top \bB_{\bx}^\top\right)$ resulting from the completion of $\left\{\bx, \bB_{\bx} \bxi\right\}$ to the orthonormal basis $\left\{\bx, \bB_{\bx} \bxi, \mathbf{a}_1, \ldots, \mathbf{a}_{d-1}\right\}$ of $\mathbb{R}^{d+1}$. With these two changes of variables,
$$
\by^\top \bz=s t+\tau\left(1-s^2\right)^{1/2}, \quad \bx^\top\left(\bB_{\bx} \bxi\right)=\bx^\top\left(\bA_{\bxi} \etab\right)=\left(\bB_{\bx} \bxi\right)^\top\left(\bA_{\bxi} \etab\right)=0,
$$
and therefore,
\begin{align}
\eqref{eq:AhBhCh1}=&\;c_{d,A}(h)c_{d,B}(h)c_{d,C}(h)\nonumber\\
&\times \int_{\Sd}\int_{\Sdm}\int_{-1}^1A\lrp{\frac{1-s}{h^2}}\Bigg\{\int_{\Sdmm}\int\int_{t^2+\tau^2<1}C\lrp{\frac{1-t}{h^2}}B\lrp{\frac{1-st-\tau(1-s^2)^{1/2}}{h^2}}\nonumber\\
&\times \bar f(t\bx+\tau\bB_{\bx}\bxi+(1-t^2-\tau^2)^{1/2}\bA_{\bxi}\etab)(1-t^2-\tau^2)^{(d-3)/2}\,\rd t\,\rd \tau\,\sigma_{d-2}(\rd \etab)\Bigg\}\nonumber\\
&\times (1-s^2)^{d/2-1} \bar f(\bx)\, \rd s\,\sigma_{d-1}(\rd \bxi)\,\sigma_d(\rd \bx). \label{eq:AhBhCh2}
\end{align}

Consider now the change of variables $r=\frac{1-s}{h^2}$ and then
$$
\left\{\begin{array}{l}
\rho=\frac{1-t}{h^2}, \\
\theta=\frac{\tau}{h\left[\rho\left(2-h^2 \rho\right)\right]^{1/2}},
\end{array} \quad\left|\frac{\partial(t, \tau)}{\partial(\rho, \theta)}\right|=h^3\left[\rho\left(2-h^2 \rho\right)\right]^{1/2}\right.
$$

With these changes of variables, $\tau=h \theta\left[\rho\left(2-h^2 \rho\right)\right]^{1/2}, t=1-h^2 \rho$ and, as a result:
\begin{align*}
1-s^2 & =h^2 r\left(2-h^2 r\right) \\
1-t^2 & =h^2 \rho\left(2-h^2 \rho\right) \\
1-t^2-\tau^2 & =\left(1-\theta^2\right) h^2 \rho\left(2-h^2 \rho\right), \\
\frac{1-s t-\tau\left(1-s^2\right)^{1/2}}{h^2} & =r+\rho-h^2 r \rho-\theta\left[r \rho\left(2-h^2 r\right)\left(2-h^2 \rho\right)\right]^{1/2}.
\end{align*}

Then:
\begin{align}
\eqref{eq:AhBhCh2}=&\;c_{d,A}(h)c_{d,B}(h)c_{d,C}(h)\nonumber\\
&\times \int_{\Sd}\int_{\Sdm}\int_{0}^{2h^{-2}}A\lrp{r}\Bigg\{\int_{\Sdmm}\int_{0}^{2h^{-2}}\int_{-1}^1C\lrp{\rho}\nonumber\\
&\times B\lrp{r+\rho-h^2 r \rho-\theta\left[r \rho\left(2-h^2 r\right)\left(2-h^2 \rho\right)\right]^{1/2}}\nonumber\\
&\times \bar f((1-h^2 \rho) \bx+h[\rho(2-h^2 \rho)]^{1/2}[\theta \bB_{\bx} \bxi+(1-\theta^2)^{1/2} \bA_{\bxi} \etab])\nonumber\\
&\times(1-\theta^2)^{\frac{d-3}{2}} h^{d-3}[\rho(2-h^2 \rho)]^{(d-3)/2} h^3[\rho(2-h^2 \rho)]^{1/2}\,\rd \theta\,\rd \rho\,\sigma_{d-2}(\rd \etab)\Bigg\}\nonumber\\
&\times h^{d-2} r^{d/2-1}(2-h^2 r)^{d/2-1} h^2 \bar f(\bx)\, \rd r\,\sigma_{d-1}(\rd \bxi)\,\sigma_d(\rd \bx). \label{eq:AhBhCh3}
\end{align}

Using the Dominated Convergence Theorem (DCT), it follows that
\begin{align}
\eqref{eq:AhBhCh3}\sim&\,h^{2d}c_{d,A}(h)c_{d,B}(h)c_{d,C}(h)\nonumber\\
&\times \int_{\Sd}\int_{\Sdm}\int_{0}^{\infty}A\lrp{r}\Bigg\{\int_{\Sdmm}\int_{0}^{\infty}\int_{-1}^1C\lrp{\rho} B\lrp{r+\rho-2\theta(r \rho)^{1/2}}\bar f(\bx)\nonumber\\
&\times(1-\theta^2)^{(d-3)/2}(2\rho)^{d/2-1}\,\rd \theta\,\rd \rho\,\sigma_{d-2}(\rd \etab)\Bigg\}(2r)^{d/2-1} \bar f(\bx)\, \rd r\,\sigma_{d-1}(\rd \bxi)\,\sigma_d(\rd \bx)\nonumber\\
=&\;h^{2d}c_{d,A}(h)c_{d,B}(h)c_{d,C}(h)\om{d-1}\om{d-2}R(f)2^{d-2}\nonumber\\
&\times \int_{0}^{\infty}A\lrp{r}\Bigg\{\int_{0}^{\infty}\int_{-1}^1C\lrp{\rho} B\lrp{r+\rho-2\theta(r \rho)^{1/2}}\nonumber\\
&\times(1-\theta^2)^{(d-3)/2}\rho^{d/2-1}\,\rd \theta\,\rd \rho\,\Bigg\}r^{d/2-1} \, \rd r.\label{eq:AhBhCh4}
\end{align}

For the calculation of \eqref{eq:AhBhCh1} in the case
$d=1$, similar steps can be followed. To begin with, consider the change of variable:
$$
\bz=s \bx+\left(1-s^2\right)^{1/2} \bB_{\bx} \xi, \quad \sigma_1(\rd \bz)=\left(1-s^2\right)^{d/2-1}\, \rd s \sigma_{0}(\rd \xi)
$$
where $s \in(-1,1), \xi \in \{-1,1\}$ and $\bB_{\bx}=\left(\mathbf{b}_1\right)_{2\times 1}$ is the vector resulting from the completion of $\bx$ to the orthonormal basis $\left\{\bx, \mathbf{b}_1\right\}$ of $\mathbb{R}^{2}$. Consider also the change of variable:
$$
\by=t \bx+\left(1-t^2\right)^{1/2} \tilde{\bB}_{\bx} \eta, \quad \sigma_1(\rd \by)=\left(1-t^2\right)^{d/2-1} \rd t \sigma_{0}(\rd \eta)
$$
where $t \in(-1,1), \eta \in \{-1,1\}$ and $\tilde{\bB}_{\bx}=\left(\tilde{\mathbf{b}}_1\right)_{2 \times 1}$ is the vector resulting from the completion of $\bx$ to the orthonormal basis $\left\{\bx, \tilde{\mathbf{b}}_1\right\}$ of $\mathbb{R}^{2}$. Notice that $\bB_{\bx} \xi=\tilde{\bB}_{\bx} \eta$ or  $\bB_{\bx} \xi=-\tilde{\bB}_{\bx} \eta$. Therefore:
\begin{align}
\eqref{eq:AhBhCh1}=&\;c_{1,A}(h)c_{1,B}(h)c_{1,C}(h)\int_{\Sdo}\int_{\{-1,1\}}\int_{-1}^1A\lrp{\frac{1-s}{h^2}}\Bigg\{\int_{\{-1,1\}}\int_{-1}^1 C\lrp{\frac{1-t}{h^2}}\nonumber\\
&\times B\lrp{\frac{1-st-(1-t^2)^{1/2}(1-s^2)^{1/2}(\bB_{\bx} \xi)^\top(\tilde{\bB}_{\bx} \eta)}{h^2}}\nonumber\\
&\times \bar f(t \bx+\left(1-t^2\right)^{1/2} \tilde{\bB}_{\bx} \eta)(1-t^2)^{-1/2}\,\rd t\,\sigma_{0}(\rd \eta)\Bigg\}\nonumber\\
&\times (1-s^2)^{-1/2} \bar f(\bx)\, \rd s\,\sigma_{0}(\rd \xi)\,\sigma_1(\rd \bx). \label{eq:AhBhCh21}
\end{align}

Consider now the change of variables $r=\frac{1-s}{h^2}$ and $\rho=\frac{1-t}{h^2}$. With these changes of variable, $s=1-rh^2$ and $t=1-h^2 \rho$ and, as a result:
\begin{align*}
	1-s^2 & =h^2 r\left(2-h^2 r\right), \\
	1-t^2 & =h^2 \rho\left(2-h^2 \rho\right).
\end{align*}
Therefore:
\begin{align}
	\eqref{eq:AhBhCh21}=&\;h^2c_{1,A}(h)c_{1,B}(h)c_{1,C}(h) \int_{\Sdo}\int_{\{-1,1\}}\int_{0}^{2h^{-2}}A\lrp{r}\Bigg\{\int_{\{-1,1\}}\int_{0}^{2h^{-2}} C\lrp{\rho}\nonumber\\
	&\times B\lrp{r+\rho-h^2r\rho-(r\rho(2-h^2r)(2-h^2\rho))^{1/2}(\bB_{\bx} \xi)^\top(\tilde{\bB}_{\bx} \eta)}\nonumber\\
	&\times \bar f((1-h^2 \rho) \bx+\left(h^2 \rho\left(2-h^2 \rho\right)\right)^{1/2} \tilde{\bB}_{\bx} \eta)(\rho\left(2-h^2 \rho\right))^{-1/2}\,\rd \rho\,\sigma_{0}(\rd \eta)\Bigg\}\nonumber\\
	&\times ( r\left(2-h^2 r\right))^{-1/2} \bar f(\bx)\, \rd r\,\sigma_{0}(\rd \xi)\,\sigma_1(\rd \bx)\nonumber\\
	=&\;h^2c_{1,A}(h)c_{1,B}(h)c_{1,C}(h) \int_{\Sdo}\int_{\{-1,1\}}\int_{0}^{2h^{-2}}A\lrp{r}\Bigg\{\int_{0}^{2h^{-2}} C\lrp{\rho}\nonumber\\
	&\times \Big[ B\lrp{r+\rho-h^2r\rho+(r\rho(2-h^2r)(2-h^2\rho))^{1/2}} \bar f((1-h^2 \rho) \bx+\left(h^2 \rho\left(2-h^2 \rho\right)\right)^{1/2} {\bB}_{\bx} \xi)\nonumber\\
    &+ B\lrp{r+\rho-h^2r\rho-(r\rho(2-h^2r)(2-h^2\rho))^{1/2}} \bar f((1-h^2 \rho) \bx-\left(h^2 \rho\left(2-h^2 \rho\right)\right)^{1/2} {\bB}_{\bx} \xi)\Big]\nonumber\\
	&\times (\rho\left(2-h^2 \rho\right))^{-1/2}\,\rd \rho\Bigg\} ( r\left(2-h^2 r\right))^{-1/2} \bar f(\bx)\, \rd r\,\sigma_{0}(\rd \xi)\,\sigma_1(\rd \bx). \label{eq:AhBhCh22}
\end{align}

Using the DCT, it follows that
\begin{align}
	\eqref{eq:AhBhCh22}=&\;2h^2c_{1,A}(h)c_{1,B}(h)c_{1,C}(h)R(f) \int_{\Sdo}\int_{0}^{\infty}A\lrp{r}\Bigg\{\int_{0}^{\infty} C\lrp{\rho}\nonumber\\
	&\times \Big[ B\lrp{r+\rho+2(r\rho)^{1/2}} \bar f(\bx)+ B\lrp{r+\rho-2(r\rho)^{1/2}} \bar f( \bx)\Big](2\rho)^{-1/2}\,\rd \rho\Bigg\}\nonumber\\
	&\times  (2r)^{-1/2} \bar f(\bx)\, \rd r\,\sigma_{0}(\rd \xi)\,\sigma_1(\rd \bx)\nonumber\\
	=&\;h^2c_{1,A}(h)c_{1,B}(h)c_{1,C}(h) \int_{0}^{\infty}A\lrp{r}\Bigg\{\int_{0}^{\infty} C\lrp{\rho}\nonumber\\
	&\times \Big[ B\lrp{r+\rho+2(r\rho)^{1/2}} + B\lrp{r+\rho-2(r\rho)^{1/2}} \Big] \rho^{-1/2}\,\rd \rho\Bigg\} r^{-1/2} \, \rd r.
	 \label{eq:AhBhCh23}
\end{align}

Therefore, using \eqref{eq:AhBhCh4} and \eqref{eq:AhBhCh23} it follows that
\begin{align*}
\eqref{eq:AhBhCh1}=&h^{2d}c_{d,A}(h)c_{d,B}(h)c_{d,C}(h) R(f)  \gamma_d \int_0^{\infty} A(r) r^{d/2-1}\left\{\int_0^{\infty} \rho^{d/2-1} C(\rho) \varphi_d(B, r, \rho)\,\rd \rho\right\}\,\rd r, \\
\varphi_d(L, r, \rho)= & \begin{cases}L\left(r+\rho-2(r \rho)^{1/2}\right)+L\left(r+\rho+2(r \rho)^{1/2}\right), & d=1, \\
\int_{-1}^1\left(1-\theta^2\right)^{(d-3)/2} L\left(r+\rho-2 \theta(r \rho)^{1/2}\right)\,\rd \theta, & d \geq 2,\end{cases} \\
\gamma_d= & \begin{cases}1, & d=1, \\
\om{d-1} \om{d-2} 2^{d-2}, & d \geq 2 .\end{cases}
\end{align*}
\end{proof}

\color{black}
\begin{proposition}
    \label{prop:AhBhChDh}
Let $\bar{f}:\mathbb{R}^{d+1}\backslash\{\zero\}\to\R$ be continuous differentiable and $A,B,C,D:\R_{\geq 0}\to\R$ be arbitrary functions satisfying condition \ref{A2}. Under \ref{A3}, it follows that
\begin{align*}
\int_{\Sd}\int_{\Sd}&(A_h*B_h)(\bx,\by)(C_h*D_h)(\bx,\by)\bar f(\bx)\bar f(\by)\,\sigma_d(\rd \bx)\,\sigma_d(\rd \by)\\
&\sim h^{3d}c_{d,A}(h)c_{d,B}(h)c_{d,C}(h) c_{d,D}(h) R(f)  \tilde{\gamma}_d \\
&\quad\times \int_0^{\infty}\Bigg\{ \int_0^{\infty} A(r) r^{d/2-1}\varphi_d(B, r, c)\,\rd r\int_0^{\infty} C(\rho) \rho^{d/2-1}\varphi_d(D, \rho, c)\,\rd \rho\Bigg\}c^{d/2-1}\, \rd c, \\
\tilde{\gamma}_d= & \begin{cases}2^{-1/2}, & d=1, \\
\om{d-1}\om{d-2}^22^{(3d-6)/2}, & d \geq 2 .\end{cases}
\end{align*}
\end{proposition}

\begin{proof}[Proof of Proposition \ref{prop:AhBhChDh}]
First, note that
\begin{align}
    \int_{\Sd}\int_{\Sd}&(A_h*B_h)(\bx,\by)(C_h*D_h)(\bx,\by)\bar f(\bx)\bar f(\by)\,\sigma_d(\rd \bx)\,\sigma_d(\rd \by)\nonumber\\   
    =&\int_{\Sd}\int_{\Sd}\left\{\int_{\Sd}A_h(\bx,\bv)B_h(\by,\bv)\,\sigma_d(\rd \bv)\right\}\left\{\int_{\Sd}C_h(\bx,\bw)D_h(\by,\bw) \,\sigma_d(\rd \bw)\right\}\nonumber\\
    &\times \bar f(\bx)\bar f(\by)\,\sigma_d(\rd \bx)\,\sigma_d(\rd \by)\nonumber\\
    =&\;c_{d,A}(h)c_{d,B}(h)c_{d,C}(h)c_{d,D}(h)\nonumber\\
    &\times \int_{\Sd}\int_{\Sd}\left\{\int_{\Sd}A\lrp{\frac{1-\bx^\top\bv}{h^2}}B\lrp{\frac{1-\by^\top\bv}{h^2}}\,\sigma_d(\rd \bv)\right\}\nonumber\\
    &\times \left\{\int_{\Sd}C\lrp{\frac{1-\bx^\top\bw}{h^2}}D\lrp{\frac{1-\by^\top\bw}{h^2}}\,\sigma_d(\rd \bw)\right\}\bar f(\bx)\bar f(\by)\,\sigma_d(\rd \bx)\,\sigma_d(\rd \by).\label{eq:AhBhChDh1}
\end{align}

Again, the computation of \eqref{eq:AhBhChDh1} will be divided into the cases $d \geq 2$ and $d=1$. To begin with, suppose that $d \geq 2$. Consider the change of variable:
$$
\bv=s \bx+\left(1-s^2\right)^{1/2} \bB_{\bx} \bxi, \quad \sigma_d(\rd \bv)=\left(1-s^2\right)^{d/2-1}\, \rd s \sigma_{d-1}(\rd \bxi)
$$
where $s \in(-1,1), \bxi \in \Sdm$ and $\bB_{\bx}=\left(\mathbf{b}_1, \ldots, \mathbf{b}_d\right)_{(d+1) \times d}$ is the semi-orthonormal matrix resulting from the completion of $\bx$ to the orthonormal basis $\left\{\bx, \mathbf{b}_1, \ldots, \mathbf{b}_d\right\}$ of $\mathbb{R}^{d+1}$.
Therefore:
\begin{align}
\int_{\Sd}&A\lrp{\frac{1-\bx^\top\bv}{h^2}}B\lrp{\frac{1-\by^\top\bv}{h^2}}\,\sigma_d(\rd \bv)\nonumber\\
=&\;\int_{\Sdm}\int_{-1}^1A\lrp{\frac{1-s}{h^2}}B\lrp{\frac{1-s \by^\top\bx -\left(1-s^2\right)^{1/2} \by^\top(\bB_{\bx} \bxi)}{h^2}}\left(1-s^2\right)^{d/2-1}\, \rd s\, \sigma_{d-1}(\rd \bxi). \label{eq:AhBh1}
\end{align}
Consider now the change of variable $r=\frac{1-s}{h^2}$. With this change of variable, $s=1-rh^2$ and, as a result:
\begin{align*}
1-s^2 & =h^2 r\left(2-h^2 r\right).
\end{align*}
Therefore:
\begin{align}
\eqref{eq:AhBh1}=&\;h^{d}\int_{\Sdm}\int_{0}^{2h^{-2}}A\lrp{r}B\lrp{\frac{1-(1-rh^2) \by^\top\bx -\left(h^2 r\left(2-h^2 r\right)\right)^{1/2} \by^\top(\bB_{\bx} \bxi)}{h^2}}\nonumber\\
    &\times [r\left(2-h^2 r\right)]^{d/2-1}\,\rd r\, \sigma_{d-1}(\rd \bxi). \label{eq:AhBh2}
\end{align}
 Define $\bz_{\bx,\by}\defin (1-(\by^\top \bx)^2)^{-1/2}\by^\top\bB_{\bx}$ and consider the change of variable
$$
\bxi=s \bz_{\bx,\by}+\left(1-s^2\right)^{1/2} \bA_{\bx,\by} \bzeta, \quad \sigma_{d-1}(\rd \bxi)=\left(1-s^2\right)^{(d-3)/2}\, \rd s \sigma_{d-2}(\rd \bzeta)
$$
where $s \in(-1,1), \bzeta \in \Sdmm$ and $\bA_{\bx,\by}=\left(\mathbf{a}_1, \ldots, \mathbf{a}_d\right)_{(d+1) \times (d-1)}$ is the semi-orthonormal matrix. Therefore:
\begin{align}
\eqref{eq:AhBh2}=&\;h^{d}\om{d-2}\int_{-1}^1\int_{0}^{2h^{-2}}A\lrp{r} B\lrp{\frac{1-(1-rh^2) \by^\top\bx -\left(h^2 r\left(2-h^2 r\right)\right)^{1/2} s (1-(\by^\top \bx)^2)^{1/2}}{h^2}}\nonumber\\
&\times [r\left(2-h^2 r\right)]^{d/2-1}\left(1-s^2\right)^{(d-3)/2}\, \rd r\,  \rd s. \label{eq:AhBh3}
\end{align}

Similarly, we have that
\begin{align}
\int_{\Sd}&C\lrp{\frac{1-\bx^\top\bw}{h^2}}D\lrp{\frac{1-\by^\top\bw}{h^2}}\,\sigma_d(\rd \bw)\nonumber\\
&=\,h^{d}\om{d-2}\int_{-1}^1\int_{0}^{2h^{-2}}C\lrp{\rho} D\lrp{\frac{1-(1-\rho h^2) \by^\top\bx -\left(h^2 \rho\left(2-h^2 \rho\right)\right)^{1/2} t (1-(\by^\top \bx)^2)^{1/2}}{h^2}}\nonumber\\
&\times [\rho\left(2-h^2 \rho\right)]^{d/2-1}\left(1-t^2\right)^{(d-3)/2}\, \rd \rho\,  \rd t.\label{eq:ChDh1}
\end{align}
Consider the change of variable:
$$
\by=u \bx+\left(1-u^2\right)^{1/2} \bC_{\bx} \vartb, \quad \sigma_d(\rd \by)=\left(1-u^2\right)^{d/2-1} \,\rd u \,\sigma_{d-1}(\rd \vartb)
$$
where $u \in(-1,1), \vartb \in \Sdm$ and $\bC_{\bx}=\left(\mathbf{c}_1, \ldots, \mathbf{c}_d\right)_{(d+1) \times d}$ is the semi-orthonormal matrix resulting from the completion of $\bx$ to the orthonormal basis $\left\{\bx, \mathbf{c}_1, \ldots, \mathbf{c}_d\right\}$ of $\mathbb{R}^{d+1}$. Using  \eqref{eq:AhBh3} and \eqref{eq:ChDh1}, it follows that:
\begin{align}
\eqref{eq:AhBhChDh1}=&\;h^{2d}c_{d,A}(h)c_{d,B}(h)c_{d,C}(h)c_{d,D}(h)\om{d-2}^2\nonumber\\
    &\times \int_{\Sdm}\int_{-1}^1\int_{\Sd}\Bigg\{\int_{-1}^1\int_{0}^{2h^{-2}}A\lrp{r}B\lrp{\frac{1-(1-rh^2) u -\left(h^2 r\left(2-h^2 r\right)\right)^{1/2} s (1-u^2)^{1/2}}{h^2}}\nonumber\\
    &\times [r\left(2-h^2 r\right)]^{d/2-1}\left(1-s^2\right)^{(d-3)/2}\,\rd r\,  \rd s \Bigg\}\nonumber\\
    &\times \Bigg\{\int_{-1}^1\int_{0}^{2h^{-2}}C\lrp{\rho}D\lrp{\frac{1-(1-h^2 \rho) u-\left(h^2 \rho\left(2-h^2 \rho\right)\right)^{1/2} t (1-u^2)^{1/2}}{h^2}}\nonumber\\
    &\times [\rho\left(2-h^2 \rho\right)]^{d/2-1}\left(1-t^2\right)^{(d-3)/2}\,\rd \rho\, \rd t \Bigg\} \nonumber\\
    &\times \left(1-u^2\right)^{d/2-1}\bar f(u \bx+\left(1-u^2\right)^{1/2} \bC_{\bx} \vartb)\bar f(\bx)\, \sigma_d(\rd \bx)\,\rd u\, \sigma_{d-1}(\rd \vartb). \label{eq:AhBhChDh2}
\end{align}
Consider now the change of variable
$c=\frac{1-u}{h^2}$. With this change of variable,
$u=1-ch^2$, and as a result:
\begin{align*}
1-u^2 & =h^2 c\left(2-h^2 c\right).
\end{align*}
Therefore:
\begin{align}
\eqref{eq:AhBhChDh2}=&\;h^{3d}c_{d,A}(h)c_{d,B}(h)c_{d,C}(h)c_{d,D}(h)\om{d-2}^2\nonumber\\
    &\times \int_{\Sdm}\int_{-1}^1\int_{\Sd}\Bigg\{\int_{-1}^1\int_{0}^{2h^{-2}}A\lrp{r}\nonumber\\
    &\times B\lrp{\frac{1-(1-rh^2) (1-ch^2) -\left(h^2 r\left(2-h^2 r\right)\right)^{1/2} s (h^2 c\left(2-h^2 c\right))^{1/2}}{h^2}}\nonumber\\
    &\times [r\left(2-h^2 r\right)]^{d/2-1}\left(1-s^2\right)^{(d-3)/2}\,\rd r\,  \rd s \Bigg\} \Bigg\{\int_{-1}^1\int_{0}^{2h^{-2}}C\lrp{\rho}\nonumber\\
    &\times D\lrp{\frac{1-(1-h^2 \rho) (1-ch^2)-\left(h^2 \rho\left(2-h^2 \rho\right)\right)^{1/2} t (h^2 c\left(2-h^2 c\right))^{1/2}}{h^2}}\nonumber\\
    &\times [\rho\left(2-h^2 \rho\right)]^{d/2-1}\left(1-t^2\right)^{(d-3)/2}\,\rd \rho\, \rd t \Bigg\} \nonumber\\
    &\times \left(c\left(2-h^2 c\right)\right)^{d/2-1}\bar f((1-ch^2) \bx+\left(h^2 c\left(2-h^2 c\right)\right)^{1/2} \bC_{\bx} \vartb)\bar f(\bx)\nonumber\\
    &\times \, \,\sigma_d(\rd \bx)\,\rd c \,\sigma_{d-1}(\rd \vartb).\label{eq:AhBhChDh6}
\end{align}
Using the DCT, it follows that
\begin{align}
\eqref{eq:AhBhChDh6}\sim&\,h^{3d}c_{d,A}(h)c_{d,B}(h)c_{d,C}(h)c_{d,D}(h)\om{d-2}\nonumber\\
    &\times \int_{\Sdm}\int_{0}^{\infty}\int_{\Sd}\Bigg\{\int_{-1}^{1}\int_{0}^{\infty}A\lrp{r} B\lrp{r+c-2s(rc)^{1/2}} (2r)^{d/2-1}\left( 1-s^2\right)^{(d-3)/2}\nonumber\\
    &\times \, \rd r\,  \rd s \Bigg\}\Bigg\{\int_{-1}^{1}\int_{0}^{\infty}C\lrp{\rho}D\lrp{\rho+c-2t(\rho c)^{1/2} }(2\rho)^{d/2-1}\left(1-t^2\right)^{(d-3)/2} \nonumber\\
    &\times\,\rd \rho\, \rd t\Bigg\} \left(2c\right)^{d/2-1}\bar f( \bx)^2\,\sigma_d(\rd \bx)\, \rd c \,\sigma_{d-1}(\rd \vartb)\nonumber\\
    =&\;h^{3d}c_{d,A}(h)c_{d,B}(h)c_{d,C}(h)c_{d,D}(h)\om{d-2}^2\om{d-1}R(f)2^{(3d-6)/2}\nonumber\\
    &\times \int_{0}^{\infty}\Bigg\{\int_{-1}^{1}\int_{0}^{\infty}A\lrp{r} B\lrp{r+c-2s(rc)^{1/2}}r^{d/2-1}\left( 1-s^2\right)^{(d-3)/2} \, \rd r\,  \rd s \Bigg\}\nonumber\\
    &\times \Bigg\{\int_{-1}^{1}\int_{0}^{\infty}C\lrp{\rho}D\lrp{\rho+c-2t(\rho c)^{1/2} }\rho^{d/2-1}\left(1-t^2\right)^{(d-3)/2} \,\rd \rho\, \rd t \Bigg\} c^{d/2-1}\, \rd c. \label{eq:AhBhChDh7}
\end{align}

For the calculation of \eqref{eq:AhBhChDh1} in the case $d=1$, similar steps can be followed. To begin with, consider the change of variable:
$$
\bv=s \bx+\left(1-s^2\right)^{1/2} \bB_{\bx} \xi, \quad \sigma_1(\rd \bv)=\left(1-s^2\right)^{d/2-1}\, \rd s \sigma_{0}(\rd \xi)
$$
where $s \in(-1,1), \xi \in \{-1,1\}$ and $\bB_{\bx}=\left(\mathbf{b}_1\right)_{2\times 1}$ is the vector resulting from the completion of $\bx$ to the orthonormal basis $\left\{\bx, \mathbf{b}_1\right\}$ of $\mathbb{R}^{2}$. Therefore:
\begin{align}
\int_{\Sdo}&A\lrp{\frac{1-\bx^\top\bv}{h^2}}B\lrp{\frac{1-\by^\top\bv}{h^2}}\,\sigma_d(\rd \bv)\nonumber\\
=&\;\int_{\{-1,1\}}\int_{-1}^1A\lrp{\frac{1-s}{h^2}}B\lrp{\frac{1-s \by^\top\bx -\left(1-s^2\right)^{1/2} \by^\top(\bB_{\bx} \xi)}{h^2}}\left(1-s^2\right)^{-1/2}\, \rd s\, \sigma_{0}(\rd \xi). \label{eq:AhBh11}
\end{align}

Consider now the change of variable $r=\frac{1-s}{h^2}$. With this change of variable, $s=1-rh^2$ and, as a result:
\begin{align*}
1-s^2 & =h^2 r\left(2-h^2 r\right).
\end{align*}
Therefore:
\begin{align}
\eqref{eq:AhBh11}=&\;h\int_{\{-1,1\}}\int_{0}^{2h^{-2}}A\lrp{r}B\lrp{\frac{1-(1-rh^2) \by^\top\bx -\left(h^2 r\left(2-h^2 r\right)\right)^{1/2} \by^\top(\bB_{\bx} \xi)}{h^2}}\nonumber\\
    &\times [r\left(2-h^2 r\right)]^{-1/2}\,\rd r\, \sigma_{0}(\rd \xi)\nonumber\\
    =&h\int_{0}^{2h^{-2}}A\lrp{r}\Bigg\{B\lrp{\frac{1-(1-rh^2) \by^\top\bx -\left(h^2 r\left(2-h^2 r\right)\right)^{1/2} \by^\top\bB_{\bx}}{h^2}}\nonumber\\
    &+B\lrp{\frac{1-(1-rh^2) \by^\top\bx +\left(h^2 r\left(2-h^2 r\right)\right)^{1/2} \by^\top\bB_{\bx}}{h^2}}\Bigg\}\nonumber\\
    &\times [r\left(2-h^2 r\right)]^{-1/2}\,\rd r.\label{eq:AhBh21}
\end{align}

Similarly,
\begin{align}
\int_{\Sdo}&C\lrp{\frac{1-\bx^\top\bw}{h^2}}D\lrp{\frac{1-\by^\top\bw}{h^2}}\,\sigma_d(\rd \bw)\nonumber\\
=&h\int_{0}^{2h^{-2}}C\lrp{\rho}\Bigg\{D\lrp{\frac{1-(1-\rho h^2) \by^\top\bx -\left(h^2 \rho\left(2-h^2 \rho\right)\right)^{1/2} \by^\top\tilde{\bB}_{\bx}}{h^2}}\nonumber\\
    &+D\lrp{\frac{1-(1-\rho h^2) \by^\top\bx +\left(h^2 \rho\left(2-h^2 \rho\right)\right)^{1/2} \by^\top\tilde{\bB}_{\bx}}{h^2}}\Bigg\}\nonumber\\
    &\times [\rho\left(2-h^2 \rho\right)]^{-1/2}\,\rd \rho,\label{eq:ChDh21}
\end{align}
where $\tilde{\bB}_{\bx}=\left(\tilde{\mathbf{b}}_1\right)_{2\times 1}$ is the vector resulting from the completion of $\bx$ to the orthonormal basis $\left\{\bx, \tilde{\mathbf{b}}_1\right\}$ of $\mathbb{R}^{2}$. Consider the change of variable:
$$
\by=u \bx+\left(1-u^2\right)^{1/2} \bC_{\bx} \vartheta, \quad \sigma_1(\rd \by)=\left(1-u^2\right)^{-1/2}\,\rd u \,\sigma_{0}(\rd \vartheta)
$$
where $u \in(-1,1), \vartheta \in \{-1,1\}$ and $\bC_{\bx}=\left(\mathbf{c}_1\right)_{2 \times 1}$ is the vector resulting from the completion of $\bx$ to the orthonormal basis $\left\{\bx, \mathbf{c}_1\right\}$ of $\mathbb{R}^{2}$. Using  \eqref{eq:AhBh21} and \eqref{eq:ChDh21}, it follows that:
\begin{align}
\eqref{eq:AhBhChDh1}=&\;h^{2}c_{1,A}(h)c_{1,B}(h)c_{1,C}(h)c_{1,D}(h)\nonumber\\
    &\times \int_{\{-1,1\}}\int_{-1}^1\int_{\Sdo}\Bigg[\int_{0}^{2h^{-2}}A\lrp{r}\nonumber\\
    &\times \Bigg\{B\lrp{\frac{1-(1-rh^2) u -\left(h^2 r\left(2-h^2 r\right)\right)^{1/2} \left(1-u^2\right)^{1/2} (\bC_{\bx} \vartheta)^\top\bB_{\bx}}{h^2}}\nonumber\\
    &+B\lrp{\frac{1-(1-rh^2) u +\left(h^2 r\left(2-h^2 r\right)\right)^{1/2} \left(1-u^2\right)^{1/2} (\bC_{\bx} \vartheta)^\top\bB_{\bx}}{h^2}}\Bigg\}\nonumber\\
    &\times [r\left(2-h^2 r\right)]^{-1/2}\,\rd r\Bigg]\nonumber\\
    &\times \Bigg[\int_{0}^{2h^{-2}}C\lrp{\rho}\Bigg\{D\lrp{\frac{1-(1-\rho h^2) u -\left(h^2 \rho\left(2-h^2 \rho\right)\right)^{1/2} \left(1-u^2\right)^{1/2} (\bC_{\bx} \vartheta)^\top\tilde{\bB}_{\bx}}{h^2}}\nonumber\\
    &+D\lrp{\frac{1-(1-\rho h^2) u +\left(h^2 \rho\left(2-h^2 \rho\right)\right)^{1/2} \left(1-u^2\right)^{1/2} (\bC_{\bx} \vartheta)^\top\tilde{\bB}_{\bx}}{h^2}}\Bigg\}\nonumber\\
    &\times [\rho\left(2-h^2 \rho\right)]^{-1/2}\,\rd \rho\Bigg]\bar f(\bx)\bar f(u \bx+\left(1-u^2\right)^{1/2} \bC_{\bx} \vartheta)\left(1-u^2\right)^{-1/2}\nonumber\\
    &\times \,\sigma_1(\rd \bx)\,\rd u\,\sigma_0(\rd \vartheta)\nonumber\\
    =&\;h^{2}c_{1,A}(h)c_{1,B}(h)c_{1,C}(h)c_{1,D}(h)\nonumber\\
    &\times \int_{-1}^1\int_{\Sdo}\Bigg[\int_{0}^{2h^{-2}}A\lrp{r}\Bigg\{B\lrp{\frac{1-(1-rh^2) u -\left(h^2 r\left(2-h^2 r\right)\right)^{1/2} \left(1-u^2\right)^{1/2} }{h^2}}\nonumber\\
    &+B\lrp{\frac{1-(1-rh^2) u +\left(h^2 r\left(2-h^2 r\right)\right)^{1/2} \left(1-u^2\right)^{1/2} }{h^2}}\Bigg\}\nonumber\\
    &\times [r\left(2-h^2 r\right)]^{-1/2}\,\rd r\Bigg]\nonumber\\
    &\times \Bigg[\int_{0}^{2h^{-2}}C\lrp{\rho}\Bigg\{D\lrp{\frac{1-(1-\rho h^2) u -\left(h^2 \rho\left(2-h^2 \rho\right)\right)^{1/2} \left(1-u^2\right)^{1/2}}{h^2}}\nonumber\\
    &+D\lrp{\frac{1-(1-\rho h^2) u +\left(h^2 \rho\left(2-h^2 \rho\right)\right)^{1/2} \left(1-u^2\right)^{1/2}}{h^2}}\Bigg\}\nonumber\\
    &\times [\rho\left(2-h^2 \rho\right)]^{-1/2}\,\rd \rho\Bigg]\bar f(\bx) \{\bar f(u \bx-\left(1-u^2\right)^{1/2} \bC_{\bx} )+\bar f(u \bx+\left(1-u^2\right)^{1/2} \bC_{\bx} )\}\nonumber\\
    &\times \left(1-u^2\right)^{-1/2} \,\sigma_1(\rd \bx)\,\rd u.\label{eq:AhBhChDh21}
\end{align}
Consider now the change of variable
$c=\frac{1-u}{h^2}$. With this change of variable,
$u=1-ch^2$, and as a result:
\begin{align*}
1-u^2 & =h^2 c\left(2-h^2 c\right).
\end{align*}
Therefore:
\begin{align}
\eqref{eq:AhBhChDh21}=&\;h^{3}c_{1,A}(h)c_{1,B}(h)c_{1,C}(h)c_{1,D}(h)\nonumber\\
    &\times \int_{0}^{2h^{-2}}\int_{\Sdo}\Bigg[\int_{0}^{2h^{-2}}A\lrp{r}\Bigg\{B\lrp{\frac{1-(1-rh^2) (1-ch^2) -\left(h^2 r\left(2-h^2 r\right)\right)^{1/2} \left(h^2 c\left(2-h^2 c\right)\right)^{1/2} }{h^2}}\nonumber\\
    &+B\lrp{\frac{1-(1-rh^2) (1-ch^2) +\left(h^2 r\left(2-h^2 r\right)\right)^{1/2} \left(h^2 c\left(2-h^2 c\right)\right)^{1/2} }{h^2}}\Bigg\}\nonumber\\
    &\times [r\left(2-h^2 r\right)]^{-1/2}\,\rd r\Bigg]\nonumber\\
    &\times \Bigg[\int_{0}^{2h^{-2}}C\lrp{\rho}\Bigg\{D\lrp{\frac{1-(1-\rho h^2) (1-ch^2) -\left(h^2 \rho\left(2-h^2 \rho\right)\right)^{1/2} \left(h^2 c\left(2-h^2 c\right)\right)^{1/2}}{h^2}}\nonumber\\
    &+D\lrp{\frac{1-(1-\rho h^2) (1-ch^2) +\left(h^2 \rho\left(2-h^2 \rho\right)\right)^{1/2} \left(h^2 c\left(2-h^2 c\right)\right)^{1/2}}{h^2}}\Bigg\}\nonumber\\
    &\times [\rho\left(2-h^2 \rho\right)]^{-1/2}\,\rd \rho\Bigg]\nonumber\\
    &\times \bar f(\bx) \{\bar f((1-ch^2) \bx-\left(h^2 c\left(2-h^2 c\right)\right)^{1/2} \bC_{\bx} )+\bar f((1-ch^2) \bx+\left(h^2 c\left(2-h^2 c\right)\right)^{1/2} \bC_{\bx} )\}\nonumber\\
    &\times \left(c\left(2-h^2 c\right)\right)^{-1/2} \,\sigma_1(\rd \bx)\,\rd c.\label{eq:AhBhChDh22}
\end{align}

Using the DCT, it follows that
\begin{align}
\eqref{eq:AhBhChDh22}\sim &\,h^{3}c_{1,A}(h)c_{1,B}(h)c_{1,C}(h)c_{1,D}(h) R(f) 2^{-1/2}\int_{0}^\infty\Bigg[\int_{0}^{\infty}A\lrp{r}\Big\{B\lrp{r+c-2(rc)^{1/2}}\nonumber\\
    &+B\lrp{r+c+2(rc)^{1/2}}\Big\} r^{-1/2}\,\rd r\Bigg]\nonumber\\
    &\times \Bigg[\int_{0}^{\infty}C\lrp{\rho}\Big\{D\lrp{\rho+c-2(\rho c)^{1/2}}+D\lrp{\rho+c+2(\rho c)^{1/2}}\Big\} \rho^{-1/2}\,\rd \rho\Bigg]\nonumber\\
    &\times c^{-1/2} \,\rd c.\label{eq:AhBhChDh23}
\end{align}

Therefore, using \eqref{eq:AhBhChDh7} and \eqref{eq:AhBhChDh23} it follows that
\begin{align*}
\eqref{eq:AhBhChDh1}\sim &\,h^{3d}c_{d,A}(h)c_{d,B}(h)c_{d,C}(h) c_{d,D}(h) R(f)  \tilde{\gamma}_d \\
&\times \int_0^{\infty}\Bigg\{ \int_0^{\infty} A(r) r^{-1/2}\varphi_d(B, r, c)\,\rd r\int_0^{\infty} C(\rho) \rho^{-1/2}\varphi_d(D, \rho, c)\,\rd \rho\Bigg\}c^{-1/2}\, \rd c, \\
\tilde{\gamma}_d= & \begin{cases}2^{-1/2}, & d=1, \\
\om{d-2}^2\om{d-1}2^{(3d-6)/2}, & d \geq 2 .\end{cases}
\end{align*}
\end{proof}

\section{Proofs for the Euclidean case}
\label{sec:euclidean}

Here we provide a brief sketch for the proof of \eqref{eq:ANeu}. Denote $V=W*W-2W$. The derivative of $V_h$ with respect to $h$ becomes
$$
\frac{\partial}{\partial h}V_h(\bx)=\frac{\partial}{\partial h}\big\{h^{-d}V(\bx/h)\big\}=-dh^{-d-1}V(\bx/h)-h^{-d}\D V(x/h)^\top\bx/h^2=h^{-1}\lambda_h(\bx),
$$
where $\lambda(\bx)=-dV(\bx)-\bx^\top\D V(\bx)$. Hence, the derivative of the cross-validation criterion \eqref{eq:cveu} is
$$
{\rm CV}'(h)=-dn^{-1}h^{-d-1}R(W)+\textstyle{\binom{n}{2}}^{-1}h^{-1}\displaystyle\sum_{1\leq i<j\leq n}\lambda_h(\bX_i-\bX_j).
$$

As in the directional case, in order to derive the limit distribution of the cross-validation bandwidth we need to analyze the asymptotic behavior of ${\rm Var}\{{\rm CV}'(h_A)\}$, where $h_A=c_An^{-1/(d+4)}$, with $c_A=[dR(W)/\{\mu_2(W)^2R(\nabla^2f)\}]^{1/(d+4)}$. Again, using Equation \eqref{eq:varU}, such a variance can be exactly written as
\begin{align*}
{\rm Var}\{{\rm CV}'(h_A)\}&=\frac{2}{n(n-1)}h_A^{-2}{\rm Var}\{\lambda_{h_A}(\bX_1-\bX_2)\}\\
&\quad+\frac{4(n-2)}{n(n-1)}h_A^{-2}{\rm Cov}\{\lambda_{h_A}(\bX_1-\bX_2),\lambda_{h_A}(\bX_1-\bX_3)\}.
\end{align*}
The previous formula involves three expectations, which we analyze next.

It can be shown that the vector moments $\bmu_j(\lambda)\defin\int_{\mathbb R^d}\bx^{\otimes j}\lambda(\bx)\rd\bx$ satisfy $\bmu_j(\lambda)=0$ for $j=0,1,2,3$ and $\bmu_4(\lambda)=24\boldsymbol{\mathcal S}_{d,4}\bmu_2(W)^{\otimes 2}$. This allows writing $(\lambda_{h}*f)(\bx)\sim h^4\bmu_2(W)^{\otimes 2\top}\D^{\otimes 4}f(\bx)$ as $h\to0$. Therefore, as $h\to0$ we can express
\begin{align}
    \mathbb E\{\lambda_{h}(\bX_1-\bX_2)\}&=\int_{\mathbb R^d}(\lambda_{h}*f)(\bx)f(\bx)d\bx\sim h^4\bmu_2(W)^{\otimes 2\top}\bpsi_4,\label{eq:cveu1}\\
    \mathbb E\{\lambda_{h}(\bX_1-\bX_2)\lambda_{h}(\bX_1-\bX_3)\}&=\int_{\mathbb R^d}(\lambda_{h}*f)(\bx)^2f(\bx)d\bx\sim h^8\bmu_2(W)^{\otimes 4\top}\bpsi_{4,4}\label{eq:cveu2},
\end{align}
where $\bpsi_4=\int_{\mathbb R^d}\D^{\otimes 4}f(\bx)f(\bx)d\bx$ and $\bpsi_{4,4}=\int_{\mathbb R^d}\{\D^{\otimes 4}f(\bx)\}^{\otimes 2}f(\bx)d\bx$. Finally, we need to find the order of $\mathbb E\{\lambda_{h}(\bX_1-\bX_2)^2\}=\int_{\mathbb R^d}\lambda_{h}(\bx-\by)^2f(\bx)f(\by)d\bx d\by$.
Noting that $\lambda_{h}(\bx)^2=h^{-2d}\lambda(\bx/h)^2=h^{-d}(\lambda^2)_{h}(\bx)$ we have $\mathbb E\{\lambda_{h}(\bX_1-\bX_2)^2\}=h^{-d}\int_{\mathbb R^d}\{(\lambda^2)_{h}*f\}(\bx)f(\bx)d\bx$ so
\begin{align}\label{eq:cveu3}
  \mathbb E\{\lambda_{h}(\bX_1-\bX_2)^2\}\sim h^{-d}R(\lambda)R(f)
\end{align}
as $h\to0$.
Hence, putting together \eqref{eq:cveu1}, \eqref{eq:cveu2} and \eqref{eq:cveu3} we obtain
\begin{align*}
    {\rm Var}\{{\rm CV}'(h_A)\}&\sim 2n^{-2}h_A^{-2}\big\{ h_A^{-d}R(\lambda)R(f)  -h_A^8\big[\bmu_2(W)^{\otimes 2\top}\bpsi_4]^2\big\}\\
    &\quad+4n^{-1}h_A^{-2}\big\{h_A^8\bmu_2(W)^{\otimes 4\top}\bpsi_{4,4}-h_A^8\big[\bmu_2(W)^{\otimes 2\top}\bpsi_4]^2\big\}.
\end{align*}
Comparing the orders of each term, eventually, we obtain
\begin{align*}
   {\rm Var}\{{\rm CV}'(h_A)\}&\sim2n^{-2}h_A^{-d-2}R(\lambda)R(f)=2c_A^{-d-2}R(\lambda)R(f)n^{-(d+6)/(d+4)}.
\end{align*}
Then, reasoning analogously as in the directional case, it follows that the relative error of the cross-validation bandwidth satisfies $n^{d/(2d+8)}(\tilde h_{\rm CV}-h_{\rm MISE})/h_{\rm MISE}\stackrel{d}{\longrightarrow}\mathcal N(0,\sigma_{\rm CV}^2)$, with
\begin{align*}
    \sigma_{\rm CV}^2=2R(\lambda)R(f)(d+4)^{-2}[\mu_2(W)^{2d}R(\nabla^2 f)^d\{dR(W)\}^{d+8}]^{-1/(d+4)}.
\end{align*}

The only thing left is to prove that $R(\lambda)=R(\rho)$ and to show the exact value of the latter for the Gaussian kernel. For the former, note that
$$
2\int_{\mathbb R^d}V(\bx)\bx^\top\D V(\bx)=\int_{\mathbb R^d}\bx^\top\D\{V(\bx)^2\}\,\rd\bx=-d\int_{\mathbb R^d}V(\bx)^2\,\rd\bx
$$
(by a univariate application of integration by parts to the partial derivatives), so if we define $\rho(\bx)\defin\bx^\top\D V(\bx)=\bx^\top\D(W*W-2W)(\bx)$ then we have
\begin{align*}
    R(\lambda)&=\int_{\mathbb R^d}\{dV(\bx)+\bx^\top\D V(\bx)\}^2\,\rd\bx\\
    &=d^2\int_{\mathbb R^d}V(\bx)^2\,\rd\bx+\int_{\mathbb R^d}\{\bx^\top\D V(\bx)\}^2\,\rd\bx+2d\int_{\mathbb R^d}V(\bx)\bx^\top\D V(\bx)\,\rd\bx=R(\rho).
\end{align*}

In the following, we will compute $R(\rho)$ for the case where a Gaussian kernel is used in the kernel density estimator. Let us denote by $\phi_{\bSigma}$ the density of the $d$-variate $\mathcal N({\bf 0},\bSigma)$ distribution. When the standard Gaussian kernel $W_{\rm G}=\phi_{\bI_d}$ is used, then $V_{\rm G}=W_{\rm G}*W_{\rm G}-2W_{\rm G}=\phi_{2\bI_d}-2\phi_{\bI_d}$. Since $\D\phi_{\bSigma}(\bx)=-\phi_{\bSigma}(\bx)\bSigma^{-1}\bx$, it follows that $\D V_{\rm G}(\bx)=2\bx\phi_{\bI_d}(\bx)-\frac12\bx\phi_{2\bI_d}(\bx)$ and
$\rho_{\rm G}(\bx)=2\bx^\top\bx\phi_{\bI_d}(\bx)-\frac12\bx^\top\bx\phi_{2\bI_d}(\bx)$. Therefore,
\begin{align}
    R(\rho_{\rm G})
    &=4\int_{\mathbb R^d}(\bx^\top\bx)^2\phi_{\bI_d}(\bx)^2\,\rd\bx+\frac14\int_{\mathbb R^d}(\bx^\top\bx)^2\phi_{2\bI_d}(\bx)^2\,\rd\bx\nonumber\\
    &\quad-2\int_{\mathbb R^d}(\bx^\top\bx)^2\phi_{\bI_d}(\bx)\phi_{2\bI_d}(\bx)\,\rd\bx\label{eq:RrhoG}
\end{align}
To compute each of the terms in \eqref{eq:RrhoG} take into account that, according to Fact C.2.1 in \citet{Wand1994},
\begin{align}\label{eq:e1}
    \phi_{c_1\bI_d}(\bx)\phi_{c_2\bI_d}(\bx)=\phi_{(c_1+c_2)\bI_d}({\bf 0})\phi_{\frac{c_1c_2}{c_1+c_2}\bI_d}(\bx)=\{2\pi(c_1+c_2)\}^{-d/2}\phi_{\frac{c_1c_2}{c_1+c_2}\bI_d}(\bx).
\end{align}
And also that
\begin{align}
 \int_{\mathbb R^d}(\bx^\top\bx)^2\phi_{c\bI_d}(\bx)d\bx&=   c^2\int_{\mathbb R^d}(\bz^\top\bz)^2\phi_{\bI_d}(\bz)d\bz=c^2d(d+2),  \label{eq:e2}
\end{align}
since the raw second-order moment of a $\chi_d^2$ distribution is $2d+d^2=d(d+2)$. So, from \eqref{eq:RrhoG}, using \eqref{eq:e1} and \eqref{eq:e2} we obtain
\begin{align*}
    (2\pi)^{d/2}R(\rho_{\rm G})&=4(2)^{-d/2}\int_{\mathbb R^d}(\bx^\top\bx)^2\phi_{\frac12\bI_d}(\bx)\,\rd\bx+\frac14(4)^{-d/2}\int_{\mathbb R^d}(\bx^\top\bx)^2\phi_{\bI_d}(\bx)\,\rd\bx\nonumber\\
    &\quad-2(3)^{-d/2}\int_{\mathbb R^d}(\bx^\top\bx)^2\phi_{\frac{2}{3}\bI_d}(\bx)\,\rd\bx\\
    &=d(d+2)\Big\{(2)^{-d/2}+\frac14(4)^{-d/2}-2(3)^{-d/2}\frac{4}{9}\Big\}.
\end{align*}
Simplifying,
$R(\rho_{\rm G})=(2\pi)^{-d/2}d(d+2)(2^{-d/2}+\tfrac142^{-d}-\tfrac893^{-d/2}),$
as announced.

\section{Calculations for the von Mises--Fisher kernel}
\label{sec:calc}

\begin{lemma}\label{lemma:lambdaLGvmf}
    Let $L(s)=e^{-s}$ and $G(s)=-se^{-s}$. We have that
    $$
    \lambda_d(L)=(2\pi)^{d/2},\quad \lambda_d(L^2)=\pi^{d/2},\quad \lambda_d(G^2)=\frac{d(d+2)}{16}\pi^{d/2}.
    $$
\end{lemma}

\begin{proof}[Proof of Lemma \ref{lemma:lambdaLGvmf}]
The proof readily follows from \eqref{eq:lambdad} and the gamma function definition.%
\end{proof}

\begin{lemma} \label{lem:cGvmf}
Let $L(s)=e^{-s}$ and $G(s)=-se^{-s}$. Then:
\begin{align*}
    c_{d,L}(h)^{-1}=&\;e^{-1/h^2}c_{d}^\mathrm{vMF}(1/h^2)^{-1},\\
    c_{d,G}(h)^{-1}=&\;\frac{1}{h^2}\lrb{\frac{c_{d+2,L}(h)^{-1}}{2\pi h^2}-c_{d,L}(h)^{-1}}=\frac{e^{-1/h^2}}{h^2}\lrb{\frac{c_{d+2}^{\mathrm{vMF}}(1/h^2)^{-1}}{2\pi h^2}-c_{d}^{\mathrm{vMF}}(1/h^2)^{-1}}.
\end{align*}
\end{lemma}

\begin{proof}[Proof of Lemma \ref{lem:cGvmf}]
That $c_{d,L}(h)^{-1}=e^{-1/h^2}c_{d}^\mathrm{vMF}(1/h^2)^{-1}$ immediately follows from \eqref{eq:vmf} and $c_{d}^\mathrm{vMF}(\kappa)^{-1}=\int_{\Sd}e^{\kappa\bx^\top\bmu}\,\sigma_d(\rd \bx)$.

For the second equality, we will use that
\begin{align}
    \int_{-1}^{1}t^m (1-t^2)^{\nu-1/2} e^{\kappa t}\, \rd t=\frac{2^\nu\pi^{1 / 2} \Gamma(\nu+1/2)}{\kappa^{\nu}}\mathcal{I}_{\nu+m}(\kappa),\quad m=0,1,\label{eq:Inum}
\end{align}
with $\nu=(d-1)/2\geq 0$. Using \eqref{eq:Inum} with $m=1$ and \eqref{eq:vmf},
\begin{align*}
    c_{d,G}(h)^{-1}
    =&-\int_{\Sd} \frac{1-\bx^\top \by}{h^2} e^{-(1-\bx^\top \by)/h^2} \,\sigma_d(\rd \by)\\
    =&\;\frac{1}{h^2}\lrb{-\int_{\Sd} e^{-(1-\bx^\top \by)/h^2} \,\sigma_d(\rd \by)+\int_{\Sd} \bx^\top \by \, e^{-(1-\bx^\top \by)/h^2} \,\sigma_d(\rd \by)}\\
    =&\;\frac{1}{h^2}\lrb{-c_{d,L}(h)^{-1}+e^{-1/h^2}\om{d-1}\int_{-1}^1 t e^{t/h^2} (1-t^2)^{d/2-1}\,\rd t}\\
    =&\;\frac{1}{h^2}\lrb{-c_{d,L}(h)^{-1}+e^{-1/h^2}\om{d-1} \frac{2^{(d-1)/2}\pi^{1/2}\Gamma(d/2)}{(1/h^2)^{(d-1)/2}}\mathcal{I}_{(d+1)/2}(1/h^2)}\\
    =&-\frac{c_{d,L}(h)^{-1}}{h^2}+e^{-1/h^2}\om{d-1} \frac{2^{(d-1)/2}\pi^{1/2}\Gamma(d/2)}{h^2(1/h^2)^{(d-1)/2}}\mathcal{I}_{(d+1)/2}(1/h^2)\\
    =&-\frac{c_{d,L}(h)^{-1}}{h^2}+e^{-1/h^2}\frac{1}{2\pi h^4}\frac{(2\pi)^{(d+3)/2}}{(1/h^2)^{(d+1)/2}}\mathcal{I}_{(d+1)/2}(1/h^2)\\
    =&-\frac{c_{d,L}(h)^{-1}}{h^2}+e^{-1/h^2}\frac{c_{d+2}^{\mathrm{vMF}}(h)^{-1}}{2\pi h^4},
\end{align*}
since $c_{d}^\mathrm{vMF}(1/h^2)^{-1}=(2\pi)^{(d+1)/2}\mathcal{I}_{(d-1)/2}(1/h^2)/\lrb{(1/h^2)^{(d-1)/2}}$.
\end{proof}

\begin{lemma} \label{lem:vmfstar}
Let $L(s)=e^{-s}$ and $G(s)=-se^{-s}$, and $f_{\mathrm{vMF}}(\cdot;\bmu,\kappa)$ defined in \eqref{eq:vmf}. We have that
\begin{align*}
    (L_h*f_{\mathrm{vMF}}(\cdot;\bmu,\kappa))(\bx)=&\;
    \frac{c^\mathrm{vMF}_{d}(\kappa)c^\mathrm{vMF}_{d}(1/h^2)}{c^\mathrm{vMF}_{d}(\|\kappa\bmu+\bx/h^2\|)},\\
       (L_h*L_h)(\bx,\by)=&\;
    \frac{c^\mathrm{vMF}_{d}(1/h^2)^2}{c^\mathrm{vMF}_{d}(\|\bx+\by\|/h^2)},\\
    (G_h*L_h)(\bx,\by)=&\;c_{d}^{\mathrm{vMF}}(1/h^2)\\
    &\times \frac{(1+\bx^\top\by)\lrb{2\pi h^2 c_{d+2}^{\mathrm{vMF}}(\|\bx+\by\|/h^2)}^{-1}-c_d^{\mathrm{vMF}}(\|\bx+\by\|/h^2)^{-1}}{\lrb{2\pi h^2 c_{d+2}^{\mathrm{vMF}}(1/h^2)}^{-1}-c_{d}^{\mathrm{vMF}}(1/h^2)^{-1}}.
\end{align*}
\end{lemma}

\begin{proof}[Proof of Lemma \ref{lem:vmfstar}]
The form of the density \eqref{eq:vmf} and its normalizing constant entails that
\begin{align*}
    (L_h*f_{\mathrm{vMF}}(\cdot;\bmu,\kappa))(\bx)=&\;\int_{\Sd} L_h(\bx,\by) f_{\mathrm{vMF}}(\by;\bmu,\kappa) \,\sigma_d(\rd \by)\\
    =&\;c_{d,L}(h)e^{-1/h^2}c^\mathrm{vMF}_{d}(\kappa) \int_{\Sd} e^{\kappa\by^\top\bmu+\bx^\top \by/h^2} \,\sigma_d(\rd \by)\\
    =&\;c^\mathrm{vMF}_{d}(1/h^2)c^\mathrm{vMF}_{d}(\kappa) \int_{\Sd} e^{\|\kappa\bmu+\bx/h^2\|\by^\top(\kappa\bmu+\bx/h^2)/\|\kappa\bmu+\bx/h^2\|} \,\sigma_d(\rd \by)\\
    =&\;\frac{c^\mathrm{vMF}_{d}(\kappa)c^\mathrm{vMF}_{d}(1/h^2)}{c^\mathrm{vMF}_{d}(\|\kappa\bmu+\bx/h^2\|)}.
\end{align*}

The second statement readily follows from $(L_h*L_h)(\bx,\by)=\{L_h * L_h(\by, \cdot)\}(\bx)=\int_{\Sd} L_h(\bx,\bz) L_h(\by,\bz)\,\sigma_d(\rd \bz)$ and $L_h(\by,\bz)=f_{\mathrm{vMF}}(\bz;\by,1/h^2)$.

The third equality arises from
\begin{align*}
    (G_h*L_h)(\bx,\by)=&\;\int_{\Sd} G_h(\bx,\bz) L_h(\bz,\by) \,\sigma_d(\rd \bz)\\
    =&-c_{d,L}(h)c_{d,G}(h)e^{-2/h^2}\int_{\Sd} \frac{1-\bz^\top \by}{h^2} \exp\{\bz^\top (\bx+\by)/h^2\} \,\sigma_d(\rd \bz)\\
    =&-c_{d,L}(h)c_{d,G}(h)e^{-2/h^2}\\
    &\times \int_{-1}^1 \int_{\Sdm} \frac{1-\by^\top\{t(\bx+\by)/\|\bx+\by\|+(1-t^2)^{d/2-1}\bB_{(\bx+\by)/\|\bx+\by\|}\bxi\}}{h^2}\\
    &\times\exp\{[\|\bx+\by\|/h^2]t\} (1-t^2)^{d/2-1}\,\sigma_{d-1}(\rd \bxi)\,\rd t\\
    =&-c_{d,L}(h)c_{d,G}(h)e^{-2/h^2}\\
    &\times\om{d-1}\int_{-1}^1 \frac{1-t(1+\bx^\top\by)/\|\bx+\by\|}{h^2}
    \exp\{[\|\bx+\by\|/h^2]t\} (1-t^2)^{d/2-1}\,\rd t.
\end{align*}
We use \eqref{eq:Inum} with $\nu=(d-1)/2\geq0$ to have
\begin{align*}
    \om{d-1}\int_{-1}^1 &\frac{1-t(1+\bx^\top\by)/\|\bx+\by\|}{h^2} \exp\{[\|\bx+\by\|/h^2]t\} (1-t^2)^{d/2-1}\,\rd t\\
    &=\om{d-1}\frac{2^{(d-1)/2}\pi^{1 / 2} \Gamma(d/2)}{(\|\bx+\by\|/h^2)^{(d-1)/2}}\frac{1}{h^2}\lrb{\mathcal{I}_{(d-1)/2}(\|\bx+\by\|/h^2)-\frac{(1+\bx^\top\by)}{\|\bx+\by\|}\mathcal{I}_{(d+1)/2}(\|\bx+\by\|/h^2)}\\
    &=\frac{(2\pi)^{(d+1)/2}}{(\|\bx+\by\|/h^2)^{(d-1)/2}}\frac{1}{h^2}\lrb{\mathcal{I}_{(d-1)/2}(\|\bx+\by\|/h^2)-\frac{(1+\bx^\top\by)}{\|\bx+\by\|}\mathcal{I}_{(d+1)/2}(\|\bx+\by\|/h^2)}\\
    &=\frac{1}{h^2}\lrb{c_d^{\mathrm{vMF}}(\|\bx+\by\|/h^2)^{-1}-\frac{(1+\bx^\top\by)c_{d+2}^{\mathrm{vMF}}(\|\bx+\by\|/h^2)^{-1}}{2\pi h^2}},
\end{align*}
since $c_{d}^\mathrm{vMF}(1/h^2)^{-1}=(2\pi)^{(d+1)/2}\mathcal{I}_{(d-1)/2}(1/h^2)/\lrb{(1/h^2)^{(d-1)/2}}$. Hence, using Lemma \ref{lem:cGvmf},
\begin{align*}
    (G_h*L_h)(\bx,\by)=&-c_{d,L}(h)c_{d,G}(h)e^{-2/h^2}\\
    &\times \frac{1}{h^2}\lrb{c_d^{\mathrm{vMF}}(\|\bx+\by\|/h^2)^{-1}-\frac{(1+\bx^\top\by)c_{d+2}^{\mathrm{vMF}}(\|\bx+\by\|/h^2)^{-1}}{2\pi h^2}}\\
    =&-c_{d}^{\mathrm{vMF}}(1/h^2)\lrb{\frac{c_{d+2}^{\mathrm{vMF}}(1/h^2)^{-1}}{2\pi h^2}-c_{d}^{\mathrm{vMF}}(1/h^2)^{-1}}^{-1}\\
    &\times \lrb{c_d^{\mathrm{vMF}}(\|\bx+\by\|/h^2)^{-1}-\frac{(1+\bx^\top\by)c_{d+2}^{\mathrm{vMF}}(\|\bx+\by\|/h^2)^{-1}}{2\pi h^2}}\\
    =&\;c_{d}^{\mathrm{vMF}}(1/h^2)\\
    &\times \frac{(1+\bx^\top\by)\lrb{2\pi h^2 c_{d+2}^{\mathrm{vMF}}(\|\bx+\by\|/h^2)}^{-1}-c_d^{\mathrm{vMF}}(\|\bx+\by\|/h^2)^{-1}}{\lrb{2\pi h^2 c_{d+2}^{\mathrm{vMF}}(1/h^2)}^{-1}-c_{d}^{\mathrm{vMF}}(1/h^2)^{-1}}.
\end{align*}
\end{proof}

Lemmas \ref{lem:cGvmf} and \ref{lem:vmfstar} give a means of estimating $\mu_1(h)$, the leading term of $\mathbb{V}\mathrm{ar}\{{\rm CV}'(h)\}$ in \eqref{eq:CV1}, for $L_{\mathrm{vMF}}$ and general density $f$. By Monte Carlo, given $(\bX_i,\bY_i)\sim f\times f$ iid, $i=1,\ldots,M$, from \eqref{eq:mu11}:
\begin{align*}
    \mu_1(h)=&\;\frac{16c_{d,L}(h)^2}{h^2 c_{d,G}(h)^2}\int_{\Sd}\int_{\Sd}\{(L_h-G_h)*L_h-(L_h-G_h)\}(\bx,\by)^2f(\bx)f(\by)\,\sigmad(\rd \bx)\,\sigmad(\rd \by)\\
    \approx &\;\frac{16c_{d,L}(h)^2}{h^2 c_{d,G}(h)^2}\frac{1}{M}\sum_{i=1}^n \{L_h*L_h-G_h*L_h-L_h+G_h\}(\bX_i,\bY_i)^2,
\end{align*}
where the normalizing constants follow from Lemma \ref{lem:cGvmf} and the convolutions from Lemma \ref{lem:vmfstar}. This allows computing $\mu_1(h)$ without asymptotic expansions, and has been used to validate numerically the asymptotic expressions in Propositions \ref{prop:AhBhCh} and \ref{prop:AhBhChDh}.

\section{Additional numerical experiments}
\label{sec:add}

This section contains additional plots from the numerical experiments in Sections \ref{sec:conv} and~\ref{sec:num}.

Figure \ref{fig:rho_kappa} gives the complementary view to Figure \ref{fig:sigma2:rho}, showing $\kappa\mapsto\rho_d(\kappa)$ for fixed $d$. It shows these curves stabilize as $\kappa$ increases and that only $d=1$ has a global minimum.

\begin{figure}[htb!]
  \centering
  \includegraphics[width=0.5\linewidth]{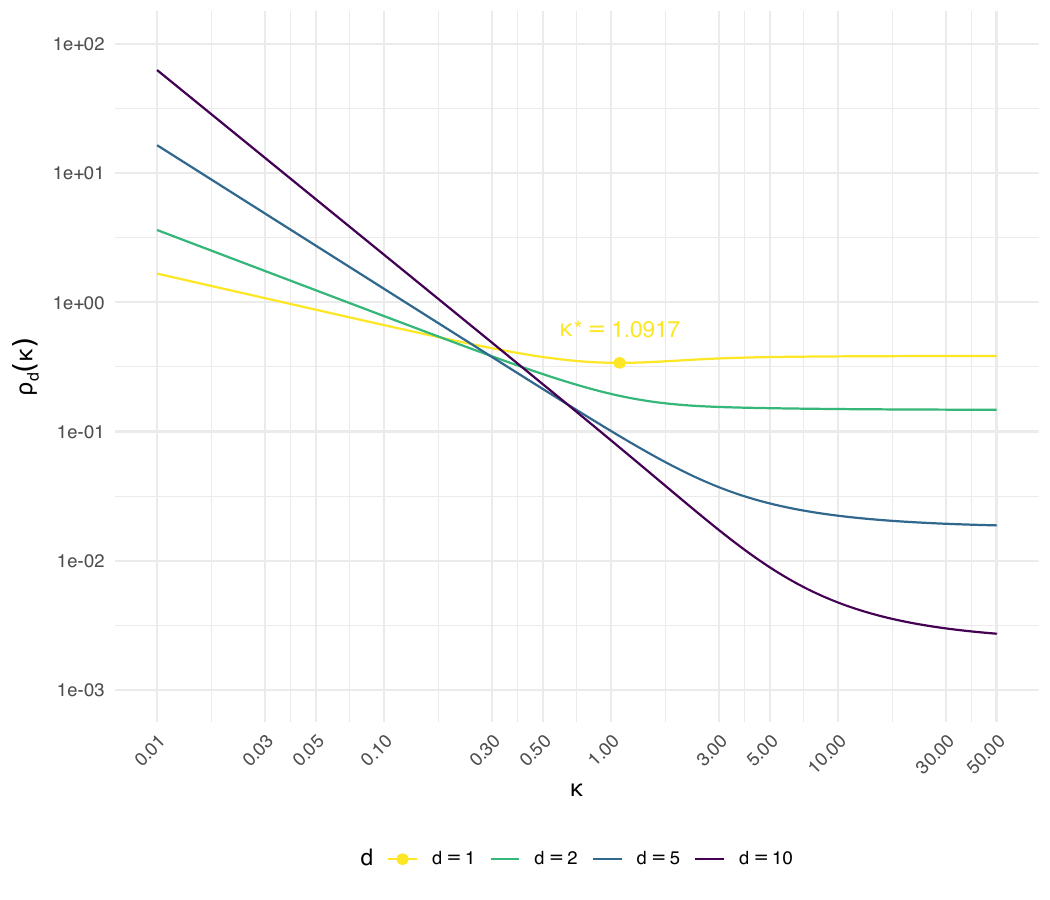}
  \caption{Curves $\kappa\mapsto\rho_d(\kappa)$ giving the contribution of the vMF density to the asymptotic variance for several dimensions. The global minimum for $d=1$ is highlighted with a dot. Both axes are $\log_{10}$-scaled.\label{fig:rho_kappa}}
\end{figure}

Figure \ref{fig:medqq:med} complements Figure \ref{fig:avg-rmse:avg} by showing the estimated median curves of the cross-validation error $R(n,d)$. It illustrates the skewed nature of the distribution of $R(n,d)$ and its convergence, also for $d=1$, toward zero. Figures \ref{fig:medqq:qq} and \ref{fig:pval} complement the results in Section \ref{sec:num} regarding the convergence to normality of $R(n,d)$. In particular, Figure \ref{fig:pval} shows how normality is not rejected for $d=7,8,9,10$, even for small sample sizes $n$ and with Monte Carlo samples of size $M=10,\!000$.

\begin{figure}[htb!]
  \centering
  \includegraphics[width=0.5\linewidth]{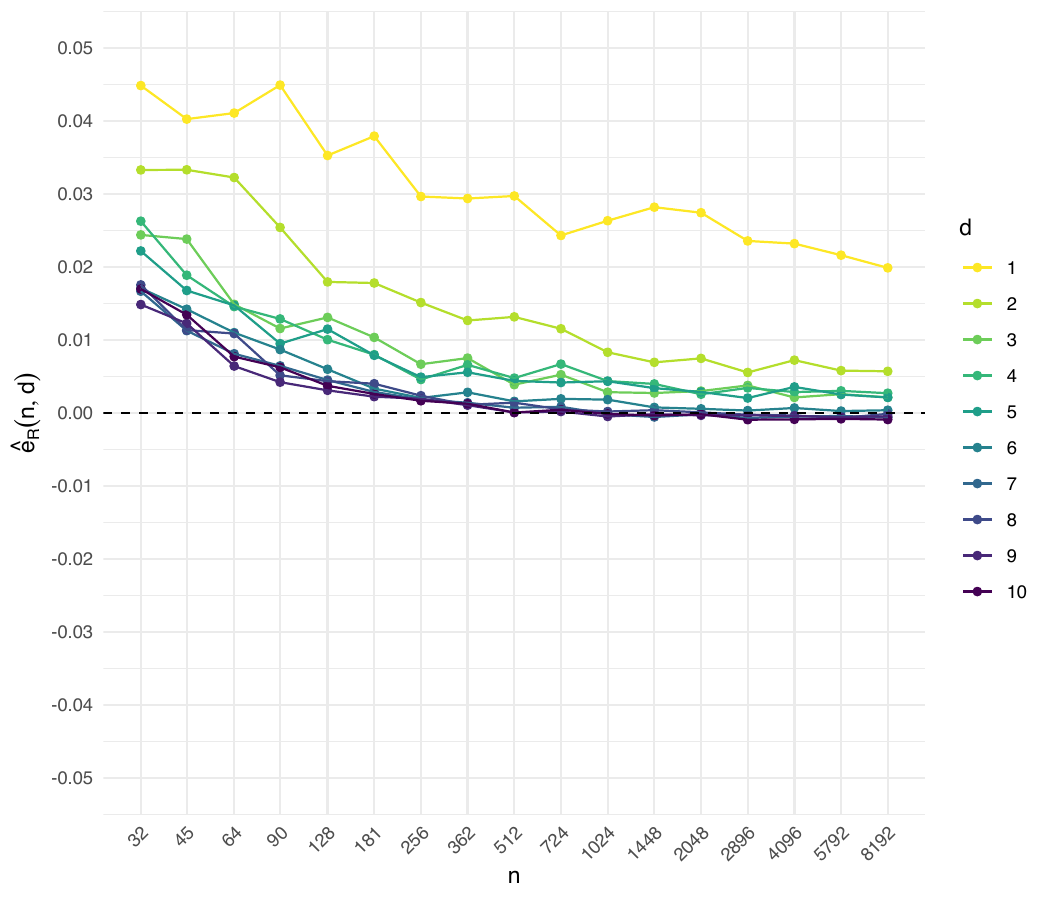}
  \caption{Median curves $n\mapsto\widehat{\mathrm{med}}(n,d)$ of $\{R^{(j)}(n,d)\}_{j=1}^M$ for dimensions $d=1,2,\ldots,10$, with $\log_2$-scale used in the horizontal axis.\label{fig:medqq:med}}
\end{figure}

\begin{figure}[htb!]
  \centering
  \includegraphics[width=\linewidth]{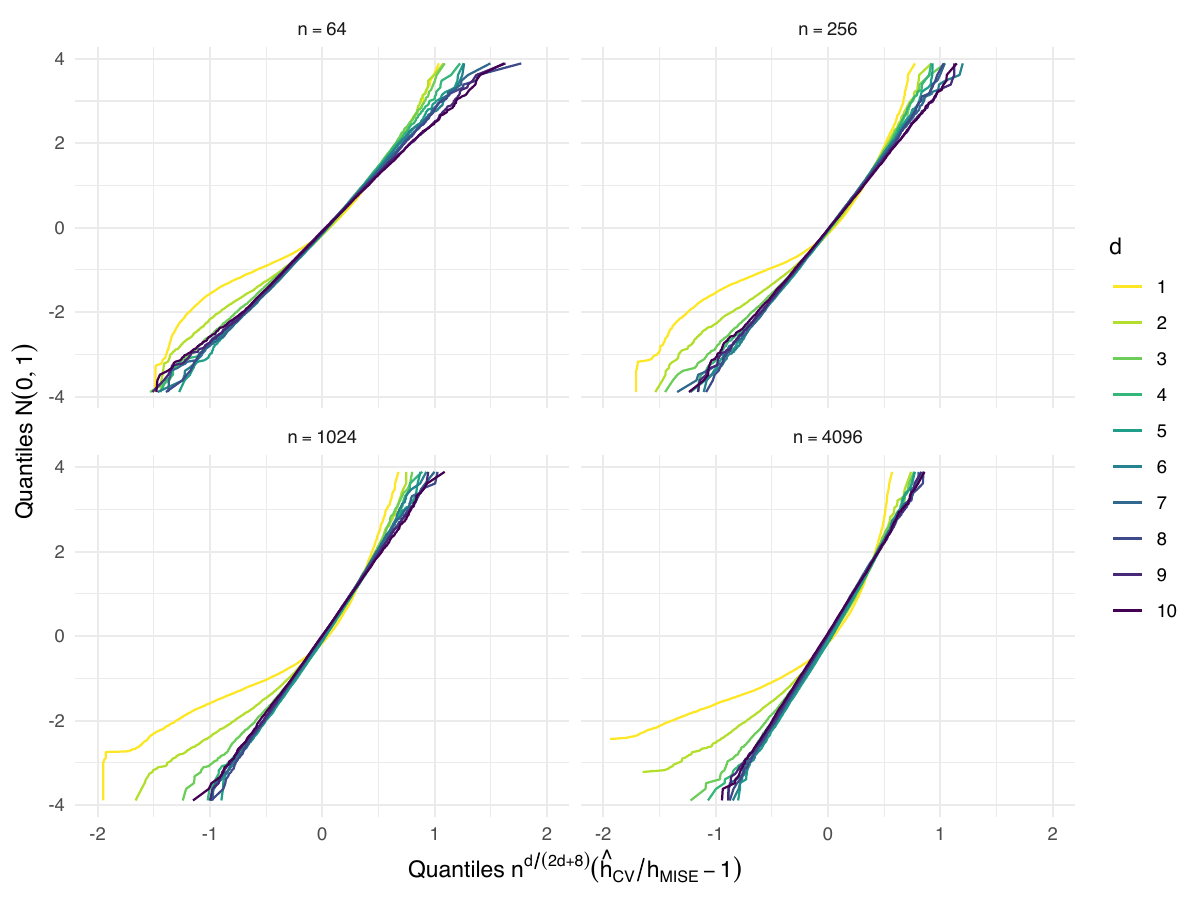}
  \caption{Normal QQ-plots with the quantiles of $\{n^{-\beta_*(d)}R^{(j)}(n,d)\}_{j=1}^M$ for $n=64,256,1024,4096$ against those of $\mathcal{N}(0,1)$. \label{fig:medqq:qq}}
\end{figure}

\begin{figure}[htb!]
  \centering
  \includegraphics[width=0.6\linewidth]{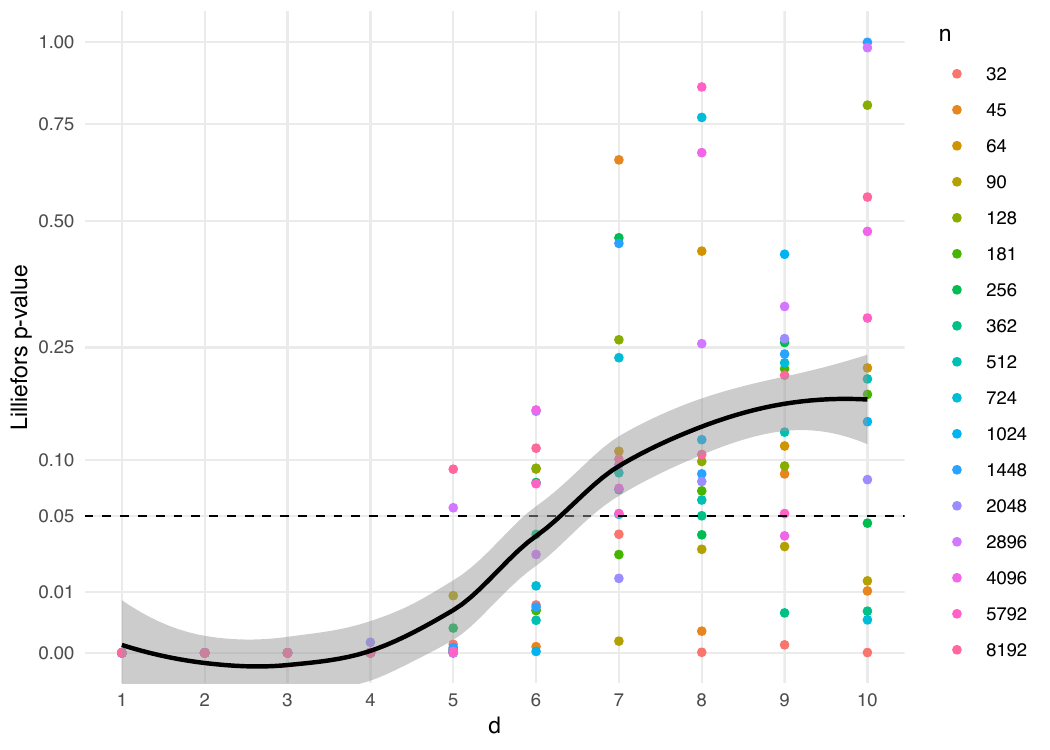}
  \caption{Testing of the normality of $n^{-\beta_*(d)}R(n,d)$ for sample sizes $n=2^\ell$, $\ell=5,6,\ldots,13$ and dimensions $d=1,2,\ldots,10$. The points represent the $p$-values of the Lilliefors test of normality applied on the sample $\{n^{-\beta_*(d)}R^{(j)}(n,d)\}_{j=1}^M$, for each pair $(n,d)$. The smooth curve is a LOESS fit to the $p$-values with a span parameter equal to $0.75$. The vertical axis is in a square root scale. \label{fig:pval}}
\end{figure}

\fi

\end{document}